\documentclass[11pt] {amsart}
\usepackage{amsmath,amsfonts,amssymb,amsxtra,latexsym,amscd,enumerate,amsthm}

\usepackage{latexsym}
\usepackage{epsfig}
\usepackage{verbatim}

\def\t{\tau}
\def\r{\mathcal{R}}
\def\g{\gamma}
\def\a{\alpha}
\def\v{\varphi}

\def\d{\delta}
\def\b{\beta}

\def\p{\Phi}
\def\f{\varphi}
\def\ep{\epsilon}
\def\es{\emptyset}
\def\l{\lambda}
\def\o{\omega}

\def\R{\mathbb{R}}
\def\r{\mathcal{R}}
\def\C{\mathcal{C}}
\def\I{\mathcal{I}}
\def\Q{\mathcal{Q}}
\def\Z{\mathbb{Z}}
\def\K{\mathcal{K}}
\def\D{\mathcal{D}}
\def\P{\mathbb{P}}
\def\p{\mathcal{P}}

\def\L{\mathcal{L}}

\def\cN{\mathcal{N}}
\def\J{\mathcal{J}}
\def\M{\mathcal{M}}

\def\S{\mathcal{S}}
\def\T{\mathcal{T}}
\def\A{\mathcal{A}}
\def\B{\mathcal{B}}
\def\TT{\mathbb{T}}
\def\N{\mathbb{N}}
\def\Z{\mathbb{Z}}
\def\V{\mathcal{V}}
\def\beq{\begin{equation}}
\def\eeq{\end{equation}}

\def\beq{\begin{equation}}
\def\eeq{\end{equation}}

\newtheorem{thm}{Theorem}

\newtheorem{claim}[thm]{Claim}

\newtheorem{obs}[thm]{Observation}

\newtheorem{d0}[thm]{Definition}
\newtheorem{o0}[thm]{Observation}

\newtheorem{l1}[thm]{Lemma}
\newtheorem{p1}[thm]{Proposition}

\begin{document}
\title[The Polynomial Carleson Operator]{The Polynomial Carleson Operator}

\author{Victor Lie}

\date{\today}

\address{Victor Lie, Department of Mathematics, Purdue, IN 46907 USA}

\thanks{The author was supported by the National Science Foundation under Grant No. DMS-1500958. The current paper was revised while the author was in residence at the Mathematical Sciences Research Institute in Berkeley, California, during the Spring 2017 semester.}

\email{vlie@purdue.edu}

\address{Institute of Mathematics of the Romanian Academy, Bucharest, RO
70700 \newline \indent  P.O. Box 1-764}

\keywords{Time-frequency analysis, Carleson's Theorem, higher order wave-packet analysis.}

\subjclass[2000]{42A20, 42A50.}

\dedicatory{Dedicated to Elias Stein on the\\ occasion of his $80^{\textrm{th}}$ birthday celebration.}

\maketitle

\begin{abstract}
We prove affirmatively the one dimensional case of a conjecture of Stein regarding the $L^p$-boundedness of the Polynomial Carleson operator, for $1<p<\infty$.

The proof is based on two new ideas: i) developing
a framework for \emph{higher-order wave-packet analysis} that is consistent with the
time-frequency analysis of the (generalized) Carleson operator, and ii) a new tile discretization
of the time-frequency plane that has the major consequence of
\emph{eliminating the exceptional sets} from the analysis of the Carleson operator.
As a further consequence, we are able to provide the full $L^p$ boundedness range and prove directly -- without interpolation techniques -- the strong $L^2$ bound for the (generalized) Carleson operator, answering
a question raised by C. Fefferman.
\end{abstract}

\section{\bf Introduction}

In this paper we will discuss the following conjecture of E. Stein regarding the behavior of the so-called Polynomial Carleson operator:

$\newline$
{\bf Conjecture (\cite{s2},\cite{sw}).}
{\it Let $G$ denote either $\TT$ or $\R$ with $G^n:=\prod_{j=1}^n G$, $n\in\N$. Further, let $\Q_{d,n}$ be the class of all real-coefficient polynomials in $n$ variables with no constant term and of degree less than or equal to $d$, $d\in\N$, and let $K$ be a suitable Calder\'on--Zygmund kernel on $G^n$. Then the Polynomial Carleson operator defined as
\beq\label{polcarl} C_{d,n}f(x):=\sup_{Q\in\Q_{d,n}}\left|
\,\int_{G^n}e^{i\,Q(y)}\,K(y)\,f(x-y)\,dy\,\right|\:
\eeq
obeys the bound
\beq\label{polcarlb}
\|C_{d,n} f\|_{L^p(G^n)} \lesssim \|f\|_{L^p(G^n)}
\eeq
for any $1<p<\infty$.}
$\newline$

The main result of our paper is:

$\newline$
{\bf Main Theorem.} {\it The above conjecture holds for $n=1$.}

\medskip
\subsection{Historical background and motivation.}
$\newline$

Before explaining the underlying motivation for Stein's conjecture, let us rewrite the expression \eqref{polcarl} for the Polynomial Carleson operator in two equivalent forms that will put matters in proper perspective.  Throughout this section, for simplicity, we will consider the case of $G = \R$.

First, notice that we can express
\beq\label{polcarlv1}
C_{d,n}f(x)=\sup_{\l}\left| T_{\l} f(x)\right|\:,
\eeq
with
\beq\label{Tlam}
T_{\l}f(x):= \int_{\R^n}e^{i\,Q_{\l}(y)}\,K(y)\,f(x-y)\,dy\:,
\eeq
where here $Q_{\l}(y)=\sum_{1\leq |\b|\leq d} \l_{\b}\, y^{\b}\in \Q_{d,n}$ is a general real-coefficient polynomial with no constant term in $n$ variables of degree at most $d$, with $\b=(\b_1,\ldots,\b_n)$ a multi-index in\footnote{Throughout the paper we use the convention $\N:=\{0,\,1,\,2,\ldots\}$.} $\N^n$ and $\l=(\l_{\b})_{\b}$ the sequence of coefficients of  $Q_{\l}$.

Now, by making the change of variable $y\,\mapsto \,x-y$, we notice that the operator $C_{d,n}$ is part of a larger class of maximal operators of the type
\beq\label{genT}
T_{*}f(x):= \sup_{\l}\left|\int_{\R^n}e^{i\,Q_{\l}(x,y)}\,K(x,y)\,f(y)\,dy\right|\:,
\eeq
where $Q_{\l},\,K\,:\R^n\times\R^n\,\rightarrow\,\mathbb C$ are such that the phase function $Q_{\l}$ is smooth and real-valued while $K$ is a suitable integral kernel that is smooth away from the main diagonal $x = y$.

Second, we note that it is possible to recast the problem of boundedness for $C_{d,n}$ without the parameter $\l$ (and thus, of course, without the corresponding supremum), at the price of losing smoothness of the phase in the $x$-parameter of $Q_{\l}$ in \eqref{genT}. Indeed, by applying the Kolmogorov--Seliverstov--Plessner linearization argument (\cite{Zyg}), one sees that the $L^p$-boundedness of $C_{d,n}$ follows from the corresponding $L^p$ bounds for an operator of the form
\beq\label{genT1}
\int_{\R^n}e^{i\,Q(x,y)}\,K(x,y)\,f(y)\,dy\:,
\eeq
where in the specific situation of $C_{d,n}$ we have $Q(x,\cdot)\in \Q_{d,n}$ a real polynomial whose coefficients are measurable functions of $x$, and $K(x,y)=K(x-y)$ with $K$ a suitable Calder\'on--Zygmund kernel on $\R^n$.

The interest in studying the Polynomial Carleson operator comes from several different directions
and with these alternative formulations of the operator in hand, we can now turn to discuss the motivations for such a study.
$\newline$

\noindent\textbf{A. Maximal singular oscillatory integrals in the Euclidean setting}.
$\newline$

The key prototypical example of a \emph{maximal} singular oscillatory integral is the so-called Carleson operator (presented as Example A.1 below). This operator arises naturally in the study of the almost-everywhere convergence of Fourier Series. This latter topic originates in the effort of 19th-century mathematics to provide a rigorous foundation for the theory of Fourier Series initiated by J. Fourier in \cite{Fou}. As a very brief historical overview, we mention the following landmark results:

Dirichlet established the convergence at all points of Fourier Series for \emph{differentiable} functions, while Du Bois Reymond subsequently showed the existence of \emph{continuous} functions whose Fourier Series diverge at a point (and in fact at any rational point). Once H. Lebesgue (\cite{Leb}) had established his theory of measure and integration---which provided the correct framework to understand the previous divergence pathologies as behavior on ``negligible'' sets---N. Luzin (\cite{Luz}) conjectured in 1913 that the Fourier Series of any $f\in L^2(\TT)$
converges to $f$ almost everywhere. In 1923, A. Kolmogorov (\cite{Kol1}) showed surprisingly that there are functions in $L^1(\TT)$ whose Fourier series diverge \emph{almost everywhere}. After decades of misbelief in light of Kolmogorov's result, L. Carleson proved in 1966 that Luzin's conjecture is in fact true (\cite{c1}), thereby setting the foundation for what is known today as time-frequency analysis.

By analogy with the approach to proving Lebesgue's differentiation theorem for $L^1$ functions via the $L^1$-weak bounds for the Hardy--Littlewood maximal function, Carleson established the almost-everywhere convergence of Fourier Series of $L^2$ functions  by providing $L^{2}$-weak bounds for the corresponding  maximal operator $\sup_{n\in\N}\,|S_n f(x)|$ derived from the sequence of partial Fourier sums $S_n f$ attached to $f$ -- which, up to admissible error terms, represents nothing else than the aforementioned Carleson operator $C := C_{1,1}$.

At this point, we can present several significant examples of operators in the literature that fit within the framework of either \eqref{genT} or \eqref{genT1}. These in turn will lead us naturally to consider the Polynomial Carleson operator:
$\newline$

\noindent\textsf{Example A.1.} \emph{Consider an operator as in \eqref{genT}, with $n=1$, $Q(x,y)= \l\,(x-y)$ and $K(x,y)=K(x-y)=\frac{1}{x-y}$.  Equivalently, in \eqref{genT1}, set $n=1$, $Q(x,y)= a(x)\cdot y$ with $a$ measurable, and $K(x-y)=\frac{1}{x-y}$.}
\medskip

In this context\footnote{We mention here that by applying a general transference principle due to Marcinkeiwicz and Zygmund one can show that $L^p$-bounds for the (generalized) Carleson operator over $\R$  or $\TT$ are equivalent.}, \eqref{genT} or \eqref{genT1} represents the Carleson operator over $\R$ whose $L^2$-weak boundedness implies and---based on Stein's maximal principle (\cite{s1})---is in fact equivalent to the affirmative answer to Luzin's conjecture.

The $L^p$ bounds, $1<p<\infty$, for the Carleson operator were established by R. Hunt in \cite{hu}.
\medskip

\noindent\textsf{Example A.2.} \emph{In \eqref{genT}, set $n\geq 1$, $Q_{\l}(x,y)= \l \cdot (x-y)$ and $K(x,y)=K(x-y)$ a Calder\'on--Zygmund kernel.  Equivalently, in \eqref{genT1}, set $n\geq 1$, $Q(x,y)= a(x)\cdot y$ with $a=(a_1,\ldots a_n)$ measurable, and $K$ as before.}
\medskip

This situation corresponds to the $n$-dimensional Carleson operator for which full $L^p$ bounds, $1<p<\infty$, were provided by Sj\"olin in \cite{sj2} and later reproved by different means in \cite{PT}.
\medskip

\noindent\textsf{Example A.3.} \emph{In \eqref{genT}, set $n=1$, $Q_{\l}(x,y)= \l \cdot (x-y)^2$, and $K(x,y)=K(x-y)=\frac{1}{x-y}$, with the obvious analog in \eqref{genT1}: $n=1$, $Q(x,y)= a(x)\cdot (x-y)^2$, and $K$ as before.}
\medskip

This case was proposed and treated by E. Stein (\cite{s2}). Unlike Carleson's theorem in \cite{c1}, whose proof relies on wave-packet analysis, this result was obtained via more standard Fourier analysis techniques, namely obtaining a good asymptotic formula for the Fourier transform of the expression $e^{i\l y^2}/y$ and appealing to $T T^{*}$-methods.
\medskip

\noindent\textsf{Example A.4.} \emph{In \eqref{genT}, set $n\geq 1$, $Q_{\l}(x,y)=\sum_{2\leq |\b|\leq d} \l_{\b}\, (x-y)^{\b}\in \Q_{d,n}$ with $d\geq 2$, and $K_{\l}(x,y)=K(x-y)$ with $K$ a standard Calder\'on--Zygmund kernel, again with the obvious analog in \eqref{genT1}.}
\medskip

This situation extends the previous setting from A.3 and was investigated by Stein and Wainger in \cite{sw}. Notice that this latter setting does not include Carleson or Sjolin's results, since \emph{no} linear term is allowed in $Q_{\l}$. The Stein--Wainger proof is based on Van der Corput estimates and again $T T^{*}$-methods.
\medskip

\noindent\textsf{Convergent point of interests:} \emph{A very natural motivating theme arises: to find a common path connecting the methods of proof and the results presented in Examples A.3 and A.4 (i.e., Stein (\cite{s2}) and  Stein-Wainger (\cite{sw})) with those of Examples A.1 and A.2 (i.e., Carleson--Hunt (\cite{c1}, \cite{hu}) and Sj\"olin (\cite{sj2})). We thus arrive naturally at the definition of the Polynomial Carleson operator in \eqref{polcarl} and the formulation of Stein's conjecture regarding its $L^p$ bounds.}\medskip

\noindent\textbf{B. Singular oscillatory integrals on nilpotent groups}. \medskip

In an extensive study regarding harmonic analysis on nilpotent Lie groups, \cite{RS1}, \cite{RS2},\cite{RS3},  Ricci and Stein proved, under the assumptions that $Q$ is a real polynomial
in both variables $(x,y) \in \R^n\times \R^n$ and $K(x,y)=K(x-y)$ with $K$ a standard Calder\'on--Zygmund kernel, the operator represented by \eqref{genT1} is bounded on $L^p(\mathbb R^n)$ for $1<p<\infty$. This of course can be regarded as \emph{a model case for our conjecture above}, in the situation in which the ``stopping times" represented by the coefficients of the monomials in $y$ in the expression $Q(x,\cdot)\in \Q_{d,n}$ are themselves polynomials in $x$.

Stein and Ricci's motivation in considering this problem relies on the fact that such operators appear naturally in three distinct but interrelated contexts:
\begin{itemize}
\item singular integrals on lower-dimensional varieties in $\mathbb R^n$ (see e.g. \cite{Sjoconv}, \cite{sw1}, \cite{Str}),
\item twisted convolution on the Heisenberg group and extensions to other nilpotent groups (see e.g. \cite{GeSt}, \cite{MPR}, \cite{Mul}), and
\item Radon transforms and their application to the study of the $\bar{\partial}$-Neumann problem (see e.g. \cite{PhSt1}, \cite{PhSt2}, \cite{GrUh}, \cite{Chhilb}, \cite{CNSW}).
\end{itemize}

For more on this, we refer the interested reader to the specific examples corresponding to each of these topics and appearing in \cite{RS2}.\medskip

\noindent \textbf{C. Connections with Radon-like transforms}.\medskip

With $n\in\N$ as before, let $\g\,:\R\,\rightarrow\,\R^n$ be a suitable (smooth) curve.  We define two operators on $\mathbb R^n$ as follows:\footnote{The function $f$ here is assumed to be in $L^1_{loc}(\R^n)$.}
\begin{itemize}
\item the maximal function along $\g$ given by
 \beq\label{maxgam}
M_{\g}f(x):=\sup_{0<\ep<1}\frac{1}{2\ep}\int_{|t|<\ep} |f(x-\g(t))|\,dt\:;
\eeq
\item the Hilbert transform along $\g$ given by\footnote{Throughout this paper we will ignore the principal value symbol.}
\beq\label{hilbgam}
H_{\g}f(x):=\int_{|t|< 1} f(x-\g(t))\,\frac{dt}{t}\:.
\eeq
\end{itemize}
The theory of singular integral operators of type \eqref{hilbgam} arose naturally in the study of solutions of constant-coefficient parabolic differential operators (see the works of \cite{Jon}, \cite{fr} and \cite{Fabe}). One specific example is the $L^2$-boundedness of \eqref{hilbgam} obtained by Fabes (\cite{Fabe}) in the case $n=2$ and $\g(t)=(t,t^2)$, by applying the method of rotations to a singular integral associated with the heat equation. This was later extended by several authors, e.g. \cite{STeWA}, \cite{Alpar},\cite{Kauf}, \cite{Hala}.

The study of maximal operators of type \eqref{maxgam} was hinted at by use of the method of rotations in connection with Poisson integrals on symmetric spaces (\cite{Stepois}). The first $L^p$ results were obtained by Nagel, Rivi\`ere, and Wainger in \cite{NRW1}, \cite{NRW2},  while a more general  Euclidean-translation-invariant theory was developed by Stein and Wainger (\cite{sw1}) in the case of one-dimensional submannifolds. All of the above results relied fundamentally on 1) Fourier methods via the Plancherel formula and 2) a suitable non-degeneracy curvature condition on $\g$ via the method of stationary phase.

With these settled, the interest naturally shifted to the case of ``variable" curves\footnote{Or, more generally, submanifolds.}  $\g$. Thus, in this new context, one is concerned with operators of the form
\beq\label{maxgamvar}
M_{\g}f(x):=\sup_{0<\ep<1}\frac{1}{2\ep}\int_{|t|<\ep} |f(x-\g(x,t))|\,dt\:,
\eeq
and the associated singular integral form
\beq\label{hilbgamvar}
H_{\g}f(x):=\int_{|t|< 1} f(x-\g(x,t))\,\frac{dt}{t}\:.
\eeq
This more general situation brings many complications for which new methods needed to be developed; in particular, one finds oneself in a non-translation-invariant setting,  suggesting that one needs to go beyond Fourier-analytic tools. A first step in this direction was made by Nagel, Stein and Wainger in \cite{NSWvarcurv}, where they obtained an $L^2$ result in the special case of some smooth variable curves $\g$.\footnote{Note that here it is essential that $\g(x,t)$ be smooth not just in the $t$-parameter but also in the $x$-parameter.}  Their approach relied on $T\,T^{*}$  methods. This work was greatly extended in the deep study of \cite{CNSW}. Other extensions to more general contexts such as nilpotent Lie groups or integral operators arising from the study of boundary-value problems in connection with the $\bar{\partial}$-Neumann problem for strongly pseudo-convex domains were already discussed in the ``Singular oscillatory integrals on nilpotent groups'' subsection above. All of these results relied on various curvature and smoothness conditions.

In an effort to unify and extend many of the above themes one could aim to:

\noindent \emph{i)  require minimal or no smoothness in the $x$-parameter, or}

\noindent \emph{ii) preserve smoothness but drop the curvature condition in the $t$-parameter.}\medskip

The various possible combinations of the presence of one or both of the above items give rise to a new class of problems, which generally are significantly more involved than the problem described above and for which there is presently no satisfactory answer. To understand the relevance and difficulty of some of these classes of problems we list here several important examples; for simplicity we only focuss on the case $n=2$ and hence $x=(x_1,x_2)\in\R^2$:
\medskip

\noindent\textsf{Example C.1.} \emph{The case $\g(x,t)=(t, v(x)\,t)$.}
\medskip

This is one of the most striking examples. Let us assume first that we only have \emph{ii)} above and thus presume $v$ is sufficiently smooth. For $v$ \emph{analytic}, Bourgain (\cite{Bolip}) proved $L^p$ bounds on \eqref{maxgamvar},\footnote{Strictly speaking his result is for $p=2$, but the extension to the case $1<p\leq \infty$ is more or less standard.} while the analogous result for the Hilbert transform was proved in a slightly more general context by Stein and Street in \cite{SS}.

Assuming now that both  \emph{ii)} and \emph{i)} are present, the story is as follows: If $v(x)=v(x_1,x_2)= v(x_1)$ is a function of only \emph{one variable} that is only assumed to be \emph{measurable}, then the $L^2$ boundedness of \eqref{hilbgamvar} is \emph{equivalent} to Carleson's theorem on the pointwise convergence of Fourier Series discussed in Example A.1..  More general $L^p$ bounds---but still not within the fully expected range $1<p<\infty$---were only recently obtained (see \cite{Bat} and \cite{BT}).  Regarding the general setting of genuinely two-variable vector fields $v$, it is a well known fact that mere measurability, or even $\alpha$-H\"older continuity with any $\alpha < 1$, is not enough to guarantee any $L^p$ bounds\footnote{Excepting of course the trivial case $p=\infty$ for the operator \eqref{maxgamvar}.} for either \eqref{maxgamvar} or \eqref{hilbgamvar}. The difficult and long-standing open problem of whether or not Lipschitz regularity\footnote{With suitable smallness condition on $\|v\|_{Lip}$.} of $v$ is enough to imply any non-trivial $L^p$ bounds for \eqref{maxgamvar} is often referred to as the \emph{Zygmund conjecture}.  The analogous problem for the Hilbert transform \eqref{hilbgamvar} was raised by Stein and is currently also widely open. As of today, the best general\footnote{\textit{I.e.}, with no extra assumption that $v$ be essentially a Lipschitz perturbation of a single-variable vector field.} regularity result is due to Lacey and Li (\cite{LL1}), who via time-frequency analysis proved---using only measurability assumptions on $u$---$L^p$ control for $p>2$ over the Hilbert transform restricted to annuli. As a consequence of this last result they also showed that if one assumes that $v$ has $C^{1+\ep}$ regularity then $H_{\g}$ is bounded on $L^2(\R^2)$. For more on this, we invite the reader to consult \cite{LL1} and \cite{PGTZ} and the bibliography therein.
\medskip

\noindent\textsf{Example C.2.} \emph{The case $\g(x,t)=(t,\,v(x)\,t^2)$.}
\medskip

In this situation we completely remove item \emph{ii)}, reimposing a non-trivial curvature in $t$. If $v$ is only assumed to be measurable, then $L^p$ bounds with $2<p\leq \infty$ are known to be true for \eqref{maxgamvar} (\cite{MRi}) and to fail for \eqref{hilbgamvar} (\cite{Kar}). If $v$ Lipschitz, then $L^p$ bounds for the full range $1<p\leq \infty$ hold
for both $M_{\g}$  (\cite{GHLR}) and  $H_{\g}$  (\cite{PGTZ}). Notice again that if $v(x)=v(x_1)$ is a measurable function of only \emph{one variable}, then the $L^2$-boundeness of  \eqref{hilbgamvar} is equivalent to Stein's result (\cite{s2}) discussed in Example A.3. above.
\medskip

\noindent\textsf{Example C.3.} \emph{The case $\g(x,t)=(t,\,\sum_{\b\leq d} v_\b(x_1)\,t^\b)$ with $d\in \N$, $d\geq 2$ and $v_{\b}$ measurable functions.}
\medskip

This represents a natural attempt to unify Examples C.1 and C.2 in terms of the $t$-variable behavior, at the price of restricting the $x$-dependence of the $v_{\beta}$'s to only the first variable.  Based on our comments above, one can easily see now that the $L^2$ bounds of  \eqref{hilbgamvar} in this setting are in fact \emph{equivalent} to the Polynoimal Carleson Conjecture stated the beginning of our paper for the case $n=1$ and $p=2$.

\subsection{Further motivation.}
$\newline$

In this subsection we focus on the classes of symmetries for various relevant operators.

Fix $f\in L^2(\R)$ and define the following symmetries:
\begin{itemize}
\item Translations:
\beq\label{Tr}
\t_{y}f(x):= f(x-y)\,,\:\:\:a\in\R\,;
\eeq
\item Dilations:
\beq\label{Dil}
D_{\l}f(x):= \l^{\frac{1}{2}}\,f(\l x)\,,\:\:\:\l\in\R_{+}\,;
\eeq
\item Generalized modulations of order $j$, $j\in\N$:
\beq\label{Modj}
M_{j,a_j}f(x):=e^{i a_j x^j}\:f(x)\,\:\:\:a_j\in\R\:.
\eeq
\end{itemize}

\subsubsection{Hilbert Transform}

As is well known, the classical Hilbert transform over $\R$, defined as
\beq\label{Hilb}
Hf(x):=\int_{\R} f(x-y)\,\frac{dy}{y}\:,
\eeq
is the only $L^2$-bounded linear operator (up to linear combinations with the identity operator) that commutes with translations and dilations, that is:

\noindent 1) $H\,\t_y=\t_y H$;

\noindent 2) $H D_{\l}= D_{\l} H$.

The $L^p$-boundedness, $1<p<\infty$, of the Hilbert transform is due to M. Riesz (\cite{Rie}),  and along the years many other proofs have been found. A particularly suggestive approach studies the action of the Hilbert transform over a wavelet system, using as an intermediate step the existence of wavelet systems that form bases for $L^2(\R)$. Recall that a wavelet system may be generated by the discrete action of dilation and translation symmetries on a single function, that is $\{D_{2^j} \t_{k} \f \}_{k,j\in\Z}$ with $\f$ a suitable smooth function on $\R$.

It is precisely this symmetry of the Hilbert transform with respect to the translation and dilation actions generating wavelet bases that give symbolic value to the wavelet-based study of the Hilbert transform: \emph{an operator's symmetries are to be involved in the analysis of its boundedness properties}. In this sense, one can view the wavelet theory as a dyadic framework for Calder\'on--Zygmund theory.

\subsubsection{Carleson operator}

Over $\R$, the Carleson operator
\beq\label{carl}
C f(x)= C_{1,1}f(x):=\sup_{a\in\R}\left|
\,\int_{\R}e^{i\,a\, y}\,\frac{1}{y}\,f(x-y)\,dy\,\right|\,,
\eeq
can be rewritten as\footnote{For notational simplicity we refer \textbf{to} modulations of order one as simply ``modulations'', and instead of $M_{1,a}$ we simply write $M_a$.}
\beq\label{carlrew}
C f(x)= \sup_{a\in\R}\left| M_{a}^{*}\,H\,M_{a}\,f(x)\right|\,,
\eeq
where throughout the paper we denote the adjoint of an operator $T$ by $T^*$.

It is now easily observed that the Carleson operator is a maximal (sublinear) operator that commutes with translations and dilations and is invariant under modulations of order one.  That is, beyond commuting with translations and dilations  as in 1) and 2) above (with $H$ replaced by $C$), the Carleson operator obeys the further symmetry

\noindent 3) $C M_a = C$.

This suggests that any method one chooses to prove the $L^2$ (weak) boundedness of $C$ should remain invariant under such symmetries, in particular under the standard modulation symmetry. This heuristic principle was in fact followed in disguise by Carleson in \cite{c1}, when he performed an analysis of the time-frequency ``adapted" Fourier coefficients of the input function of $C$.

The same heuristic principle was later expressed explicitly by C. Fefferman in his influential new proof of Carleson's result, \cite{f}, where he introduced the wave-packet discretization of the Carleson operator.  Mirroring the wavelet approach in the Hilbert transform setting, Fefferman used the elementary building blocks consisting of wave-packets, that is objects of the form $\{M_{a 2^{-j}} D_{2^j} \t_{k} \f \}_{a,k,j\in\Z}$, where again $\f$ is a suitable smooth function over $\R$.

Consequently, from the above description of the work in \cite{c1} and  \cite{f}, we see that the natural framework for standard time-frequency analysis lies within wave-packet theory, which in turn relies on the action of the three relevant symmetries: dilations, translations and (standard) modulation.

\subsubsection{Bilinear Hilbert Transform}\label{BHT}

Consider the Bilinear Hilbert transform $B$, defined \textit{a priori} for Schwartz functions $f,\,g\in \S(\R)$ by
\beq\label{bht}
B(f,g)(x):=\int_{\R} f(x-t)\,f(x+t)\,\frac{dt}{t}\,.
\eeq
This bilinear operator appeared in Calder\'on's study of the Cauchy integral on Lipschitz curves (\cite{Cal}). In this context, Calder\'on conjectured that $B$ maps boundedly $L^p\times L^q\,\rightarrow L^r$ whenever $1<p,\,q<\infty$  and $\frac{1}{p}+\frac{1}{q}=\frac{1}{r}$.

One of the key insights in approaching this problem is to realize that the Bilinear Hilbert transform shares many similarities with the Carleson operator above. Indeed, in addition to the by now standard symmetries of commutation with translation and dilation, one also has the modulation symmetry given by
\beq\label{bhtmod}
B(M_a f,M_a g)= M_{2a} B(f,g)\,.
\eeq
Applying thus the heuristic principle from above, one expects wave-packet analysis to play a key role in this problem. The confirmation of this fact came in \cite{lt1} and \cite{lt2}, where Lacey and Thiele proved that Calderon's conjecture holds under the supplementary restriction $r>\frac{2}{3}$. (Despite sustained effort, the remaining case $\frac{1}{2}<r\leq \frac{2}{3}$ is still open.)   Using the time-frequency tools developed in these papers, they were able to give a third, concise, proof of Carleson's Theorem (\cite{lt3}).

\subsubsection{Trilinear Hilbert Transform}

After the Lacey--Thiele breakthrough, a series of papers (e.g. \cite{MTT1}, \cite{MTT2}, \cite{MPTT}, \cite{dtt}, \cite{dltt}) extended the modern time-frequency framework to many other classes of multilinear operators motivated by ergodic theory, nonlinear scattering theory, and other fields. However, in all these papers, the underlying common feature is that any of the treated operators are at least ``morally'' invariant under translations, dilations, and linear modulations. For this reason, these problems could be successfully addressed by the standard wave-packet theory developed for treating the Carleson operator and later the Bilinear Hilbert transform.

However, the situation changes if one investigates the boundedness of the so-called Trilinear Hilbert Transform
\beq\label{tht}
T(f,g,h)(x):=\int_{\R} f(x+t)\,g(x+2t)\,h(x+3t)\frac{dt}{t}\,.
\eeq
The main question is whether $T$ maps $L^p\times L^q\times L^r\,\rightarrow\,L^s$ boundedly with the expected H\"older homogeneity condition $\frac{1}{p}+\frac{1}{q}+\frac{1}{r}=\frac{1}{s}$ and $1<p,\,q,\,r<\infty$ with, say, $s\geq 1$. No significant progress has been made on this problem. A primary source of difficulty is that $T$ has \emph{more symmetries} than those already mentioned: in addition to the translation, dilation, and linear modulation symmetries, $T$ also obeys a \emph{quadratic} modulation symmetry; that is, for $a\in\R$:
\beq\label{thtsym}
T(M_{2,3a}f,\, M_{2, -3a}g,\, M_{2,a}h)(x):= M_{2,a}\,T(f,\,g,\, h)(x)\,.
\eeq
Thus, according to our symmetry principle, the standard wave-packet theory is not efficient in this setting since the (linear) wave-packet framework is not invariant under quadratic modulations. Indeed, all previous attempts to approach this problem with only linear wave-packet theory have failed. Thus, developing a higher-order wave-packet theory, that in particular includes quadratic wave packets seems a natural enterprise toward a better understanding of this problem.

\subsubsection{Polynomial Carleson operator}\label{SymPCO}

Recall the one-dimensional Polynomial Carleson operator (of degree $d\in\N$)
\beq\label{polcarlrec} C_{d,1}f(x):=\sup_{Q\in\Q_{d,1}}\left|
\,\int_{\R}e^{i\,Q(y)}\,\frac{1}{y}\,f(x-y)\,dy\,\right|\:.
\eeq
We immediately notice that our Polynomial Carleson operator enjoys translation, dilation, and linear modulation invariance, thus obeying all the preliminary conditions that point towards a wave-packet methodology in the treatment of this operator.
However, one further notices that if $d\geq 2$ then beyond the previous symmetries, the Polynomial Carleson operator $C_{d,1}$ is further invariant under the action of higher-order modulations (see \eqref{Modj}) given by $\{M_{j,a_j}\}_{j\in\{2,\ldots, d\}}$. Thus, based on our earlier considerations, it seems natural that a successful approach to Stein's conjecture on the $L^p$-boundedness of the Polynomial Carleson operator should involve higher-order wave-packet theory. As we will see, this is indeed the case---in our proof of the one-dimensional case of this conjecture we develop a new way of representing and understanding the time-frequency representation and interaction of higher-order wave packets.

We stress that our analysis of the Polynomial Carleson operator (including here the Quadratic Carleson operator partially treated in \cite{q}) represents the first step in the present literature in passing from the (standard) linear to the higher-order wave-packet approach. With respect to the hierarchy of \emph{symmetry complexity}, the Polynomial Carleson operator is one level up relative to the standard Carleson operator or the Bilinear Hilbert transform. Also notice that if we fix $d=2$ ($n=1$) the Quadratic Carleson operator $C_{2,1}$ obeys similar symmetry invariances with the Trilinear Hilbert transform and thus can be regarded as an intermediate milestone between our understanding of the Bilinear and the Trilinear Hilbert transform.

We end this section with a word of caution: while the \emph{symmetry complexity} paradigm serves as a helpful heuristic in understanding the level of difficulty and the nature of the approach involved in bounding certain operators, this hierarchy need not be taken \textit{ad litteram}.  Indeed, the deeper structure of a given operator may reveal several other subtleties that significantly impact the difficulty of addressing the operator's boundedness.  For example, such subtleties likely render the problem of the boundedness Trilinear Hilbert transform extremely difficult and in particular possibly more challenging than the boundedness problem solved in this paper for the Polynomial Carleson operator $C_{d,1}$. This is the case even though, for large $d\in\N$, $(d\geq 3)$, the Polynomial Carleson operator has \emph{more symmetries} than the Trilinear Hilbert transform.

\subsection{Intermediate results}
$\newline$

Having motivated the conjecture on the Polynomial Carleson operator from various perspectives in the preceding sections, we now briefly recapitulate the results to date bearing directly on special cases of the conjecture itself; the reader will note that some of these partial results have already been discussed above.

As mentioned previously, Stein's conjecture on the Polynomial Carleson operator can be regarded:

\begin{itemize}
\item in the case $n=1$ as the extension of the celebrated Carleson--Hunt (\cite{c1}, \cite{hu}) Theorem that $C_{1,1}$ is bounded from $L^p$ to $L^p$ as long as $1<p<\infty$; and
\item in the case of general $n$ as the extension of Sj\"olin's result (\cite{sj2}); see also \cite{PT}.
\end{itemize}

As described above, under the crucial limiting assumption that the supremum in \eqref{polcarl} be taken over polynomials with \emph{no linear term}, a special case of the conjecture was established in work of Stein (\cite{s2}, for dimension $n=1$ and quadratic polynomials) and Stein--Wainger (\cite{sw}, for general dimensions $n$ and polynomials of arbitrary degree $d$).  We again note that, due to the absence of linear terms in the phase of the kernel, the Stein and Stein--Wainger results do not contain the Carleson--Hunt result.\footnote{Note further that the operators considered by Stein and Wainger have no (generalized) modulation symmetry; thus, based on the symmetry complexity heuristic discussed earlier, one expects the analysis of such operators to be significantly simpler.}

Finally, in \cite{q}, we made a significant advance by proving the $L^2$-weak boundedness of the full Quadratic Carleson operator $C_{2,1}$---incorporating polynomials with linear terms---in dimension one. As expected, in our proof we developed a new approach to the time-frequency analysis of the quadratic phase, relying on the so-called \textit{relational} perspective introduced there in Section 2. In developing the framework for quadratic wave-packet analysis, we adapted our approach to the insights developed by Fefferman in his reproof of Carleson's theorem (\cite{f}).

\subsection{Insights in our proof}
$\newline$

Passing now to the mathematical aspects of the present paper, we mention here the two main ideas on which our proof is based:
\begin{itemize}
\item  Development of the \emph{proper framework for the higher-order wave-packet theory} that in our context needs to be adapted to the time-frequency analysis of the (Polynomial) Carleson operator.

\item  A \emph{new discretization} of the family of time-frequency tiles arising in the decomposition of our operator.  This discretization has as a major implication the \emph{elimination of exceptional sets} from the analysis of the Carleson operator. This latter fact has in turn two main consequences: i) it yields boundedness for the complete range of exponents for the one-dimensional case of Stein's conjecture,  and ii) it provides for the first time a direct proof---without recourse to interpolation---of the $L^2$-boundedness of the Carleson operator, thus answering an open question raised by C. Fefferman in \cite{f}.
\end{itemize}

Beyond these facts, there will be several other points in our approach (see e.g. Section \ref{prop12}) that extend the intuition and methods developed in \cite{q} for treating the particular case $d,\,p=2$.
These latter methods were further influenced by the powerful geometric
and combinatorial ideas presented in \cite{f}.

This being said, we briefly elaborate on the two main ideas mentioned earlier:

Regarding the higher-order wave-packet framework, we develop a tile decomposition of the time-frequency plane into Heisenberg well-localized  ``curved regions'' representing area-one neighborhoods of polynomials in the class $\Q_{d,1}$. The precise geometry of the tiles appears as a manifestation of the so-called {\it relational} perspective introduced in \cite{q} and is directly related with a good control over the inner product---see e.g. equation \eqref{IBT} below---of the ``smaller pieces'' (operators) into which $\Q_{d,1}$ is decomposed. Indeed, as the name suggests, this perspective stresses the importance of \emph{interactions} between objects rather than simply treating them independently (for further details, see Section 2 in \cite{q}). At this point it is worth mentioning that our time-frequency representation of the tiles recovers, from a completely different angle, the more general uncertainty principle discussed in relation with differential operators by C. Fefferman in \cite{Fefunc}.

With respect to tile discretization,  we design a new procedure of partitioning the family of tiles that relies on a refined definition of the concept of \emph{mass} of a tile, recursive stopping-time arguments, and a very delicate combinatorial procedure. Within this process a special role is played by the counting functions associated with suitable geometric configurations of tiles called ``trees.'' All previously known estimates involved the $L^{\infty}$ size of these counting functions, which in turn required one to excise the sets on which the
$L^{\infty}$ norms are too large. In particular, these ``exceptional" sets
caused a series of technical difficulties in all the earlier works regarding the $L^p$-boundedness
of the Carleson operator; these difficulties account for the lack of a direct approach to providing strong $L^2$ bounds. In the present paper one of the key insights is that we relate, via the mass parameter, the structure of the trees of tiles with the behavior of the counting functions, thereby enabling us to replace the previous $L^{\infty}$-norm estimates with weaker $BMO$-norm-type estimates.

\subsection{Structure of the paper}
$\newline$

Next, we briefly outline the structure of our paper:

\begin{itemize}
\item In Section \ref{Not} we establish various notation and present
the general procedure of
constructing our tiles.
\item In Section \ref{Discret} we elaborate on the discretization of our operator $C_{d,1}$.
\item Section \ref{Tileinteract} is dedicated
to the study of the interaction between tiles.
\item The key idea in organizing
the family of tiles and the Main Proposition are presented in
Section \ref{proofmainth}.
\item Next, in Section \ref{Reductmainprop}, we present the main definitions and reduce the Main Proposition to two auxiliary propositions, Proposition \ref{prop1} and Proposition \ref{prop2}.
\item Section \ref{prop12}---the most technical one---contains the proofs of Propositions \ref{prop1} and \ref{prop2},
while Section \ref{remark} is dedicated to some final remarks.
\item In the Appendix we include several useful results regarding the
distribution and growth of polynomials.
\end{itemize}

Finally, given that in many respects \cite{f} and \cite{q} can be regarded as a foundation
for this paper, when possible, we have chosen to preserve here the notation, definitions and
general structure of those earlier works.
\medskip

{\bf Acknowledgements.} I would like to express my deep gratitude to Charlie Fefferman for reading parts of the manuscript and providing useful feedback. Also, I would like to thank Jim Wright for helpful advice about the history of the problem and to Christoph Thiele for first mentioning me about the question regarding the behavior of the Quadratic (Polynomial) Carleson operator. Finally, I thank Pavel Zorin-Kranich for pointing out some minor typos in an earlier version of this paper.

\section{\bf Notations and construction of the tiles}\label{Not}

We start by introducing the corresponding canonical dyadic grids on\footnote{Depending on our convenience the symbol $\TT$  stands for either $[-\frac{1}{2},\,\frac{1}{2})$ - when appearing in the definition of the Polynomial Carleson operator, or $[0,1)$ - when referring to the discretization of our time-frequency plane.} $[0,1)=\TT$ and in $\R$.
Throughout the paper the letters $I$ and $J$ refer to dyadic
intervals corresponding to the grid associated with $\TT$ while the greek letters $\a^1,\ldots ,
\a^d$, with $d\in\N$ a fixed parameter, stand for dyadic intervals associated with the grid in
$\R$. All the dyadic intervals considered in this paper are of the form $[k2^{-j},\,(k+1)2^{-j})$ for appropriate $k,\,j\in\Z$.

A \emph{tile} $P$ is a $(d+1)$-tuple of dyadic
intervals, {\it i.e.}
\beq\label{tiledef}
P=[\a^1,\a^2,\ldots , \a^d,I],\:\:\textrm{s.t.}\:\:
|\a^j|=|I|^{-1},\:\:j\in\{1,\ldots , d\}\,.
\eeq
For notational simplicity we will often refer to $P=[\a^1,\a^2,\ldots , \a^d,I]$ as $P=[\vec{\a},I]$ where here $\vec{\a}=(\a^1,\a^2,\ldots , \a^d)$.

The collection of all tiles $P$ will be denoted by $\P$.

Now, for each tile $P=[\vec{\a},I]$ we will associate a
\textit{geometric time-frequency representation}, denoted with $\hat{P}$. The
exact procedure is described in several steps:
\begin{itemize}
\item for $I$ above, we set
$x_{I}=(x_I^1, x_I^2,\ldots , x^d_{I})\in \TT^d$ to be the $d-$tuple
defined inductively as follows: $x_I^1, x_I^2$ are the endpoints of the interval $I$,
then, if $d\geq 3$, $x_I^3=\frac{x_I^1+ x_I^2}{2}$ is the mid-point of $I$, next, if $d\geq 4$, $x_I^4=\frac{x_I^1+ x_I^3}{2}$ is the mid-point of the left half of $I$, next, if $d\geq 5$,
$x_I^5=\frac{x_I^3+ x_I^2}{2}$ is the mid-point of the right half of $I$ and so on until we reach the $d$-th
coordinate.
\item recalling that $\Q_{d}$ stands for the class of all real polynomials of degree at most $d$, we make the following conventions: If not
specified, $q$ will always designate an element of $\Q_{d-1}$, while
$Q$ will refer to an element of $\Q_{d}$. When
appearing together in a proof $q$ will designate the derivative of $Q$.
\item we define
$$\Q_{d-1}(P):=\{q\in \Q_{d-1}\:|\:q(x^j_{I})\in \a^j\:\:\forall\:j\in\{1,\ldots ,
d\} \}\:,$$
and set the notation
\beq\label{convbel}
 q\in P\:\:\:\textrm{iff}\:\:\:q\in \Q_{d-1}(P)\,.
\eeq
\item with all these done, we define
\beq\label{gtile} \hat{P}:=\{(x,q(x))\:|\:x\in I\:
\&\: q\in P \}\:. \eeq
\end{itemize}

 The collection of all geometric tiles $\hat{P}$ will be
denoted with $\hat{\P}$.

For each tile $P=[\vec{\a},I]=[\a^1,\a^2,\ldots , \a^d,I]\in\P$ we associate the
``\textit{central polynomial}" $q_P\in\Q_{d-1}$ given by the Lagrange
interpolation polynomial:
\beq\label{Lagrange}
 q_P(y):=\sum_{j=1}^{d}\frac{\prod_{k=1\atop{k\not=j}}^{d}(y-x_I^k)}{\prod_{k=1\atop{k\not=j}}^{d}(x_I^j-x_I^k)}\:c(\a^j)\:.
\eeq

Now, if $I$ is any (dyadic) interval we denote by $c(I)$ the
center of $I$. Let $I_r$ be the ``right brother" of I, that is, the interval having the properties:
$c(I_r)=c(I)+|I|$ and $|I_r|=|I|$; similarly, the ``left brother" of
$I$ will be denoted $I_l$ with $c(I_l)=c(I)-|I|$ and $|I_l|=|I|$. If
$a>0$ is some real number, by $aI$ we mean the interval with the
same center $c(I)$ and with length $|aI|=a|I|$; the same conventions
apply to intervals $\{\a^k\}_k$.

In the following we will also work with dilates of our tiles: for
$a>0$ and $P=[\a^1,\a^2,\ldots , \a^d,I]$ we set
$aP:=[a\a^1,a\a^2,\ldots , a\a^d,I]$. Similarly, we write
$$a\hat{P}:=\widehat{aP}=\{(x,q(x))\:|\:x\in I\: \&\: q\in\Q_{d-1}(aP) \}\:.$$
Also, if $\p\subseteq \P$ then by convention
$a\p:=\left\{aP\:|\:P\in\p\right\}$; similarly, if
$\hat{\p}\subseteq \hat{\P}$ then
$\widehat{a\p}:=\left\{\widehat{aP}\:|\:P\in\p\right\}$.

For $P=[\vec{\a},I]$ we denote the collection of its
neighbors by
$$N(P)=\{P'=[\vec{\a'},I]\:|\:{\a'^k}\in\{\a^k,\:\a^k_r,\:\a^k_l\}\:\:\:\forall\:k\in\{1,\ldots ,
d\}\}\,.$$

Assume $P=[\vec{\a},I_P]$. We define
\beq\label{int1}
I_{P^{*}}=[c(I_P)+\frac{17}{2}|I_P|,\,c(I_P)+\frac{3}{2}|I_P|)\cup [c(I_P)-\frac{17}{2}|I_P|,\,c(I_P)-\frac{3}{2}|I_P|)
\eeq
and let
\beq\label{partint1}
I_{P^{*}}=\bigcup_{j=1}^{14}I_{P^{*}}^j
\eeq
be the partition of $I_{P^{*}}$ into dyadic intervals of length $|I_P|$.

Also we let
\beq\label{int2}
I_{\tilde{P}}=\tilde{I}_{P}:=17I_P\,.
\eeq

In some situations, for notational simplicity, we will abuse the notation and identify\footnote{There should be no confusion as the precise meaning should be clear from the context.} $P=[\vec{\a},I_P]\in\P$ with its correspondent representation $\widehat{P}\in \widehat{\P}$. Similarly, we will often identify $P^{*}$ with its geometric representation $\widehat{P}^{*}$ where
\beq\label{pstar}
\widehat{P}^{*}:=\{(x,q(x))\:|\:x\in I_{P^{*}}\: \&\: q\in P\}\,.
\eeq

Throughout the paper $p$ will be the index of the Lebesgue
space $L^p$ and, unless otherwise mentioned, will obey $1<p<\infty$. Also,
$p'$ will be its H\"older conjugate ({\it i.e.}
$\frac{1}{p}+\frac{1}{p'}=1$), while $p^*=^{def}\min(p,p')$.

For $f\in L^p(\TT)$, we denote by
$$Mf(x)=\sup_{x\in I}\frac{1}{|I|}\int_{I}|f|$$ the Hardy-Littlewood
maximal function associated to $f$.

For $A,\:B>0 $ we say $A\lesssim B$ (resp. $A\gtrsim B$) if there exists an
absolute constant $C>0$ such that $A<CB$ (resp. $A>CB$); if the constant $C$
depends on some quantity $\d>0$ then we may write $A\lesssim_{\d}B$.
If $C^{-1}A<B<CA$ for some (positive) absolute constant $C$
then we write $A\approx B$. Also we write $A>>_{d} B$ iff there exists $c(d)>(100d)^{100d}$ such that
$A> c(d) B$.

As in \cite{q}, for $x\in \R$ we set $\left\lceil
x\right\rceil:=\frac{1}{1+|x|}$.

Through out the paper, the parameters $\eta=\eta(d)$, $c(d)$ designate positive numbers depending on $d$ while $c$ stands for a large positive number; all these are allowed to change from line to line.

\section{\bf Discretization}\label{Discret}

Our goal in this section is to present the discretization of our Polynomial Carleson operator on the torus\footnote{For both historical lineage continuity (see \cite{c1} and \cite{f}) as well as argumentation clarity we present our proof on the tours as opposed to the real line. However the latter situation follows similarly with no major modifications.} defined by
\beq\label{polcarlrecc} C_{d,1}f(x):=\sup_{Q\in\Q_{d,1}}\left|
\,\int_{\TT}e^{i\,Q(y)}\,\cot \pi y\,f(x-y)\,dy\,\right|\:.
\eeq
In what follows, for notational simplicity, we will refer to the operator $C_{d,1}$ as $T$.

In direct connection with the discussion in Section \ref{SymPCO}, we start by emphasizing the groups symmetries of $T$ as displayed in the following relation:
\beq \label{symgr}
Tf(x)=\sup_{a_1,\ldots , a_d \in \R}|M_{1,a_1}\ldots , M_{d,a_d} H{M^{*}_{1,a_1}}\ldots , {M^{*}_{d,a_d}}f(x)|=\sup_{Q\in\Q_d}|T_{Q}f(x)|\:,
\eeq
where $\{M_{j,a_j}\}_{j\in \{1,\ldots , d\}}$ is the family of
generalized modulations, defined in \eqref{Modj}, $H$ is the periodic Hilbert transform, and
\beq \label{TQ}
T_{Q}f(x)=\int_{\TT}{\cot \pi y\:e^{i\,(Q(x)-Q(x-y))}\,f(x-y)\,dy}\,,
\eeq
 with
 $Q\in\Q_d$ given by $Q(y)=\sum_{j=1}^{d}a_j\: y^j$.

Now, up to easily controlled (smooth) error terms, \eqref{TQ} can be written in an equivalent form as\footnote{Throughout the paper, we ignore possible absolute constants multiplying the kernel of our operators.}
\beq \label{tq}
T_{Q}f(x)=\int_{\TT}{\frac{1}{x-y}\,e^{i\,(\int_{y}^{x}q)}\:f(y)\,dy}\:,
\eeq
where, as mentioned in the previous section, $q$ stands for the derivative of $Q$.

Now linearizing the supremum in $T$, we write
\beq \label{tq1}
Tf(x)=T_{Q_x}f(x)=\int_{\TT}{\frac{1}{x-y}\,e^{i\,(\int_{y}^{x}q_x)}\:f(y)\,dy}\:,
\eeq
where $Q_x(y):=\sum_{j=1}^{d}a_j(x)\: y^j$ with
$\{a_j(\cdot)\}_{j\in \{1,\dots d\}}$ measurable functions and $q_x$ is the derivative
of $Q_x$, that is $q_x(t)=\frac{d}{dt}Q_x(t)$ with $\int_{y}^{x}q_x=\int_{y}^{x}q_x(t)\,dt$.

Further, proceeding as in \cite{f} and \cite{q}, we define $\psi$ to
be an odd $C^{\infty}$ function such that
\beq \label{suppp}
\operatorname{supp}\:\psi\subseteq \left\{y\in
\R\:|\:2<|y|<8\right\}\,,
\eeq
 and
$$\frac{1}{y}=\sum_{k\geq 0} \psi_k(y)\:\:\:\:\:\:\:\:\:\forall\:\:0<|y|<1\:,$$
where by definition $\psi_k(y):=2^{k}\psi(2^{k}y)$ with $k\in \N$.

Using this, we deduce that
\beq \label{tdec}
Tf(x)=\sum_{k\geq 0}T_{k}f(x):=\sum_{k\geq 0}\int_{\TT}e^{i\,(\int_{y}^{x}q_x)}\,\psi_{k}(x-y)\,f(y)\,dy\:.
\eeq

Now for each $P=[\vec{\a},I]\in\P$ let
\beq \label{defEp}
E(P):=\left\{x\in I\:|\:q_x\in P\right\}\,.
\eeq
Also, if $|I|=2^{-k}$ ($k\geq0$), we
define the operators $ T_P$ on $L^2(\TT)$ by
\beq \label{deftp}
T_{P}f(x)=\left\{\int_{\TT}e^{i\,(\int_{y}^{x}q_x)}\,\psi_{k}(x-y)\,f(y)\,dy\right\}\chi_{E(P)}(x)\:.
\eeq

As expected, if $\P_k:=\left\{P=[\vec{\a},I]\in\P\:|\:|I|=2^{-k}\right\}$, for fixed $k\in\N$, the set represented by
$$\{E(P)\}_{P\in\P_k}$$ forms a partition of $[0,1)$, and so
$$T_{k}f(x)=\sum_{P\in\P_{k}}T_{P}f(x)\:.$$

Consequently, we have
\beq \label{decT}
Tf(x)=\sum_{k\geq 0}T_{k}f(x)=\sum_{P\in\P}T_{P}f(x)\:.
\eeq

This ends our decomposition.

\begin{obs}\label{redu}
We record here two facts that will be very useful in our later reasonings:
\begin{itemize}
\item For a tile $P=[\vec{\a},I_{P}]$, based on \eqref{suppp} and \eqref{deftp}, we deduce that
\beq \label{tptp*}
\operatorname{supp}\:T_P\subseteq I_P\:\:\:\:\:\textrm{and}\:\:\:\:\:\:\operatorname{supp}\:T_P^{*}\subseteq I_{P^*}\:,
\eeq
where here $T_P^{*}$ denotes the adjoint of $T_P$.

\item Taking $D$ to be the
smallest integer larger than $100d\log_{2}(100d)$ and splitting
$$\P=\bigcup_{j=0}^{D-1}\bigcup_{k\geq0}\P_{kD+j}\,,$$
we can assume from now on that the following \textit{scale separation condition} holds:
\beq \label{separat}
\eeq
if $P_j=[\vec{\a}_j,I_j]\in\P\:$ with
$j\in\left\{1,2\right\}$ such that
$|I_1|\not=|I_2|$ then either $|I_1|\leq 2^{-D}\:|I_2|$ or $|I_2|\leq 2^{-D}\:|I_1|$.
\end{itemize}
\end{obs}

\section{\bf Quantifying the interactions between tiles}\label{Tileinteract}

In this section we will focus on the behavior of the expression
 \beq\label{IBT}
\left|\left\langle
T^{*}_{P_1}\:f,T^{*}_{P_2}\:g\right\rangle\right|\:. \eeq
Our purpose will be to show that the operator discretization in Section \ref{Discret} that has at its core the defining set $E(P)$ and is fundamentally based on the \emph{relational perspective} introduced in \cite{q} is designed such that the interaction in \eqref{IBT} is controlled by the appropriately defined normalized distance between the geometric representation of our tiles $\hat{P}_1$ and  $\hat{P}_2$ (see Lemma \ref{interact} below).

In order to realize this, we will first need to introduce some quantitative concepts that are
adapted to the information offered by the localization of
$\{T_{P_j}\}_{j\in\{1,2\}}$.

\subsection{Properties of $T_{P}$ and $T_{P}^{*}$}
$\newline$

In this section we very briefly record the time-frequency localization properties of our elementary building blocks that should be regarded as a weighted generalized wave-packet decomposition of our operator $T$.

 For $P=[\vec{\a},I]\in\P$ with $|I|=2^{-k},\:k\in
\N$, we have \beq \label{v9}
\begin{array}{rl}
        &T_{P}f(x)=\left\{\int_{\TT}\:e^{i\,(\int_{y}^{x}q_x)}\,\psi_{k}(x-y)\,f(y)\,dy\right\}\chi_{E(P)}(x)\:,  \\
    &T_{P}^{*}f(x)=\int_{\TT}\:e^{-i\,(\int_{x}^{y}q_y)}\,\psi_{k}(y-x)\,\left(\chi_{E(P)}f\right)(y)\,dy\:.
\end{array}
\eeq
As it will be better revealed as a consequence\footnote{This is the essence of relational perspective introduced in \cite{q}, namely to understand the time-frequency localization of an object depending on how it interacts - in terms of the scalar product - with similar nature exterior objects.} of Lemma \ref{interact} below, we have the following principle
\beq\label{loc}
\eeq
\begin{itemize}
\item the time-frequency
localization of $T_{P}$ is ``morally" given by the geometric representation $\hat{P}$;
\item the time-frequency localization of $T_{P}^{*}$ is
``morally" given by the geometric representation $\widehat{P^{*}}$.
\end{itemize}

\begin{o0}\label{Pshape}
Remark that, due to Lemma C in the Appendix, one may think of
$\hat{P}$ as roughly being the $|I|^{-1}$ neighborhood of the graph of the
``central polynomial" $q_P$ restricted to the spacial interval $I$.
\end{o0}

\subsection{Factors of a tile} $\newline$

In this section we introduce two important concepts that will impact our understanding of the interaction in \eqref{IBT}.

 For a tile $P=[\vec{\a},I]$ we define two quantities:

$\newline$a)$\:\:\:\:\:$     an {\it absolute} one (which may be
regarded as a self-interaction); we define the {\bf density
(analytic) factor of $P$} to be the expression \beq\label{m}
 A_{0}(P):=\frac{|E(P)|}{|I|}\:.
 \eeq

Notice that $A_{0}(P)$ determines the $L^{2}$ operator norm of $T_{P}$.

$\newline$ b)$\:\:\:\:\:\:$ a {\it relative} one (interaction of $P$ or
$\hat{P}$ with an exterior object) which is of geometric
nature. Here is the description of the concept:
$\newline$

Suppose first that we are given $q\in \Q_{d-1}$ and $J$ an
interval (not necessarily dyadic); we introduce the quantity
\beq\label{dqj}
\Delta_q(J):=\frac{\operatorname{dist}^{J}(q,0)}{|J|^{-1}}\:,
\eeq
where, for $q_{1},q_{2}\in\Q_{d-1}$,  we use the notations
$$\operatorname{dist}^{A}(q_{1},q_{2})=\sup_{y\in A}\left\{\operatorname{dist}_{y}(q_{1},q_{2})\right\}\:\:\:\:\&\:\:\:\:\operatorname{dist}_{y}(q_{1},q_{2})=\left|q_{1}(y)-q_{2}(y)\right|\:.$$

Observe that we have the \emph{monotonicity property}:

\beq\label{mon}
J_1\subseteq J_2\:\:\:\textrm{implies}\:\:\:\Delta_q(J_1)\leq \Delta_q(J_2)\:.
\eeq

 Now we define the {\bf geometric factor of $P$ with
respect to $q$} as
\beq\label{geomfact}
\left\lceil \Delta_q(P)\right\rceil\:,
\eeq
where\footnote{Recall that given $x\in\R$ we let $ \left\lceil x\right\rceil:=\frac{1}{1+ |x|}$.}
\beq\label{v1}
\Delta_q(P):=\inf_{q_1 \in P}\Delta_{q-q_1}(I_P)\:.
\eeq

\subsection{Spacial Calderon-Zygmund decompositions adapted to a polynomial}
$\newline$

In this section we want to develop a general algorithm for partitioning a given interval $J\subset\TT$ into a union of dyadic intervals having suitable, ``good" properties relative to a given polynomial $q\in \Q_{d-1}$. This decomposition will be very useful later when studying the interaction displayed in \eqref{IBT}. Our precise statement and description of the algorithm is given below:

\begin{l1}\label{CZdec}[\textsf{$q-$``good" decomposition of an interval $J$]}

 Let $J\subset\TT$ be an interval such that it can be decomposed into a finite union of dyadic intervals $\bigcup_{m} J^m$  with each $|J^m|\geq \frac{|J|}{100}$.
 Also let $q\in \Q_{d-1}$ with $d\in \N,\:d\geq 2$ be a polynomial  such that
\beq\label{qbout}
q\notin \Q_{0}\:\:\:\textrm{and}\:\:\:\Delta_q(J)>0\:.
\eeq
and
\beq\label{chlam}
0<\l\leq \Delta_q(J)\:.
\eeq
Then, there exists a partition
\beq\label{part}
J=J_{s}(q,\l)\cup J_{l}(q,\l)\:,
\eeq
and $c_1(d),\,c_2(d)>0$ such that:
\begin{itemize}
\item the $(q,\,\l)-$small component $J_{s}(q,\l)$ can be written as a union of at most $9 d$ dyadic intervals having the same length\footnote{Throughout this section, our choice of $c_1(d)$ and $c_2(d)$ will be made such that the quantities $w(J,q,\l)$ and $\eta(J,q,\l)$ represent dyadic numbers.}
\beq\label{partlengt}
w(J,q,\l):= c_1(d)\,\l^{\frac{1}{d}}\,\Delta_q(J)^{-\frac{1}{d}}\,|J|\:.
\eeq
\item defining
\beq\label{hight}
\eta(J,q,\l):= c_2(d)\,\l^{\frac{d-1}{d}}\,\Delta_q(J)^{\frac{1}{d}}\,|J|^{-1}\:,
\eeq
one has
\beq\label{low}
\{x\in J\,|\,|q(x)|<\eta(J,q,\l)\}\subseteq J_{s}(q,\l)\:.
\eeq
\item the $(q,\,\l)-$large component $J_{l}(q,\l)$ can be itself partitioned into finitely many dyadic intervals
\beq\label{parthigh}
J_{l}(q,\l)=\bigcup_{W\in CZ_{(q,\l)}(J)} W\,,
\eeq
where here we define $CZ_{(q,\l)}(J)$ as the \textbf{$(q,\,\l)-$Calderon-Zygmund decomposition of $J$}, that is, the Calderon-Zygmund interval decomposition of the interval $J\setminus J_{s}(q,\l)$ relative to the set $J_{s}(q,\l)$.
\item for each $W\in CZ_{(q,\l)}(J)$ the following key properties hold\footnote{Recall that throughout the paper the constant $c(d)>0$ is allowed to change from line to line.}:
\beq\label{infsup}
\inf_{x\in W}|q(x)|\gtrsim_{d}\,\sup_{x\in W}|q(x)|\gtrsim_{d} \eta(J,q,\l)\:,
\eeq

\beq\label{length}
|W|\geq c(d)\,w(J,q,\l)\:,
\eeq

\beq\label{cont}
\Delta_q(W)\geq c(d)\,\l\:,
\eeq
and
\beq\label{derivs}
\left\|\frac{q^{(s)}}{q}\right\|_{L^{\infty}(W)}\leq c(d)\,\frac{1}{|W|^s}\:\:\:\:\:\:\:\:\:\forall\:s\in\{0,\ldots,\,d-1\}\:.
\eeq
\end{itemize}
\end{l1}

\begin{proof}

Let us first define
 $$\M_q(J)=\{x\in J\:|\:x \textrm{ is a local minimum for } |q|\}\:.$$

From \eqref{qbout} we have that $\M_q(J)\not=\emptyset$ and thus we can assume\footnote{In particular we assume that the graph of $|q|$ is not a straight line parallel with the real axis, as otherwise the above lemma is trivial.} wlog that $\M_q(J)=\{x^m\}_m$ finite.
Notice that $\M_q(J)$ contains at most $3\,d$ points.

With the previous notations, we define the $(q,\,\l)-$small component of $J_{s}$ as given by
\beq\label{defjs}
J_{s}(q,\l):=\bigcup_{j=1}^{l} I^j\,,
\eeq
where the dyadic covering $\{I^j\}_{j\in \{1,\ldots,l\}}$ has the following properties\footnote{If $I$ interval then $\mathring{I}$ stands for the interior of $I$.}:
\begin{itemize}
\item $|I^j|=w(J,q,\l)\:\:\:\:\forall\:j\in\{1,\ldots,l\}$;

\item $3\,{\mathring{I}^j}\cap\M_q(J)\not=\emptyset$.
\end{itemize}

Observe here that $l\leq 9\,d$. Also, from our hypothesis about $J$, for a proper choice of $c_1(d)$, we have that
\beq\label{dist}
\textrm{either}\:\:\textrm{dist}(J_{s}(q,\l),\,\partial J)=0\:\:\textrm{or}\:\:\textrm{dist}(J_{s}(q,\l),\,\partial J)\geq w(J,q,\l)\,.
\eeq

Next, setting
  $$\L_q^{\eta(J,q,\l)}(J)=\{x\in J\:|\:\:|q|(x)<\eta(J,q,\l)\}\:,$$
 we apply Lemma B (see the Appendix) with $I=J$ and $\eta=\eta(J,q,\l)$ and together with \eqref{partlengt} and  \eqref{hight} (for an appropriate choice of $c_2(d)$ in \eqref{hight}) we deduce
\beq\label{sizeset}
|\L_q^{\eta(J,q,\l)}(J)|\leq w(J,q,\l)\:,
\eeq
thus proving property \eqref{low}.

We pass now to the analysis of the $(q,\,\l)-$large component $J_{l}(q,\l)$.

Firstly, we notice that based on the observation \eqref{dist} definition \eqref{parthigh} makes sense.

Next, from definitions \eqref{defjs} and \eqref{parthigh} we notice that
given any $W\in CZ_{(q,\l)}(J)$ the following hold:
\begin{itemize}
\item there exist unique consecutive points $x^m,\,x^{m+1}\in \M_q(J)$ and $x^m<x^{m+1}$ such that
\beq\label{wp}
W=[a,b]\subseteq [x^m+w(J,q,\l),\,x^{m+1}-w(J,q,\l)]\:.
\eeq
\item the interval $[x^m,\,x^{m+1}]$ can be decomposed in two intervals $L_1=[x^m,\,y^m]$ and $L_2=[y^m,\,x^{m+1}]$ such that
\beq\label{mon}
\eeq
$$|q|\:\textrm{restricted to}\:L_1\:\:\textrm{is monotone increasing}\;,$$
$$|q|\:\textrm{restricted to}\:L_2\:\:\textrm{is monotone decreasing}\;.$$
\end{itemize}
Now, from \eqref{mon} we further deduce that
\beq\label{mon1}
\inf_{x\in W}|q(x)|=\min\{|q|(a),\,|q|(b)\}\:.
\eeq
Assume wlog that $\inf_{x\in W}|q(x)|=|q|(a)$. Then letting $R_1:=[x^m,\,a]$ and $R_2=[x^m,\,b]$ we have
\beq\label{mon2}
\eeq
\begin{itemize}
\item $\|q\|_{L^{\infty}(R_1)}=|q|(a)=\inf_{x\in W}|q(x)|$;

\item $\|q\|_{L^{\infty}(R_2)}=\|q\|_{L^{\infty}(W)}$;

\item $1\leq \frac{|R_2|}{|R_1|}\leq 5$.
\end{itemize}
Now \eqref{infsup} follows from \eqref{mon2} and an application of Lemma A in the Appendix.

Relation \eqref{length} follows directly from the definition of the Calderon-Zygmund decomposition $CZ_{(q,\l)}(J)$.

Next, \eqref{cont} follows from
$$\Delta_q(W)=\frac{\operatorname{dist}^{W}(q,0)}{|W|^{-1}}=
\frac{\operatorname{dist}^{R_2}(q,0)}{|W|^{-1}}\geq \frac{1}{5}\,\Delta_q(R_2)$$
$$\geq \frac{1}{5}\,\Delta_q([x^m,\,x^m+w(J,q,\l)]) \geq \frac{\eta(J,q,\l)}{5\,w(J,q,\l)^{-1}}=c(d)\,\l\:.$$
Finally, \eqref{derivs} is a direct consequence of \eqref{infsup} and Lagrange interpolation formula displayed in the proof of Lemma A by taking in \eqref{Lag} $J=W$.
\end{proof}

\subsection{The resulting estimates} $\newline$

We conclude this section by describing how the concepts and definitions
introduced above in Section \ref{Tileinteract} relate in controlling the interaction in \eqref{IBT}.

As expected, we need to quantify the relative position of
$P_1^{*}$ with respect to $P_2^{*}$. Of course, we will only
consider the nontrivial case $I_{P_1}^{*}\cap I_{P_2}^*
\not=\es $; also, throughout this section, for notational simplicity we simply set $I_{P_1}=I_1$, $I_{P_2}=I_2$ and we suppose wlog that $|I_1|\geq |I_2|$.

\begin{d0}\label{fact} [\textsf{Geometric factor associated to a pair of tiles}]

Given two tiles $P_1$ and $P_2$, we define the {\bf geometric factor of the pair ($P_1,P_2$)} by $$\left\lceil\Delta(P_1,P_2) \right\rceil\:,$$ where\footnote{Recall notation \eqref{int2}.}
$$\Delta(P_1,P_2)=\Delta_{1,2}:=\frac{\sup_{y\in \tilde{I_1}\cap \tilde{I_2}}\{\inf_{{q_1\in P_1}\atop{{q_2\in P_2}}}\operatorname{dist}_y(q_{1},q_{2})\}}{|\tilde{I_1}\cap \tilde{I_2}|^{-1}}\:.$$
 \end{d0}

\begin{d0}\label{fact1} [\textsf{Interaction polynomial}]

For $P_1$ and $P_2$ as above, we define the $(P_1,\,P_2)$-{\bf interaction polynomial} as
\beq\label{interpoly}
q_{1,2}:=q_{P_1}-q_{P_2}\:.
\eeq
\end{d0}

\begin{d0}\label{fact2}[\textsf{Critical intersection set}]

Let now $\ep_{0}\in (0,1)$. With the notations and conventions from Lemma \ref{CZdec} we define the ($\ep_0$-){\bf critical intersection set} $I_{1,2}$ of the pair $(P_1,P_2)$ as
\beq\label{criticset}
I_{1,2}:=J_{s}(q_{1,2},\l)\:,
\eeq
for the particular values:
\begin{itemize}
\item $J=\tilde{I_1}\cap \tilde{I_2}$;

\item $\l:=\Delta(P_1,P_2)\,{\left\lceil\Delta(P_1,P_2)\right\rceil}^{1-\ep_0}$.
\end{itemize}
\end{d0}

\begin{o0}\label{geom} With these notations, using the results in the Appendix and our assumption $|I_1|\geq |I_2|$, we have that if $\left\lceil\Delta_{1,2}\right\rceil<<_{d} 1$ then
$${\left\lceil\Delta_{1,2}\right\rceil}\approx_{d}\max\left\{{\left\lceil\Delta_{q_{P_1}}(P_2)\right\rceil},\:{\left\lceil\Delta_{q_{P_2}}(P_1)\right\rceil}\right\} \approx_{d} {\left\lceil\Delta_{q_{1,2}}(I_2)\right\rceil} \:.$$
\end{o0}

Now using Lemma \ref{CZdec} together with the principle of
(non-)stationary phase, one deduces the following:

\begin{l1}\label{interact}[\textsf{Tile interaction control}]

Let $P_1\:,\:P_2\:\in\P$. Then, with the above notations and conventions, we
have
\beq\label{v15} \left|\int
\tilde{\chi}_{I_{1,2}^c}T_{P_1}^{*}f\:\overline{T_{P_2}^{*}g}\:\right|\lesssim_{\:n,\:d,\:\ep_0}{\left\lceil
\Delta(P_1,P_2)\right\rceil}^{n}\:\frac{\int_{E(P_1)}
|f|\int_{E(P_2)}|g|}{\max\left(|I_1|,|I_2|\right)}\:\:\:\:\:\:\forall\:n\in
\N\:, \eeq
 \beq\label{v16} \int_{I_{1,2}}|
T_{P_1}^{*}f\:\overline{T_{P_2}^{*}g}|\lesssim_d{\left\lceil
\Delta(P_1,P_2)\right\rceil}^\frac{1-\ep_0}{d}\:\frac{\int_{E(P_1)}
|f|\int_{E(P_2)}|g|}{\max\left(|I_1|,|I_2|\right)}\:, \eeq
where
$\tilde{\chi}_{I_{1,2}^c}$ is a smooth variant of the corresponding
cut-off.

Moreover, we also have
\beq\label{v17}
\left\|T_{P_1}{T}_{P_2}^{*}\right\|_2^{2}\lesssim_{d}\min\left\{\frac{|I_2|}{|I_1|},
\frac{|I_1|}{|I_2|} \right\}{\left\lceil
\Delta(P_1,P_2)\right\rceil}^{\frac{2}{d}}\:A_0(P_1)\:A_0(P_2)\:.
\eeq
\end{l1}

\begin{proof}
Assume throughout the proof that $\Delta(P_1,P_2)>>_{d} 1$ as otherwise the above statements are trivial.

Next, notice that relation \eqref{v16} is straightforward based on the definition of $I_{1,2}$ in \eqref{criticset} and on the fact that
$$|T_{P_j}^{*}f|\lesssim
\frac{\int_{E(P_j)}|f|}{|I_j|}\,\chi_{I_j^*}\:\:\:\:\:\:\:\forall\:\:j\in\{1,2\}\,,$$
which is further a consequence of \eqref{v9}.

We now turn our attention towards \eqref{v15}.

Apply the algorithm described in Lemma \ref{CZdec} for the following parameters: $J=\tilde{I_1}\cap \tilde{I_2}$ and
$\l:=\Delta(P_1,P_2)^{\ep_0}$. We then obtain the collection $CZ_{(q,\l)}(J):=\{W_r\}_r$ representing the $(q,\,\l)-$Calderon-Zygmund decomposition of $J\setminus J_{s}(q,\l)$ relative to the set $J_{s}(q,\l)$.

Let $\varphi$ be a smooth cutoff of $\chi_{I_{1,2}^c}$ such that $\varphi\geq 0$ and
\beq\label{defmoothco}
\varphi|_{J\setminus I_{1,2}^c}=1\:\:\textrm{and}\:\:\varphi|_{\frac{3}{4}\,I_{1,2}\cup (\frac{5}{4} J)^c}=0\;.
\eeq

Take now any smooth partition of unity adapted to the collection $CZ_{(q,\l)}(J)$ such that this is identically zero on the set $\frac{3}{4}\,I_{1,2}\cup (\frac{5}{4} J)^c$. Thus, wlog we may assume that
\beq\label{smooths}
\f=\sum_{r}\f_r\,,
\eeq
where here
\beq\label{smoothss}
\eeq
\begin{itemize}
\item $\f_r$ adapted to $W_r$;
\item $\|\f_r\|_{C^{\infty}}\lesssim 1$;
\item $\f_r=1$ on $W_r$;
\item $\f_r=0$ on $[0,1]\setminus\frac{5}{4}\,W_r$.
\end{itemize}

\begin{obs}\label{extlema}
It is important to notice that for appropriate choices of the $d$-depending constants in Lemma \ref{CZdec} and based on the results in the Appendix we have that the properties of the central polynomial $q_{1,2}$ on each of the $W_r$ are transferable with no modifications (up to further $d$ dependent constants) for any difference polynomial of the form $q_1-q_2$ with $q_1\in P_1$ and $q_2\in P_2$.
\end{obs}

With this we have
\begin{align*}
\int \varphi\: T_{P_1}^{*}f\:\overline{T_{P_2}^{*}g} & =\int f\:
\overline{T_{P_1}(\varphi\, T_{P_2}^{*}g)}\\
& =\int\int (f\chi_{E(P_1)})(x)\:
(\bar{g}\chi_{E(P_2)})(s)\:\K(x,s)\,dx\, ds\:,
\end{align*}

where
$$\K(x,s)=\int e^{i\,
[\int_{y}^{s}q_s-\int_{y}^{x}q_x]}\;\psi_{k_1}(x-y)\:\varphi(y)\:\psi_{k_2}(y-s)\;dy\:.$$
Here we used the convention $|I_1|=2^{-k_1}$, $|I_2|=2^{-k_2}$ with
$k_2\geq k_1$ positive integers.

Let us set $\tilde{Q}(y):=\int_{y}^{s}q_s-\int_{y}^{x}q_x$, $\tilde{q}=\tilde{Q}'$ and $u(y):=
\psi_{k_1}(x-y)\:\psi_{k_2}(y-s)$. Then, writing $e^{i\,\tilde{Q}(y)}=\left(\frac{1}{i\,\tilde{q}(y)}\,\frac{d}{dy}\right)\,(e^{i\,\tilde{Q}(y)})$ and integrating by parts $n$ times in expression
\beq\label{intpart}
\K(x,s)=\int \left[\left(\frac{1}{i\,\tilde{q}(y)}\,\frac{d}{dy}\right)^{n}\,(e^{i\,\tilde{Q}(y)})\right]\;\varphi(y)\:u(y)\;dy\,,
 \eeq
we get
\beq\label{intpartestim}
|\K(x,s)|\lesssim_{n} \int \sum_{{a_1+\ldots+a_{n+1}=n}\atop{{b_1+\ldots+b_n=n}\atop{{a_j+b_j\leq n+1}\atop{a_j,\,b_j\in\N}}}} \prod_{j=1}^{n}\left|\left(\frac{d}{dy}\right)^{a_j}\left(\frac{1}{\tilde{q}(y)^{b_j}}\right)\right|\,
\left|\left(\frac{d}{dy}\right)^{a_{n+1}} (\varphi(y)\:u(y))\right|\,dy\:.
\eeq
Now for generic $a,\,b,\,c\in\N$, making use of Observation \ref{extlema}, we have
$$|(\frac{d}{dy})^{a}(\frac{1}{\tilde{q}^b(y)})|\lesssim_{a,b,d}
\sup_{{s\leq a}\atop{{n_1+\ldots+n_s=a}\atop{n_1,\,n_2,\,\ldots,\, n_s\in\N}}}{\frac{|\tilde{q}^{(n_1)}\ldots\tilde{q}^{(n_s)}|}{|\tilde{q}|^{s+b}}} \,,$$
$$|(\frac{d}{dy})^{b}\f(y)|\lesssim\sum_{W_r\in\mathcal{W}}\,\frac{1}{|W_r|^b}\,|\tilde{\f}_r(y)|\,,$$
and
$$|(\frac{d}{dy})^{c}u(y)|\lesssim \frac{1}{|I_2|^c}\,|\psi_{k_1}(x-y)\:\tilde{\psi}_{k_2}(y-s)|\,.$$
where here $\tilde{\f}_r$ and $\tilde{\psi}_{k_2}$ are functions with the same localization/smoothness properties as $\f_r$ and $\psi_{k_2}$ respectively.

Using now \eqref{infsup}-\eqref{derivs} we get
 \beq\label{kernel}
 |\K(x,s)|\lesssim_{n}\frac{1}{|I_1|}\,\frac{1}{|I_2|}\,\sum_{W_r}\frac{|W_r|}{\Delta_q(W_r)^{n}}\lesssim
 \frac{1}{|I_1|}\,{\left\lceil
\Delta(P_1,P_2)\right\rceil}^{n\,\ep_0}\;,
 \eeq
which proves \eqref{v15}.

For $\eqref{v17}$, we repeat the previous argument but now in the setting $\ep_0=0$ and $n=1$ and once we reach the first inequality in \eqref{kernel} we use the simpler estimate
$$\Delta_q(W_r)\geq \eta(J,q,1)\,|W_r|\approx_d \Delta(P_1,P_2)^\frac{1}{d}\,\frac{|W_r|}{|J|}\:.$$

This ends the proof of our lemma.
\end{proof}

\section{\bf The proof of the main theorem}\label{proofmainth}

\subsection{\bf A key ingredient - organizing the family of tiles}
$\newline$

In this section we will recursively partition the set of all tiles $\P$ into families of tiles with some special properties. More precisely, using induction, we will show that
\beq\label{pdensity}
\P=\bigcup_{n}\P_n\:,
\eeq
such that, roughly speaking, for each family $\P_n$
\begin{itemize}
\item the tiles inside have a uniform density factor;

\item a suitable defined counting function is under ``good" control.
\end{itemize}

\subsubsection{\textbf{Preparatives.}} To make our tile-partition precise, we need to introduce the following

\begin{d0}\label{mass} [\textsf{Mass of a tile adapted to a given environment}]

Let $\A$ be a (finite) union of dyadic intervals in $[0,1]$ and $\p$ be a finite family of tiles. For $P=[\vec{\a},I]\in\p$ with $I\subseteq\A$ we define the {\bf mass} of $P$ relative to the set of tiles $\p$ and the set $\A$ as being

\beq\label{v1} A_{\p,\A}(P):=\sup_{{P'=[\vec{\a}',I']\in\:\p}\atop{I\subseteq
I'\subseteq \A}}\frac{|E(P')|}{|I'|}\:\left\lceil
\Delta(10P,\:10P')\right\rceil^{N} \eeq where $N$ is a fixed large
natural number.
\end{d0}

Next, we introduce a qualitative concept that characterizes the
overlapping relation between tiles.

\begin{d0}\label{ord}[\textsf{Aiming for ``orderings"}]

   Let $P_j=[\vec{\a}_j,I_j]\in\P$ with $j\in\left\{1,2\right\}$. We say that
$\newline$ - $P_1\leq P_2$       iff       $\:\:\:I_1\subseteq I_2$
and      $\exists\:\:q\in P_2$ such that $q\in P_1\:,$ $\newline$ -
$P_1\trianglelefteq P_2$     iff     $\:\:\:I_1\subseteq I_2$  and
$\forall \:\:q\in P_2$ we have $q\in P_1\:.$

Also we say $P_1< P_2$ if $P_1\leq P_2$ and $|I_1|<|I_2|$. Similar statement for $\vartriangleleft$.
\end{d0}

\begin{obs}\label{ordmaxrel1}
Notice that $\leq$ is not an order relation while $\trianglelefteq$ it is. Also $P_1<P_2$ implies $2 P_1\vartriangleleft 2P_2$.
 \end{obs}

In the following two definitions we elaborate on the ``pseudo-ordering" $\leq$:

\begin{d0}\label{ordmaxrel}[\textsf{Maximal/minimal tiles within a given family}]
\begin{enumerate}
\item If $\p\subseteq\P$ is a family of tiles with some prescribed properties, we say that $P\in\p$ is \textbf{maximal} (relative to $\p$) iff
\beq\label{maxx}
\forall\:P'\in\p\:\:\textrm{s.t.}\:\:P\leq P'\:\:\textrm{we have}\:\:P=P'\:.
\eeq
\item Similarly, if $\p\subseteq\P$ is a family of tiles, we say that $P\in\p$ is \textbf{minimal} (relative to $\p$) iff
\beq\label{minn}
\forall\:P'\in\p\:\:\textrm{s.t.}\:\:P\geq P'\:\:\textrm{we have}\:\:P=P'\:.
\eeq
\end{enumerate}
\end{d0}

\begin{d0}\label{incomp}[\textsf{Incomparable/negligible family of tiles}]

We say that $\p\subset\P$ is an \textbf{incomparable} family of tiles iff
\beq \label{incr}
\forall\:P_1,\,P_2\in\p\: \textrm{we have}\:P_1\nleq P_2\:\textrm{and}\:P_2\nleq P_1\;.
\eeq
Also we call $\p\subset\P$ \textbf{negligible} if $\p$ can be written as a union of at most $c(d)$ families of incomparable tiles.
\end{d0}

We end this subsection with the following observation that connects the geometric statement $P_1\leq P_2$ with the analytic behavior of the polynomials belonging to $\{P_j\}_j$.

\begin{obs}\label{qdis} As a consequence of Definition \ref{ord} and Lemma C in the Appendix the following holds:

If $P_1=[\vec{\a}_1, I_{P_1}],\,P_2=[\vec{\a}_2, I_{P_2}]\in \P$ such that $P_1\leq P_2$ then there exists $c(d)\in (0,\,10\,(2d)^d]$ such that
\begin{itemize}
\item $\exists\:q_1\in P_1$ with
\beq \label{p11}
\sup_{q_2\in P_2} \|q_2-q_1\|_{L^{\infty}(\tilde{I}_{P_2})}\leq c(d)\,|I_{P_2}|^{-1}\:.
\eeq
\item $\forall\:q_1\in P_1$ we have
\beq \label{p12}
\sup_{q_2\in P_2} \|q_2-q_1\|_{L^{\infty}(\tilde{I}_{P_1})}\leq c(d)\,|I_{P_1}|^{-1}\:.
\eeq
\end{itemize}
\end{obs}

$\newline$
\subsubsection{\textbf{Partitioning $\P$ - the inductive algorithm}}\label{parttil}
$\newline$

In this section we present an inductive algorithm of partitioning\footnote{In an earlier version of our paper we presented a slightly modified partitioning algorithm. For more on this, please see the Remarks section, third item.} our set of tiles into
\beq\label{decP}
\P=\bigcup_{n\in\N}\P_n\,,
\eeq
with each $\P_n$ being a set of tiles of mass $n$ relative to certain space regions. Our algorithm will
be based on a stopping time process involving the John-Nirenberg inequality that is correlated with the level set analysis of various counting functions. This process is constructive and it is based on an ascending induction over $n$.
$\newline$

\noindent\textbf{Step 1 Construction of the family $\P_1$}
$\newline$

This construction will be done in two stages:
\begin{itemize}
\item Stage 1.1 - we define a sequence of nested sets $\{A_1^{k}\}_{k\in\N}$ such that
\begin{itemize}
\item $A_1^k$ is a finite union of maximal disjoint dyadic intervals;

\item we have a good control on the $L^{\infty}$ norm of a suitable ``counting function of order one" adapted to $A_1^k$.
\end{itemize}

\item Stage 1.2 - for each set $A_1^k$ we define a corresponding family of tiles $\p_1[A_1^{k}]$ with the following two key properties:
\begin{itemize}
\item  $\p_1[A_1^{k}]$ is a \textit{convex} family of tiles, that is, if $P_1\leq P\leq P_2$ with $P_1,\,P_2\in \p_1[A_1^{k}]$ then $P\in \p_1[A_1^{k}]$.

\item  each tile $P=[\vec{\a},I]\in \p_1[A_1^{k}]$ has the properties
\beq\label{massone}
\begin{array}{cc}
I\subseteq A_1^{k}\:\:\&\:\:I\nsubseteq A_1^{k+1}\,,\\\\
2^{-1}<A_{\P, A_1^k}(P)\leq 1\:.
\end{array}
\eeq
\end{itemize}
\end{itemize}

Now before effectively starting our construction we introduce the following

 \begin{d0}\label{dom}
 Let $\A=\bigcup \A_j$  and $\B=\bigcup_{k}\B_k$ be two sets such that both $\{\A_j\}_j$ and $\{\B_k\}_k$ are some collections of maximal (disjoint) dyadic intervals.

 We say that
\beq\label{dom1}
\A\prec\B\,,
\eeq
 iff each $\A_j$ is contained in some $\B_k$.

Moreover, given an absolute constant $c>0$, we write
\beq\label{dom2}
\A\prec_{c}\B\,,
\eeq
 iff $\A\prec\B$ and for any $\B_k$ the following holds:
\beq\label{dom3}
|\bigcup_{\A_j\subseteq\B_k}\A_j|\leq e^{-c}\,|\B_k|\,.
\eeq
\end{d0}

This being said, we are ready to initiate the following:
$\newline$

\noindent\textbf{Stage 1.1 Construction of the sets $\{A_1^k\}_{k\geq 1}$}
$\newline$

\noindent \textbf{1.1.1 Construction of the set $A_1^1$}

\begin{itemize}
\item We start by simply defining the set
\beq\label{set10}
A_1^{0}:=[0,1]\:.
\eeq
\item Let
\beq\label{max}
\p_{1}^{max}[A_1^{0}]:=\left\{P=[\vec{\a},I]\in\P\,|\,P\:\textrm{maximal}\:\:\&\:\:\frac{|E(P)|}{|I|}> 2^{-1}\right\}\:.
\eeq
Notice that $\p_{1}^{max}[A_1^{0}]$ is formed by \textit{disjoint} or incomparable tiles, that is
\beq\label{disjo}
\forall\:P\not=P'\in \p_{1}^{max}[A_1^{0}]\:\:\:\Rightarrow\:\:\:P\nleq P'\:\:\textrm{and}\:\:P'\nleq P\,.
\eeq
\item Next, we define
\beq\label{c1}
\C_{1}[A_1^{0}]:=\sum_{P\in\p_{1}^{max}[A_1^{0}]} \chi_{E(P)}\,,
\eeq
and notice, based on \eqref{disjo}, that
\beq\label{c1sm}
\C_{1}[A_1^{0}](x)\leq 1\:\:\:\forall\:x\in[0,1]\;.
\eeq
\item Collect the time intervals of the maximal tiles in $\p_{1}^{max}[A_1^{0}]$ into the set $\I_{1}^{max}[A_1^{0}]$, that is
\beq\label{i1}
\I_{1}^{max}[A_1^{0}]:=\{I\,|\,P=[\vec{\a},I]\in\p_{1}^{max}[A_1^{0}]\}\,.
\eeq
\item Define the \textit{counting function of order one adapted to $A_1^{0}$} as
\beq\label{countf}
\cN_{1}[A_1^{0}]:=\sum_{I\in\I_{1}^{max}[A_1^{0}]}\chi_{I}\,.
\eeq
Notice that $\cN_{1}[A_1^{0}]$ verifies the relation
\beq\label{defbmoc}
\|\cN_{1}[A_1^{0}]\|_{BMO_{C}}:=\sup_{J\textrm{dyadic}\atop{J\subseteq [0,1]}}\,\frac{\sum_{I\subseteq J\atop{I\in\I_{1}^{max}[A_1^{0}]}}|I|}{|J|}\leq 2\:.
\eeq
\item Setting now
\beq\label{defbmod}
\|\cN_{1}[A_1^{0}]\|_{BMO_{D}}:=\sup_{J\textrm{dyadic}\atop{J\subseteq [0,1]}}
\frac{1}{|J|}\,\int_{J}\left|\cN_1[A_1^{0}]-\frac{\int_{J} \cN_1[A_1^{0}]}{|J|}\right|\:,
\eeq
we deduce that
\beq\label{defbmodc}
\|\cN_1[A_1^{0}]\|_{BMO_{D}}\leq 2\|\cN_{1}[A_1^{0}]\|_{BMO_{C}}\,.
\eeq
\item Applying now the John-Nirenberg inequality, we have\footnote{Throughout the section the constant $c>10^{10}$ is an absolute constant that is allowed to change from line to line.}
\beq\label{JN}
\left|\left\{x\in J\,|\,\left|\cN_1[A_1^{0}](x)-\frac{\int_{J} \cN_1[A_1^{0}]}{|J|}\right|>\gamma\right\}\right|\lesssim |J|\,e^{-c\:\frac{\g}{\|\cN_1[A_1^{0}]\|_{BMO_{D}}}}\:.
\eeq
\item Using \eqref{defbmodc} in \eqref{JN}, for $\g> c\,\|\cN_1[A_1^{0}]\|_{BMO_{C}}$, we have that
\beq\label{JNP}
|\{x\in J\,|\,\sum_{I\subseteq J\atop{I\in\I_{1}^{max}[A_1^{0}]}}\chi_{I}(x)>\gamma\}|\lesssim |J|\,e^{-c}\:.
\eeq
\item Finally, notice that
\beq\label{A11}
A_1^{1}:=\{x\in [0,1]\,|\, \cN_{1}[A_1^{0}](x)> c\,\|\cN_{1}[A_1^{0}]\|_{BMO_{C}}\}\:,
\eeq
consists of a finite union of disjoint dyadic intervals.

Applying now \eqref{JNP} we deduce that
\beq\label{kA11}
A_1^{1}\prec_{c}[0,1]\:.
\eeq
\end{itemize}

$\newline$
\noindent \textbf{1.1.2 Construction of the set $A_1^{k}$ with $k\geq 1$}
$\newline$

We apply an inductive argument. Since the first step was already verified, we assume that as the byproduct of the step $k-1$ we obtained a set
$$A_1^{k-1}\,,$$
that can be represented as a finite union of disjoint dyadic intervals.

\begin{itemize}
\item As before, we start by identifying the collection of maximal tiles
\beq\label{max1}
\p_{1}^{max}[A_1^{k-1}]:=\left\{P=[\vec{\a},I]\in\P\,\big|\,\begin{array}{cc} P\:\textrm{maximal}\\ I\subseteq A_1^{k-1}\end{array}\:\:\&\:\:\frac{|E(P)|}{|I|}> 2^{-1}\right\}\:.
\eeq

\item Set the collection of time-intervals of maximal tiles as
\beq\label{i2}
\I_{1}^{max}[A_1^{k-1}]:=\{I\,|\,P=[\vec{\a},I]\in\p_{1}^{max}[A_1^{k-1}]\:\}\,.
\eeq

\item Define the \textit{counting function of order one adapted to $A_1^{k-1}$} as
\beq\label{countf1}
\cN_{1}[A_1^{k-1}]:=\sum_{I\in\I_{1}^{max}[A_1^{k-1}]}\chi_{I}\,,
\eeq
and applying the same reasonings as in the previous situation, notice that $\cN_{1}[A_1^{k-1}]\in BMO_{D}$ and moreover that $\|\cN_{1}[A_1^{k-1}]\|_{BMO_{C}}\leq 2$.

\item Applying now the the John-Nirenberg inequality, for $\g> c\,\|\cN_1[A_1^{k-1}]\|_{BMO_{C}}$, we have
\beq\label{JNP1}
|\{x\in J\,|\,\sum_{I\subseteq J\atop{I\in\I_{1}^{max}[A_1^{k-1}]}}\chi_{I}(x)>\gamma\}|\lesssim |J|\,e^{-c}\:.
\eeq

\item Conclude that the set
\beq\label{A11}
A_1^{k}:=\{x\in [0,1]\,|\, \cN_{1}[A_1^{k-1}](x)> c\,\|\cN_1[A_1^{k-1}]\|_{BMO_{C}} \}\:,
\eeq
can be written as a finite union of disjoint dyadic intervals with
\beq\label{kA11}
A_1^{k}\prec_{c} A_1^{k-1}\,.
\eeq
\end{itemize}

This process will end in a finite number of steps since the family $\P$ is finite.

\begin{obs}\label{bmobehav}
1) Define $\I_{1}^{max}:=\bigcup_{k\geq 0}\I_{1}^{max}[A_1^{k}]$ and let the \textit{global counting function of order one} be
\beq\label{gcntf1}
\cN_{1}:=\sum_{I\in\I_{1}^{max}}\chi_{I}\,.
\eeq

Notice that as a consequence of the above construction we have
\beq\label{bmoN}
\|\cN_{1}\|_{BMO_{C}}\lesssim \max_{k}\|\cN_{1}[A_1^{k}]\|_{BMO_{C}}\lesssim 1\:.
\eeq

2)  For any $0\leq l\leq k$, we have that $A_1^{k}\subseteq A_1^{l}$ with
\beq\label{boundA}
A_1^{k}\prec_{(k-l)\:c}\,A_1^{l}\:.
\eeq
\end{obs}

$\newline$

\noindent\textbf{Stage 1.2  Construction of the sets $\{\p_1[A_1^{k}]\}_{k\geq 1}$}
$\newline$

As mentioned above, we will associate to each of the sets within
$$\{A_1^k\}_{k\geq 0}\,,$$
constructed at Stage 1, a corresponding collection of tiles $\p_1[A_1^{k}]$.

Our construction process follows an ascending induction pattern.
\begin{itemize}
\item for $k=0$, define
\beq\label{setp1a1}
\p_{1}[A_1^{0}]:=\{P=[\vec{\a},I]\in\P\,|\,I\nsubseteq A_{1}^{1}\:\:\&\:\:A_{\P,A_1^{0}}(P)\in (2^{-1}, 2^{0}]\}\:.
\eeq

\item for general $k\in\N$, we set
\beq\label{setp1ak}
\p_{1}[A_1^{k}]:=\left\{P=[\vec{\a},I]\in\P\,|\,\begin{array}{rl} I\nsubseteq A_{1}^{k+1}\:,I\subseteq A_{1}^{k}\:\:\:\:\\
A_{\P,A_{1}^{k}}(P)\in (2^{-1}, 2^{0}]
 \end{array}\right\}\:.
\eeq

\item Finally, we now define
\textit{the collection of tiles of mass (of order) $1$} as
\beq\label{PP1}
\P_1:=\bigcup_{k\geq 0}\p_{1}[A_1^{k}]\:.
\eeq
\end{itemize}

Here the construction of the $1-$mass set ends.

$\newline$
\noindent\textbf{Step 2. Constructing the family $\P_2$}
$\newline$

First let us notice that from Step 1 we realized the following partition of the interval $[0,1]$:
\beq\label{part}
[0,1]=\bigcup_{k\geq 0} A_{1}^{k}\setminus A_1^{k+1}\;.
\eeq

Given $j_1\in\N$, we fix the set\footnote{Throughout the section we assume that all the quantities/objects considered here are non-trivial (i.e. the sets of intervals or tiles are non-empty and the norms involved are non-zero). If that is not the case then the algorithm stops.} $A_{1}^{j_1}\setminus A_1^{j_1+1}$. In what follows we will adapt the reasonings described at Step 1 to this specific set in order to construct the corresponding set of tiles $\P_2[A_{1}^{j_1}\diamond  A_1^{j_1+1}]$.

 As before, we will have two stages:
\begin{itemize}
\item Stage 2.1 - define a finite sequence of nested sets $\{A_2^{k}[A_{1}^{j_1}\diamond  A_{1}^{j_1+1}]\}_{k}$ such that
\begin{itemize}
\item  $A_2^{k}[A_{1}^{j_1}\diamond  A_{1}^{j_1+1}]$ is a finite union of maximal dyadic intervals;

\item  $A_2^{k}[A_{1}^{j_1}\diamond  A_{1}^{j_1+1}]\prec A_{1}^{j_1}$;

\item for any $J\subseteq A_{1}^{j_1+1}$ dyadic we either have $J\subset A_2^{k}[A_{1}^{j_1}\diamond  A_{1}^{j_1+1}]$
or $J\cap A_2^{k}[A_{1}^{j_1}\diamond  A_{1}^{j_1+1}]=\emptyset$.
\end{itemize}
\item Stage 2.2 - for each set $A_2^{k}[A_{1}^{j_1}\diamond  A_{1}^{j_1+1}]$ construct a corresponding family of tiles $\p_2(A_2^{k}[A_{1}^{j_1}\diamond  A_{1}^{j_1+1}])$ such that
\begin{itemize}
\item  $\p_2(A_2^{k}[A_{1}^{j_1}\diamond  A_{1}^{j_1+1}])$ is a convex family of tiles;

\item  each tile $P=[\vec{\a},I]\in \p_2(A_2^{k}[A_{1}^{j_1}\diamond  A_{1}^{j_1+1}])$ has the properties
\beq\label{masstwo}
\begin{array}{cc}
I\subseteq  A_{1}^{j_1}\:\:\&\:\:I\nsubseteq  A_{1}^{j_1+1}\:\:\:\textrm{and in fact also}\,\\\\
I\subseteq A_2^{k}[A_{1}^{j_1}\diamond  A_{1}^{j_1+1}]\:\:\&\:\:I\nsubseteq A_2^{k+1}[A_{1}^{j_1}\diamond  A_{1}^{j_1+1}]\,,\\\\
2^{-2}<A_{\P, A_2^{k}[A_{1}^{j_1}\diamond  A_{1}^{j_1+1}]}(P)\leq 2^{-1}\:.
\end{array}
\eeq
\end{itemize}
\end{itemize}

With these, we have:
$\newline$

\noindent\textbf{Stage 2.1 Construction of the sets $\{A_2^{k}[A_{1}^{j_1}\diamond  A_{1}^{j_1+1}]\}_{k\geq 1}$}
$\newline$

\noindent \textbf{2.1.1 Construction of the set $A_2^{1}[A_{1}^{j_1}\diamond  A_{1}^{j_1+1}]$}

\begin{itemize}
\item We start by setting
\beq\label{set102}
A_2^{0}[A_{1}^{j_1}\diamond  A_{1}^{j_1+1}]:=A_{1}^{j_1}\:.
\eeq
\item Let
\beq\label{max02}
\p_{2}^{max}[A_2^{0}[A_{1}^{j_1}\diamond  A_{1}^{j_1+1}]]:=\left\{\begin{array}{cl} P=[\vec{\a},I]\in\P\\
I\subseteq A_2^{0}[A_{1}^{j_1}\diamond  A_{1}^{j_1+1}]\\ I\nsubseteq  A_{1}^{j_1+1}
\end{array} \,\big|\,\begin{array}{cl} P\:\textrm{maximal}\\\frac{|E(P)|}{|I|}> 2^{-2}\end{array}\right\}\:.
\eeq
\item Define
\beq\label{c102}
\C_{2}[A_2^{0}[A_{1}^{j_1}\diamond  A_{1}^{j_1+1}]]:=\sum_{P\in \p_{2}^{max}[A_2^{0}[A_{1}^{j_1}\diamond  A_{1}^{j_1+1}]]} \chi_{E(P)}\,,
\eeq
and deduce that
\beq\label{c1sm02}
\C_{2}[A_2^{0}[A_{1}^{j_1}\diamond  A_{1}^{j_1+1}]](x)\leq 2^2\:\:\:\forall\:x\in[0,1]\;.
\eeq
\item Set
\beq\label{i102}
\I_{2}^{max}[A_2^{0}[A_{1}^{j_1}\diamond  A_{1}^{j_1+1}]]:=\{I\,|\,P=[\vec{\a},I]\in\p_{2}^{max}[A_2^{0}[A_{1}^{j_1}\diamond  A_{1}^{j_1+1}]]\}\,.
\eeq
\item Define the \textit{counting function of order one adapted to $A_2^{0}[A_{1}^{j_1}\diamond  A_{1}^{j_1+1}]$} as
\beq\label{countf02}
\cN_{2}[A_2^{0}[A_{1}^{j_1}\diamond  A_{1}^{j_1+1}]]:=\sum_{I\in\I_{2}^{max}[A_2^{0}[A_{1}^{j_1}\diamond  A_{1}^{j_1+1}]]}\chi_{I}\,.
\eeq
and deduce that
\beq\label{defbmoc02}
\|\cN_{2}[A_2^{0}[A_{1}^{j_1}\diamond  A_{1}^{j_1+1}]]\|_{BMO_{C}}:=\sup_{J\textrm{dyadic}\atop{J\subseteq [0,1]}}\,\frac{\sum_{I\subseteq J\atop{I\in \I_{2}^{max}[A_2^{0}[A_{1}^{j_1}\diamond  A_{1}^{j_1+1}]]}}|I|}{|J|}\leq 2^2\:,
\eeq
and
\beq\label{defbmodc02}
\|\cN_{2}[A_2^{0}[A_{1}^{j_1}\diamond  A_{1}^{j_1+1}]]\|_{BMO_{D}}\leq 2\|\cN_{2}[A_2^{0}[A_{1}^{j_1}\diamond  A_{1}^{j_1+1}]]\|_{BMO_{C}}\,.
\eeq
\item  Applying now John-Nirenberg deduce that the set\footnote{While the value of $c$ here is irrelevant as long as we assume is bounded from below away from the origin by $c_0=10^{10}$ say, the meaning of $2c$ (or later $n c$) in this context is to point out that we assume $2c>2c_0$.}
\beq\label{JNP2}
A_2^{1}[A_{1}^{j_1}\diamond  A_{1}^{j_1+1}]:=\{x\in A_2^{0}[A_{1}^{j_1}\diamond  A_{1}^{j_1+1}]\,|\,\cN_{2}[A_2^{0}[A_{1}^{j_1}\diamond  A_{1}^{j_1+1}]](x)>2 c\,2^2\}\,,
\eeq
is a finite union of disjoint dyadic intervals with
\beq\label{kA02}
A_2^{1}[A_{1}^{j_1}\diamond  A_{1}^{j_1+1}]\prec_{2c} A_2^{0}[A_{1}^{j_1}\diamond  A_{1}^{j_1+1}]\,.
\eeq
\end{itemize}

$\newline$
\noindent \textbf{2.1.2 Construction of the sets $A_2^{k}[A_{1}^{j_1}\diamond  A_{1}^{j_1+1}]$ with $k\geq 1$}
$\newline$

We assume by induction that, for some $k\in\N$, $k\geq 1$, we have constructed a set $A_2^{k-1}[A_{1}^{j_1}\diamond  A_{1}^{j_1+1}]$ that can be represented as a finite union of maximal (disjoint) dyadic intervals.

\begin{itemize}
\item Repeating the above construction we define the collection of maximal tiles
\beq\label{max122}
\p_{2}^{max}[A_2^{k-1}[A_{1}^{j_1}\diamond  A_{1}^{j_1+1}]]=\left\{\begin{array}{cl} P=[\vec{\a},I]\in\P\\
I\subseteq A_2^{k-1}[A_{1}^{j_1}\diamond  A_{1}^{j_1+1}]\\ I\nsubseteq  A_{1}^{j_1+1}
\end{array}\big|\,\begin{array}{cc} P\:\textrm{maximal}\\ \frac{|E(P)|}{|I|}> 2^{-2}\end{array}\right\}\:.
\eeq

\item Set the collection of time-intervals of maximal tiles as
\beq\label{i222}
\I_{2}^{max}[A_2^{k-1}[A_{1}^{j_1}\diamond  A_{1}^{j_1+1}]:=\{I\,|\,P=[\vec{\a},I]\in\p_{2}^{max}[A_2^{k-1}[A_{1}^{j_1}\diamond  A_{1}^{j_1+1}]]\:\}\,.
\eeq

\item Define the \textit{counting function of order two adapted to $A_2^{k-1}[A_{1}^{j_1}\diamond  A_{1}^{j_1+1}]$} as
\beq\label{countf122}
\cN_{2}[A_2^{k-1}[A_{1}^{j_1}\diamond  A_{1}^{j_1+1}]]:=\sum_{I\in\I_{1}^{max}[A_2^{k-1}[A_{1}^{j_1}\diamond  A_{1}^{j_1+1}]]}\chi_{I}\,,
\eeq
and applying the same reasonings as in the previous situation, notice that $\cN_{2}[A_2^{k-1}[A_{1}^{j_1}\diamond  A_{1}^{j_1+1}]]\in BMO_{D}$ and moreover that $\|\cN_{2}[A_2^{k-1}[A_{1}^{j_1}\diamond  A_{1}^{j_1+1}]]\|_{BMO_{C}}\leq 2^2$.

\item Define now
\beq\label{defbmoc12}
A_2^{k}[A_{1}^{j_1}\diamond  A_{1}^{j_1+1}]:=\{x\in A_2^{k-1}[A_{1}^{j_1}\diamond  A_{1}^{j_1+1}]\,|\,\cN_{2}[A_2^{k-1}[A_{1}^{j_1}\diamond  A_{1}^{j_1+1}]](x)>2c\,2^2\}
\eeq
and remark that
\begin{equation}\label{a2}
A_2^{k}[A_{1}^{j_1}\diamond  A_{1}^{j_1+1}]\prec_{2c} A_2^{k-1}[A_{1}^{j_1}\diamond  A_{1}^{j_1+1}]\:.
\end{equation}
\end{itemize}

This process will end in a finite number of steps.

$\newline$

\noindent\textbf{Stage 2.2  Construction of the sets $\{\p_2[A_2^{k}[A_{1}^{j_1}\diamond  A_{1}^{j_1+1}]]\}_{k\in \N}$}
$\newline$

\begin{itemize}
\item For $k\in\N$, we simply define
\beq\label{setp1a122}
\p_{2}[A_2^{k}[A_{1}^{j_1}\diamond  A_{1}^{j_1+1}]]:=\left\{\begin{array}{cl} P\in\P\setminus\P_1\\P=[\vec{\a},I] \end{array}\,\big|\,\begin{array}{rl} I\subseteq A_2^{k}[A_{1}^{j_1}\diamond  A_{1}^{j_1+1}]\:\\
I\nsubseteq A_2^{k+1}[A_{1}^{j_1}\diamond  A_{1}^{j_1+1}]\:\:\&\:\:I\nsubseteq A_{1}^{j_1+1}\:\\
A_{\P,A_2^{k}[A_{1}^{j_1}\diamond  A_{1}^{j_1+1}]}(P)\in (2^{-2}, 2^{-1}] \end{array}\right\}\:.
\eeq

\item Next, we let
\beq\label{setp1a122t}
\P_2[A_{1}^{j_1}\diamond  A_{1}^{j_1+1}]:=\bigcup_{k\in\N}\p_{2}[A_2^{k}[A_{1}^{j_1}\diamond  A_{1}^{j_1+1}]]\:.
\eeq

\item Finally, we now define
\textit{the collection of tiles of mass (of order) $2$} as
\beq\label{PP22}
\P_2:=\bigcup_{j\in\N}\P_2[A_{1}^{j}\diamond  A_{1}^{j+1}]\:.
\eeq
\end{itemize}

Here the construction of the $2-$mass set ends.

$\newline$
\noindent\textbf{Step $n$. Constructing the family $\P_n$}
$\newline$

Suppose now following the above algorithm (via induction), that for  $j_1,\,j_2,\ldots j_{n-1}\in \N$ and $n\geq 3$, we have constructed the sets
\begin{equation}\label{part}
A_{n-1}^{j_{n-1}}[A_{n-2}^{j_{n-2}}\diamond  A_{n-2}^{j_{n-2}+1}]\ldots[A_{1}^{j_1}\diamond  A_{1}^{j_1+1}]\,.
\end{equation}

 As before, for the $n^{th}$ step we will have two stages:
\begin{itemize}
\item Stage n.1 - define a finite sequence of nested sets
$$\{A_{n}^{k}[A_{n-1}^{j_{n-1}}\diamond A_{n-1}^{j_{n-1}+1}]\ldots[A_{1}^{j_1}\diamond  A_{1}^{j_1+1}]\}_{k}\,,$$
such that
\begin{itemize}

\item  $A_{n}^{k}[A_{n-1}^{j_{n-1}}\diamond A_{n-1}^{j_{n-1}+1}]\ldots[A_{1}^{j_1}\diamond  A_{1}^{j_1+1}]$ is a finite union of maximal (disjoint) dyadic intervals;

\item  $A_{n}^{k}[A_{n-1}^{j_{n-1}}\diamond A_{n-1}^{j_{n-1}+1}]\ldots[A_{1}^{j_1}\diamond  A_{1}^{j_1+1}]\prec A_{n-1}^{j_{n-1}}[A_{n-2}^{j_{n-2}}\diamond A_{n-2}^{j_{n-2}+1}]\ldots[A_{1}^{j_1}\diamond  A_{1}^{j_1+1}]$;

\item for any\footnote{If $s=n-1$, we set $A_{n-s}^{j_{n-s}+1}[A_{n-s-1}^{j_{n-s-1}}\diamond A_{n-s-1}^{j_{n-s-1}+1}]\ldots[A_{1}^{j_1}\diamond  A_{1}^{j_1+1}]:= A_{1}^{j_1+1}$.} $J\subseteq A_{n-s}^{j_{n-s}+1}[A_{n-s-1}^{j_{n-s-1}}\diamond A_{n-s-1}^{j_{n-s-1}+1}]\ldots[A_{1}^{j_1}\diamond  A_{1}^{j_1+1}]$ dyadic interval and any $s\in\{1,\ldots,n-1\}$ we either have
    $$J\subset A_{n}^{k}[A_{n-1}^{j_{n-1}}\diamond A_{n-1}^{j_{n-1}+1}]\ldots[A_{1}^{j_1}\diamond  A_{1}^{j_1+1}]\,\:\:\:\:\:\textrm{or}$$
$$J\cap A_{n}^{k}[A_{n-1}^{j_{n-1}}\diamond A_{n-1}^{j_{n-1}+1}]\ldots[A_{1}^{j_1}\diamond  A_{1}^{j_1+1}]=\emptyset\,.$$
\end{itemize}

\item Stage n.2 - for each set $A_{n}^{k}[A_{n-1}^{j_{n-1}}\diamond A_{n-1}^{j_{n-1}+1}]\ldots[A_{1}^{j_1}\diamond  A_{1}^{j_1+1}]$ construct a corresponding family of tiles
    $$\p_n(A_{n}^{k}[A_{n-1}^{j_{n-1}}\diamond A_{n-1}^{j_{n-1}+1}]\ldots[A_{1}^{j_1}\diamond  A_{1}^{j_1+1}])$$
    such that

\begin{itemize}
\item  $\p_n(A_{n}^{k}[A_{n-1}^{j_{n-1}}\diamond A_{n-1}^{j_{n-1}+1}]\ldots[A_{1}^{j_1}\diamond  A_{1}^{j_1+1}])$ is a convex family;

\item  each tile $P=[\vec{\a},I]\in \p_n(A_{n}^{k}[A_{n-1}^{j_{n-1}}\diamond A_{n-1}^{j_{n-1}+1}]\ldots[A_{1}^{j_1}\diamond  A_{1}^{j_1+1}])$ has the properties
\beq\label{massn}
\begin{array}{cc}
I\subseteq A_{n}^{k}[A_{n-1}^{j_{n-1}}\diamond A_{n-1}^{j_{n-1}+1}]\ldots[A_{1}^{j_1}\diamond  A_{1}^{j_1+1}] \,,\\\\
I\nsubseteq A_{n}^{k+1}[A_{n-1}^{j_{n-1}}\diamond A_{n-1}^{j_{n-1}+1}]\ldots[A_{1}^{j_1}\diamond  A_{1}^{j_1+1}]\:\textrm{and also}\\\\
I\nsubseteq A_{n-s}^{j_{n-s}+1}[A_{n-s-1}^{j_{n-s-1}}\diamond A_{n-s-1}^{j_{n-s-1}+1}]\ldots[A_{1}^{j_1}\diamond  A_{1}^{j_1+1}]\:\:\forall\:\:s\in\{1,\ldots, n-1\}\,,\\\\
2^{-n}<A_{\P, A_{n}^{k}[A_{n-1}^{j_{n-1}}\diamond A_{n-1}^{j_{n-1}+1}]\ldots[A_{1}^{j_1}\diamond  A_{1}^{j_1+1}]}(P)\leq 2^{-n+1}\:.
\end{array}
\eeq
\end{itemize}
\end{itemize}

$\newline$
\noindent \textbf{Stage n.1 Construction of the sets $\{A_{n}^{k}[A_{n-1}^{j_{n-1}}\diamond A_{n-1}^{j_{n-1}+1}]\ldots[A_{1}^{j_1}\diamond  A_{1}^{j_1+1}]\}_{k\in\N}$}
$\newline$

\begin{itemize}
\item As at the Step 1, we will proceed by induction:
\begin{itemize}
\item When $k=0$, we simply set
\beq\label{setn0n}
A_n^{0}[A_{n-1}^{j_{n-1}}\diamond A_{n-1}^{j_{n-1}+1}]\ldots[A_{1}^{j_1}\diamond  A_{1}^{j_1+1}]:=A_{n-1}^{j_{n-1}}[A_{n-2}^{j_{n-2}}\diamond  A_{n-2}^{j_{n-2}+1}]\ldots[A_{1}^{j_1}\diamond  A_{1}^{j_1+1}]\:.
\eeq
\item For $k\geq 1$ we assume we have constructed
$$A_{n}^{k-1}[A_{n-1}^{j_{n-1}}\diamond A_{n-1}^{j_{n-1}+1}]\ldots[A_{1}^{j_1}\diamond  A_{1}^{j_1+1}]\,.$$
\end{itemize}
\item Next we identify the collection of maximal tiles
\beq\label{maxn1}
\begin{array}{cc}
\p_{n}^{max}[A_{n}^{k-1}[A_{n-1}^{j_{n-1}}\diamond A_{n-1}^{j_{n-1}+1}]\ldots[A_{1}^{j_1}\diamond  A_{1}^{j_1+1}]]:=\\\\
\left\{P=[\vec{\a},I]\in\P\,\big|\,\begin{array}{cc} I\subseteq A_{n}^{k-1}[A_{n-1}^{j_{n-1}}\diamond A_{n-1}^{j_{n-1}+1}]\ldots[A_{1}^{j_1}\diamond  A_{1}^{j_1+1}]\\
I\nsubseteq A_{n-s}^{j_{n-s}+1}[A_{n-s-1}^{j_{n-s-1}}\diamond A_{n-s-1}^{j_{n-s-1}+1}]\ldots[A_{1}^{j_1}\diamond  A_{1}^{j_1+1}]\:\forall\: s<n\\
P\:\textrm{maximal}\:\:\&\:\:\frac{|E(P)|}{|I|}> 2^{-n}
\end{array}\right\}\:.
\end{array}
\eeq
\item Define
\beq\label{c102n}
\begin{array}{cc}
\C_{n}[A_{n}^{k-1}[A_{n-1}^{j_{n-1}}\diamond A_{n-1}^{j_{n-1}+1}]\ldots[A_{1}^{j_1}\diamond  A_{1}^{j_1+1}]]\\\\
:=\sum_{P\in \p_{n}^{max}[A_{n}^{k-1}[A_{n-1}^{j_{n-1}}\diamond A_{n-1}^{j_{n-1}+1}]\ldots[A_{1}^{j_1}\diamond  A_{1}^{j_1+1}]]} \chi_{E(P)}\,,
\end{array}
\eeq
and deduce that
\beq\label{c1sm02n}
\C_{n}[A_{n}^{k-1}[A_{n-1}^{j_{n-1}}\diamond A_{n-1}^{j_{n-1}+1}]\ldots[A_{1}^{j_1}\diamond  A_{1}^{j_1+1}]](x)\leq 2^n\:\:\:\forall\:x\in[0,1]\;.
\eeq

\item Set the collection of time-intervals of maximal tiles as
\beq\label{in2}
\begin{array}{cc}
\I_{n}^{max}[A_{n}^{k-1}[A_{n-1}^{j_{n-1}}\diamond A_{n-1}^{j_{n-1}+1}]\ldots[A_{1}^{j_1}\diamond  A_{1}^{j_1+1}]]\\\\
:=\{I\,|\,P=[\vec{\a},I]\in\p_{n}^{max}[A_{n}^{k-1}[A_{n-1}^{j_{n-1}}\diamond A_{n-1}^{j_{n-1}+1}]\ldots[A_{1}^{j_1}\diamond  A_{1}^{j_1+1}]]\:\}\,.
\end{array}
\eeq

\item Define the \textit{counting function of order $n$ adapted to}
$$A_{n}^{k-1}[A_{n-1}^{j_{n-1}}\diamond A_{n-1}^{j_{n-1}+1}]\ldots[A_{1}^{j_1}\diamond  A_{1}^{j_1+1}]$$
as
\beq\label{countf1n}
\begin{array}{cc}
\cN_{n}[A_{n}^{k-1}[A_{n-1}^{j_{n-1}}\diamond A_{n-1}^{j_{n-1}+1}]\ldots[A_{1}^{j_1}\diamond  A_{1}^{j_1+1}]]\\\\
:=\sum_{I\in\I_{n}^{max}[A_{n}^{k-1}[A_{n-1}^{j_{n-1}}\diamond A_{n-1}^{j_{n-1}+1}]\ldots[A_{1}^{j_1}\diamond  A_{1}^{j_1+1}]]}\chi_{I}\,,
\end{array}
\eeq
and notice that
$$\cN_{n}[A_{n}^{k-1}[A_{n-1}^{j_{n-1}}\diamond A_{n-1}^{j_{n-1}+1}]\ldots[A_{1}^{j_1}\diamond  A_{1}^{j_1+1}]]\in BMO_{D}\,$$ with
\beq\label{contnn}
\|\cN_{n}[A_{n}^{k-1}[A_{n-1}^{j_{n-1}}\diamond A_{n-1}^{j_{n-1}+1}]\ldots[A_{1}^{j_1}\diamond  A_{1}^{j_1+1}]]\|_{BMO_{C}}\leq 2^n\,.
\eeq

\item Applying now the the John-Nirenberg inequality, for
$$\g> c\,\|\cN_n[A_{n}^{k-1}[A_{n-1}^{j_{n-1}}\diamond A_{n-1}^{j_{n-1}+1}]\ldots[A_{1}^{j_1}\diamond  A_{1}^{j_1+1}]]\|_{BMO_{C}}\,,$$ we have
\beq\label{JNP1n}
|\{x\in J\,|\,\sum_{I\subseteq J\atop{I\in\I_{n}^{max}[A_{n}^{k-1}[A_{n-1}^{j_{n-1}}\diamond A_{n-1}^{j_{n-1}+1}]\ldots[A_{1}^{j_1}\diamond  A_{1}^{j_1+1}]]}}\chi_{I}(x)>\gamma\}|\lesssim |J|\,e^{-c}\:.
\eeq

\item Define now the set
\beq\label{A1n1}
\begin{array}{cc}
A_{n}^{k}[A_{n-1}^{j_{n-1}}\diamond A_{n-1}^{j_{n-1}+1}]\ldots[A_{1}^{j_1}\diamond  A_{1}^{j_1+1}]\\\\
:=\{x\in [0,1]\,|\, \frac{\cN_n[A_{n}^{k-1}[A_{n-1}^{j_{n-1}}\diamond A_{n-1}^{j_{n-1}+1}]\ldots[A_{1}^{j_1}\diamond  A_{1}^{j_1+1}]](x)}{\|\cN_n[A_{n}^{k-1}[A_{n-1}^{j_{n-1}}\diamond A_{n-1}^{j_{n-1}+1}]\ldots[A_{1}^{j_1}\diamond  A_{1}^{j_1+1}]]\|_{BMO_{C}}}> n\,c\}\:,
\end{array}
\eeq
and notice that $A_{n}^{k}[A_{n-1}^{j_{n-1}}\diamond A_{n-1}^{j_{n-1}+1}]\ldots[A_{1}^{j_1}\diamond  A_{1}^{j_1+1}]$ is a finite union of maximal disjoint dyadic intervals with
\beq\label{kAn11}
A_{n}^{k}[A_{n-1}^{j_{n-1}}\diamond A_{n-1}^{j_{n-1}+1}]\ldots[A_{1}^{j_1}\diamond  A_{1}^{j_1+1}]\prec_{n\,c}A_{n}^{k-1}[A_{n-1}^{j_{n-1}}\diamond A_{n-1}^{j_{n-1}+1}]\ldots[A_{1}^{j_1}\diamond  A_{1}^{j_1+1}]\,.
\eeq
\end{itemize}

This process will end in a finite number of steps.

\begin{obs}\label{bmobehavn}
1) Define $\I_{n}^{max}:=\bigcup_{k}\bigcup_{j_1,\ldots j_{n-1}\in\N}\I_{n}^{max}[A_{n}^{k-1}[A_{n-1}^{j_{n-1}}\diamond A_{n-1}^{j_{n-1}+1}]\ldots[A_{1}^{j_1}\diamond  A_{1}^{j_1+1}]]$ and let the \textit{global counting function of order $n$} be
\beq\label{gcntfn}
\cN_{n}:=\sum_{I\in\I_{n}^{max}}\chi_{I}\,.
\eeq

Notice that as a consequence of the above construction we have
\beq\label{bmoNn}
\|\cN_{n}\|_{BMO_{C}}\lesssim n\,2^{n}\:.
\eeq

2)  For any $k,l,j_1,\ldots, j_{n-1}$ with $k\geq l$ we have that
\beq\label{boundAn}
|A_{n}^{k}[A_{n-1}^{j_{n-1}}\diamond A_{n-1}^{j_{n-1}+1}]\ldots[A_{1}^{j_1}\diamond  A_{1}^{j_1+1}]|\prec_{(k-l)\,n\,c}
|A_{n}^{l}[A_{n-1}^{j_{n-1}}\diamond A_{n-1}^{j_{n-1}+1}]\ldots[A_{1}^{j_1}\diamond  A_{1}^{j_1+1}]|\:.
\eeq
\end{obs}

$\newline$
\noindent \textbf{Stage n.2. Construction of the sets $\{\p_{n}[A_{n}^{k}[A_{n-1}^{j_{n-1}}\diamond  A_{n-1}^{j_{n-1}+1}]\ldots[A_{1}^{j_1}\diamond  A_{1}^{j_1+1}]]\}_{k\in\N}$}
$\newline$

\begin{itemize}
\item For $k\in\N$,  we define
\beq\label{setp1a122n}
\begin{array}{cc}
\p_{n}[A_{n}^{k}[A_{n-1}^{j_{n-1}}\diamond  A_{n-1}^{j_{n-1}+1}]\ldots[A_{1}^{j_1}\diamond  A_{1}^{j_1+1}]]:=\\\\
\left\{\begin{array}{cl} P\in\P\setminus\bigcup_{j=1}^{n-1}\P_j\\P=[\vec{\a},I] \end{array}\,\big|\,\begin{array}{rl} I\subseteq A_{n}^{k}[A_{n-1}^{j_{n-1}}\diamond A_{n-1}^{j_{n-1}+1}]\ldots[A_{1}^{j_1}\diamond  A_{1}^{j_1+1}]\:\\
I\nsubseteq A_{n}^{k+1}[A_{n-1}^{j_{n-1}}\diamond  A_{n-1}^{j_{n-1}+1}]\ldots[A_{1}^{j_1}\diamond  A_{1}^{j_1+1}]\:\\
I\nsubseteq A_{n-s}^{j_{n-s}+1}[A_{n-s-1}^{j_{n-s-1}}\diamond A_{n-s-1}^{j_{n-s-1}+1}]\ldots[A_{1}^{j_1}\diamond  A_{1}^{j_1+1}]\:\forall\: s<n\\
A_{\P,A_{n}^{k}[A_{n-1}^{j_{n-1}}\diamond  A_{n-1}^{j_{n-1}+1}]\ldots[A_{1}^{j_1}\diamond  A_{1}^{j_1+1}]}(P)\in (2^{-n}, 2^{-n+1}] \end{array}\right\}\:.
\end{array}
\eeq

\item Next, we set
\beq\label{pna}
\begin{array}{cc}
\P_n[A_{n-1}^{j_{n-1}}\diamond  A_{n-1}^{j_{n-1}+1}]\ldots[A_{1}^{j_1}\diamond  A_{1}^{j_1+1}]:=\\\\
\bigcup_{k}\p_{n}[A_{n}^{k}[A_{n-1}^{j_{n-1}}\diamond A_{n-1}^{j_{n-1}+1}]\ldots[A_{1}^{j_1}\diamond  A_{1}^{j_1+1}]]\:.
\end{array}
\eeq
\item Finally, we define \textit{the collection of tiles of mass (of order) $n$}
\beq\label{Pnd}
\P_n:=\bigcup_{j_1,\ldots j_{n-1}\in\N}\P_n[A_{n-1}^{j_{n-1}}\diamond  A_{n-1}^{j_{n-1}+1}]\ldots[A_{1}^{j_1}\diamond  A_{1}^{j_1+1}]\:.
\eeq
\end{itemize}

Here the construction of the $n-$mass set ends.

Remark that, from the above algorithm, we have
\beq\label{P}
\P=\bigcup_{n\geq0}\P_n\:.
\eeq

This ends the partition of our set $\P$.

\subsection{\bf Main Proposition; ending the proof}
$\newline$

In what follows, we will state the key result on which our theorem is based. The proof of this proposition will be postponed for the next sections. With the notations form the previous section, we have
$\newline$

\noindent {\bf Main Proposition.}
{\it Fix $n\in\N$. Then there exist a constant $\eta=\eta(d)\in(0,\frac{1}{2})$ depending only
on $d$ such that
$$\left\|T^{{\P}_n}f\right\|_{p}\lesssim_{p,d} \,2^{-n\,\eta(1-\frac{1}{p^*})}\,\left\|f\right\|_p\:,$$
for all $f\in L^p(\TT)$.}
$\newline$

If we believe this for the moment, then we trivially have
$$\left\|Tf\right\|_{p}\leq\sum_n \left\|T^{{\P}_n}f\right\|_{p}\lesssim_{p,d} \sum_n 2^{-n\,\eta(1-\frac{1}{p^*})}\,\left\|f\right\|_p
\lesssim_{p,d} \,\left\|f\right\|_p\:.$$

\section{\bf Reduction of the main proposition}\label{Reductmainprop}

In this section, we will present the strategy needed to prove our main proposition.

\subsection{\bf Preparatives} We first introduce the main concepts that will play the central role in the analysis
of the collections of tiles $\P_n$, $n\in\N$.

 \begin{d0}\label{tree}[\textsf{Tree}]

 We say that a set of tiles $\p\subset\P$ is a \textbf{tree} (relative to $\leq$) with \emph{top} $P_0$ if
the following conditions\footnote{To avoid the boundary problems arising
from the useage of a single dyadic
grid and from the definition of our tiles, we will often involve in our reasonings a dilation factor (of the tiles).}
are satisfied: $\newline
1)\:\:\:\:\:\forall\:\:P\in\p\:\:\:\Rightarrow\:\:\:\:\frac{3}{2}P\lneq
10 P_0 $ $\newline 2)\:\:\:\:\:$if $P\in\p$ and $P'\in N(P)$ such
that $\frac{4}{3}P'\lneq 10 P_{0}$ then $P'\in\p$ $\newline
3)\:\:\:\:\:$if $P_1,\:P_2\: \in\p$ and $P_1\leq P \leq P_2$ then
$P\in\p\:.$
\end{d0}

\begin{d0}\label{stree}[\textsf{Sparse tree}]

Let $C>0$ be an absolute constant. We say that a set of tiles $\p\subset\P$ is a \textbf{$C-$sparse tree} if
 $\p$ is a tree and for any $P=[\vec{\a},I]\in\p$ we have
 \beq\label{cmstree}
 \sum_{{P'=[\vec{\a}',I']\in\p}\atop{I'\subseteq I}}|I'|\leq C\,|I|\:.
 \eeq
In our later reasonings, the  specific value of the constant $C$ will be of no relevance\footnote{All the constants $C$'s appearing in this context will be bounded by a positive absolute constant possibly depending only on $d$.} and thus we will simply refer to a $C-$sparse tree as a \textbf{sparse tree}.
 \end{d0}

\begin{d0}\label{infforest} [$L^{\infty}-$\textsf{forest}]

Fix $n\in\N$. We say that $\p\subseteq\P_n$ is an $L^\infty$-\textbf{forest} of generation $n$ iff the folowing two conditions hold:
\begin{enumerate}

\item $\p$ is a collection of \textbf{separated} trees, {\it i.e.}
$$\p=\bigcup_{j\in\N}\p_j\,$$
with each $\p_j$ a tree with top $P_j=[\vec{\a}_j,I_j]$ and such that
 \beq\label{wellsep}
 \forall\:\:k\not=j\:\:\&\:\:\forall\:\:P\in\p_j\:\:\:\:\:\:\:\:2P\nleq 2P_k\:.
 \eeq

\item the $\p-$\textbf{counting function}
\beq\label{cfunct}
 \cN_{\p}(x):=\sum_{j}\chi_{I_j}(x)
 \eeq
 obeys the estimate $\|\cN_{\p}\|_{L^{\infty}}\lesssim 2^n$.
\end{enumerate}

 Further on, if $\p\subseteq\P_n$ only consists of sparse separated trees then
we refer at $\p$ as a \textbf{sparse $L^\infty$-forest}.
 \end{d0}

 \begin{d0}\label{forest}[$BMO-$\textsf{forest}]

  A set  $\p\subseteq\P_n$ is called a  $BMO$-\textbf{forest} of generation $n$ or
 just simply a \textbf{forest}\footnote{When the context is clear we may no
 longer specify the order of the generation.} iff the following holds
\begin{enumerate}
\item $\p$ may be written as
\beq\label{for1}
\p=\bigcup_{j\in\N}\p_j\,,
\eeq
with each $\p_j$ an $L^\infty$-\textbf{forest} (of generation $n$);

\item for any $P\in\p_j$ and $P'\in\p_k$ with $j,k\in\N$,  $j<k$ we either have $I_P\cap I_{P'}=\emptyset$ or\footnote{The base $2$ here has no relevance. One could replace it with any $c>1$ so that \eqref{for2} transforms into $|I_{P'}|\leq c^{j-k}\,|I_{P}|$. More generally, it is in fact enough for the collection $\{I_{P_k}\}_k$ to obey a Carleson packing condition.}
\beq\label{for2}
|I_{P'}|\leq 2^{j-k}\,|I_{P}|\:.
\eeq
\end{enumerate}

 As before, if $\p\subseteq\P_n$ only consists of sparse $L^\infty$-forests, then,
we refer at $\p$ as a \textbf{sparse forest}.
\end{d0}

\begin{obs} Notice that if $\p\subseteq\P_n$ is a forest then, due to \eqref{for2} above, the counting function
\beq\label{cfunctbmo}
 \cN_{\p}:=\sum_{j}\cN_{\p_j}\,,
\eeq
obeys  the estimate
\beq\label{cfunctbmos}
\|\cN_{\p}\|_{BMO_{C}}\lesssim 2^n\,,
\eeq
hence the alternative name of the $BMO$-forest.

Also notice that if $\p\subseteq\P_n$ is a collection of separated trees then $\p$ is automatically a  ($BMO$-)forest.
\end{obs}

$\newline$

Now we can state the main results of this section; their proofs will
be postponed until Section 7.

\begin{p1}\label{prop1}[\textsf{Control over a sparse forest}]

Let $\p\subseteq\P_n$ be a sparse forest. Then there exists $\eta=\eta(d)\in(0,\frac{1}{2})$, depending only on the degree $d$,
such that for any $1<p<\infty$ we have
\beq\label{pro1}
\left\|T^{\p}\right\|_{p}\lesssim_{p,d} 2^{-n\,\eta\,(1-\frac{1}{p^*})}\:.
\eeq
\end{p1}

\begin{p1}\label{prop2}[\textsf{Control over a (general) forest}]

Let $\p\subseteq\P_n$ be a forest. Then there exists $\eta=\eta(d)\in(0,\frac{1}{2})$, depending only on the degree $d$,
such that for any $1<p<\infty$ we have
\beq\label{pro2}
\left\|T^{\p}\right\|_{p}\lesssim_{p,d} 2^{-n\,\eta\,(1-\frac{1}{p^*})}\:.
\eeq
\end{p1}

$\newline$
\subsection{\bf Reduction of the Main Proposition to Proposition \ref{prop2}}\label{redmp}
$\newline$

\textsf{Aim}: \emph{In this section our goal is to show that, for a fixed $n$, the set $\P_n$ can be roughly decomposed into a union of $\approx\,n$ forests.}

We start by recalling \eqref{setp1a122n} - \eqref{Pnd} and thus we have that
\beq\label{PPn1}
\P_n:=\bigcup_{k\in\N}\bigcup_{j_1,\ldots j_{n-1}\in\N} \p_{n}[A_{n}^{k}[A_{n-1}^{j_{n-1}}\diamond A_{n-1}^{j_{n-1}+1}]\ldots[A_{1}^{j_1}\diamond  A_{1}^{j_1+1}]]\:.
\eeq

We now make the following

\begin{claim}\label{CL}
For each $j_1,\ldots, j_{n-1}, k\in\N$ the set
\beq\label{clo}
\p_{n}[A_{n}^{k}[A_{n-1}^{j_{n-1}}\diamond A_{n-1}^{j_{n-1}+1}]\ldots[A_{1}^{j_1}\diamond  A_{1}^{j_1+1}]]
\eeq
can be decomposed in a union of at most $c\,n$ $L^\infty$-forests (of generation $n$)
\beq\label{cl1}
\{\p_{n}^{s}[A_{n}^{k}[A_{n-1}^{j_{n-1}}\diamond A_{n-1}^{j_{n-1}+1}]\ldots[A_{1}^{j_1}\diamond  A_{1}^{j_1+1}]]\}_{s\in\{1,\ldots , c\,n\}}\:,
\eeq
where here $c>0$ is some absolute constant.
\end{claim}

\begin{obs}
Notice that if we believe our claim for the moment, then denoting with
\beq\label{PPn1ob}
\P_n^s:=\bigcup_{k\in\N}\bigcup_{j_1,\ldots j_{n-1}\in\N} \p_{n}^s[A_{n}^{k}[A_{n-1}^{j_{n-1}}\diamond A_{n-1}^{j_{n-1}+1}]\ldots[A_{1}^{j_1}\diamond  A_{1}^{j_1+1}]]\:,
\eeq
we have that $\P_{n}^s$ is a BMO-forest. Indeed, this follows from the key condition \eqref{kAn11} in our construction of tiles and from the fact that each $\p_{n}^s[A_{n}^{k}[A_{n-1}^{j_{n-1}}\diamond A_{n-1}^{j_{n-1}+1}]\ldots[A_{1}^{j_1}\diamond  A_{1}^{j_1+1}]]$ is an $L^{\infty}$-forest.

Thus, since
 $$\P_n:=\bigcup_{s=0}^{c\,n} \P_{n}^s\:,$$
 we conclude that $\P_n$ can be written as a union of at most $c\,n$ forests as desired.
\end{obs}

We start by recalling the construction from the Step $n$ in our previous section.  Based on that algorithm, we are given the following:
\begin{itemize}
\item the set $A_{n}^{k}[A_{n-1}^{j_{n-1}}\diamond A_{n-1}^{j_{n-1}+1}]\ldots[A_{1}^{j_1}\diamond  A_{1}^{j_1+1}]$;

\item the collection of maximal tiles $\p_{n}^{max}[A_{n}^{k}[A_{n-1}^{j_{n-1}}\diamond A_{n-1}^{j_{n-1}+1}]\ldots[A_{1}^{j_1}\diamond  A_{1}^{j_1+1}]]$;

\item the collection of time intervals associated with the set of maximal tiles $\I_{n}^{max}[A_{n}^{k}[A_{n-1}^{j_{n-1}}\diamond A_{n-1}^{j_{n-1}+1}]\ldots[A_{1}^{j_1}\diamond  A_{1}^{j_1+1}]]$;

\item the counting function $\cN_{n}[A_{n}^{k}[A_{n-1}^{j_{n-1}}\diamond A_{n-1}^{j_{n-1}+1}]\ldots[A_{1}^{j_1}\diamond  A_{1}^{j_1+1}]]$;

\item the set of tiles $\p_{n}[A_{n}^{k}[A_{n-1}^{j_{n-1}}\diamond A_{n-1}^{j_{n-1}+1}]\ldots[A_{1}^{j_1}\diamond  A_{1}^{j_1+1}]]$.
\end{itemize}

We let $\bar{\p}_{n}^{max}[A_{n}^{k}[A_{n-1}^{j_{n-1}}\diamond A_{n-1}^{j_{n-1}+1}]\ldots[A_{1}^{j_1}\diamond  A_{1}^{j_1+1}]]$ be the set of all the maximal elements $P\in\p_{n}^{max}[A_{n}^{k}[A_{n-1}^{j_{n-1}}\diamond A_{n-1}^{j_{n-1}+1}]\ldots[A_{1}^{j_1}\diamond  A_{1}^{j_1+1}]]$ such that $\frac{|E(P)|}{|I_P|}> 2^{-n}$,  $I_P\subset A_{n}^{k}[A_{n-1}^{j_{n-1}}\diamond A_{n-1}^{j_{n-1}+1}]\ldots[A_{1}^{j_1}\diamond  A_{1}^{j_1+1}]$ and $I_P\nsubseteq A_{n}^{k+1}[A_{n-1}^{j_{n-1}}\diamond A_{n-1}^{j_{n-1}+1}]\ldots[A_{1}^{j_1}\diamond  A_{1}^{j_1+1}]$. Notice that
\beq\label{pmxc}
\begin{array}{cc}
\bar{\p}_{n}^{max}[A_{n}^{k}[A_{n-1}^{j_{n-1}}\diamond A_{n-1}^{j_{n-1}+1}]\ldots[A_{1}^{j_1}\diamond  A_{1}^{j_1+1}]]\\\\
=\p_{n}^{max}[A_{n}^{k}[A_{n-1}^{j_{n-1}}\diamond A_{n-1}^{j_{n-1}+1}]\ldots[A_{1}^{j_1}\diamond  A_{1}^{j_1+1}]]\\\\
\setminus\p_{n}^{max}[A_{n}^{k+1}[A_{n-1}^{j_{n-1}}\diamond A_{n-1}^{j_{n-1}+1}]\ldots[A_{1}^{j_1}\diamond  A_{1}^{j_1+1}]]\:.
\end{array}
\eeq
Moreover, from the previous tile partition algorithm and recalling \eqref{countf1n}, we remark that defining $$\bar{\I}_{n}^{max}[A_{n}^{k}[A_{n-1}^{j_{n-1}}\diamond A_{n-1}^{j_{n-1}+1}]\ldots[A_{1}^{j_1}\diamond  A_{1}^{j_1+1}]]$$
$$:=\{I\,|\,P=[\vec{\a},I]\in\bar{\p}_{n}^{max}[A_{n}^{k}[A_{n-1}^{j_{n-1}}\diamond A_{n-1}^{j_{n-1}+1}]\ldots[A_{1}^{j_1}\diamond  A_{1}^{j_1+1}]]\}$$
one has from \eqref{contnn}, \eqref{A1n1} and \eqref{pmxc} that the counting function
\beq\label{nkb}
\bar{\cN}_n[A_{n}^{k}[A_{n-1}^{j_{n-1}}\diamond A_{n-1}^{j_{n-1}+1}]\ldots[A_{1}^{j_1}\diamond  A_{1}^{j_1+1}]]:=\sum_{I\in\bar{\I}_{n}^{max}[A_{n}^{k}[A_{n-1}^{j_{n-1}}\diamond A_{n-1}^{j_{n-1}+1}]\ldots[A_{1}^{j_1}\diamond  A_{1}^{j_1+1}]]}\chi_{I}\,,
\eeq
obeys
\beq\label{contc}
\|\bar{\cN}_n[A_{n}^{k}[A_{n-1}^{j_{n-1}}\diamond A_{n-1}^{j_{n-1}+1}]\ldots[A_{1}^{j_1}\diamond  A_{1}^{j_1+1}]]\|_{L^{\infty}}\leq c\,n\,2^{n}\,.
\eeq

Fix throughout this section the values of
$k,\,j_1,\ldots,j_{n-1},n\in\N$. In what follows, for notational simplicity, we will drop the dependence on the expression $A_{n}^{k}[A_{n-1}^{j_{n-1}}\diamond A_{n-1}^{j_{n-1}+1}]\ldots[A_{1}^{j_1}\diamond  A_{1}^{j_1+1}]$.

The  main challenge in proving our claim is to create ``spaces" ({\it i.e.} separation) among trees inside our family $\p_n$. But for this, we will need first to create the tree-structures. Thus, our first step is to 'stick' every tile $P\in\p_n$ to a top (maximal tile with respect of  $``\leq"$). For this, we will proceed as follows:

Let $\bar{\p}_{n}^{max}=\left\{\bar{P}_{j}\right\}_j$. Proceeding as in \cite{q}, we define
 \beq \label{triang}
 \bar{\p}_n:=\left\{P\in\p_n\:|\:\:\exists\:j\in \N \:s.t.\:\:\:\:\:4P\triangleleft \bar{P}_{j}\right\}
 \eeq
and further define the set
\beq \label{cn}
 \C_n:=\left\{P\in\p_n\:|\:\operatorname{there\:are\: no\:chains}\:P\lneq P_{1}\lneq\ldots\lneq
P_{n}\:\&\:\left\{P_j\right\}_{j=1}^{n}\subseteq\p_n\:\right\}\:.
\eeq
With this done, we claim that
\beq \label{contt}
\p_n\setminus\C_n\subseteq\bar{\p}_n\:.
\eeq

Indeed, assume that $P\in \p_n\setminus\C_n$. Then from \eqref{cn} we have that $\exists\:\left\{P_j\right\}_{j=1}^{n}\subseteq\p_n$ such that
\beq \label{cn1}
P\lneq P_{1}\lneq\ldots\lneq P_{n}\:.
\eeq
Then, since $P_n\in \p_n$ we must have
$$A_{\P, A_{n}^{k}[A_{n-1}^{j_{n-1}}\diamond A_{n-1}^{j_{n-1}+1}]\ldots[A_{1}^{j_1}\diamond  A_{1}^{j_1+1}]}(P)\in (2^{-n},\,2^{-n+1}]$$ and hence, from the definition of $\bar{\p}_{n}^{max}$ and Definition \ref{mass} we have
\beq \label{mascx}
\exists\:\bar{P}\in\bar{\p}_{n}^{max}\:\:\textrm{s.t.}\:\: \Delta(10P_n,\:10\bar{P})< 2^{\frac{n}{N}}\:,
\eeq
which implies
\beq \label{mascx1}
\sup_{{q_{n}\in P_n}\atop{\bar{q}\in \bar{P}}}\|q_n-\bar{q}\|_{L^{\infty}(\tilde{I}_{P_n})}\leq 20\,(2d)^d\,|I_{P_n}|^{-1}\, 2^{\frac{n}{N}}\:.
\eeq
This last relation together with \eqref{cn1} and \eqref{p12} gives us
\beq \label{mascx2}
\sup_{{q_{1}\in P_1}\atop{\bar{q}\in \bar{P}}}\|q_1-\bar{q}\|_{L^{\infty}(\tilde{I}_{P_1})}\leq 50\,(2d)^d\,|I_{P_1}|^{-1}\:.
\eeq
Appealing now to \eqref{p11} and making in an essential way use of the second item in Observation \ref{redu} we have that
\beq \label{mascx3}
\exists\:q\in P\:\:\textrm{s.t.}\:\:\sup_{\bar{q}\in \bar{P}}\|q-\bar{q}\|_{L^{\infty}(\tilde{I}_{P_1})}\leq 50\,(2d)^d\,|I_{P_1}|^{-1}<|I_{P}|^{-1}\:,
\eeq
which now implies that $4P\triangleleft \bar{P}$ thus proving \eqref{contt}.

Now, defining the set $\D_n\subseteq \C_n$ with the property $\p_n\setminus\D_n=\bar{\p}_n$,
we remark that $\D_n$ breaks up as a disjoint union of a most $n$
sets $\D_{n}^{1}\cup\D_{n}^{2}\cup\ldots\cup\D_{n}^{n}$ with each
$\D_{n}^{j}$ being - recall Definition \ref{incomp} - an incomparable family of tiles. As a consequence, $\D_n$ may be written as a union of at most $n$ sparse $L^{\infty}-$forests
and hence, assuming that Proposition \ref{prop1} holds\footnote{Notice that Proposition \ref{prop1} is just a very particular case of Proposition \ref{prop2}.}, we can erase this set from $\p_n$ without affecting our claim.

Thus, in what follows, it will be enough to limit ourselves to the set of tiles $\bar{\p}_n$
which for convenience we will re-denote it with $\p_n$.

Returning to our Claim \ref{CL}, our aim is to show that
\beq \label{pnai}
\p_n=\bigcup_{j=1}^{c\,n}\S_{nj}\:,
\eeq
with each $\S_{nj}$ an $L^{\infty}$-forest of generation $n$.

Set now
\beq \label{bps}
B(P):=\#\left\{j\:|\:4P\trianglelefteq \bar{P}_j\right\}\:\:\:\:\:\:\:\forall\:\:P\in\p_n\:.
\eeq

Notice now that based on \eqref{contc}, \eqref{triang} and \eqref{bps} we have that
\beq \label{sppn}
\p_n=\bigcup_{j=1}^{c\,n} \p_{nj}\,,
\eeq
with
\beq \label{pjsn}
\p_{nj}:=\left\{P\in\p_n\:|\:2^j\leq B(P)< 2^{j+1}\right\}\:\:\:\:\:\:\:\:\:\:\forall\:j\in \left\{0,..c\,n\right\}\:.
\eeq

In what follows, we will show that each set $\p_{nj}$ can be written as
\beq \label{decp}
\p_{nj}=\S_{nj}\cup\r_{nj}\,,
\eeq
such that
\begin{itemize}
\item $\S_{nj}$ is an $L^{\infty}$-forest of generation $n$;

\item $\r_{nj}$ is a \textit{negligible} collection of tiles.
\end{itemize}

Fix now a family $\p_{nj}$.
$\newline$

\noindent \textbf{Step 1.} \textsf{Identifying the candidates for the tops of
the future trees.}
$\newline$

For this, we define
\beq \label{maxc}
\p_{nj}^{max}:=\left\{P^{r}=[\vec{\a}_r,\,I_r]\right\}_{r\in\left\{1,\ldots ,s\right\}}\subseteq\p_{nj}
\eeq
be the set of tiles with the property that
\beq \label{maxc1}
4P^{r}\: \textrm{is maximal with respect to $\leq$ inside the set}\:4\p_{nj}\,.
\eeq
Now, in many of the further reasonings we will use
the following
$\newline$

\noindent\textbf{Four key properties}
\begin{enumerate}[(A)]
    \item $\:\:\:\: 4P^{l}\leq 4P^{m}\:\:\Rightarrow\:\:I_l=I_m\:; \\$
    \item $\:\:\:\: \forall\:P\in\p_{nj}\:\:\:\exists\:\:\:P^{l}\:\:\:\textrm{s.t.}\:\:\:\:12P\trianglelefteq 4P^{l}\:;\\$
  \item  $\:\:\:\:\operatorname{If}\:P\in\p_{nj}\:\:\operatorname{s.t.}\:\:\exists\:\:m\not=l\:\:\:\operatorname{with}\:\:\:\left\{{4P\trianglelefteq\:4P^{l}}\atop{4P\trianglelefteq\: 4P^{m}} \right.\:\:\operatorname{then}\:\:\:\left\{{4P^{m}\leq\:4P^{l}}\atop{4P^{l}\leq\: 4P^{m}} \right.\;;\\$
  \item  $\:\:\:\:\textrm{If}\: P_j=[\vec{\a},I_j]\in\P\:\:\textrm{with}
\:\:j\in\left\{1,2\right\}\: \textrm{s.t.}\:
|I_1|\not=|I_2|,\,\operatorname{then}\\\:\:|I_1|\leq 2^{-D}\:|I_2|
\:\textrm{or}\: |I_2|\leq 2^{-D}\:|I_1|\:.\\$
\end{enumerate}

\noindent\textbf{The four properties - explanations}

\begin{enumerate}[(A)]
    \item this is an immediate consequence of \eqref{maxc1} and Definition \ref{ordmaxrel};
    \item from \eqref{maxc1} we have that for any $P\in\p_{nj}$ there exists $P^l$ such that $4P\leq 4P^{l}$; now $(B)$ is a consequence of Observation \ref{ordmaxrel1};
    \item this follows from a contrapositive reasoning: if $4P^{l}$ and $4 P^{m}$ are incomparable, then, using the fact that $``\trianglelefteq"$ is an order relation (see Observation \ref{ordmaxrel1}) we deduce that $B(P)\geq 2^{j+1}$ contradicting thus the fact that $P\in\p_{nj}$;
    \item this is simply restating \eqref{separat} - see Observation \ref{redu}.
\end{enumerate}

$\newline$
\noindent \textbf{Step 2.} \textsf{Isolating the negligible family of tiles $\r_{nj}$.}
$\newline$

Our aim here is to properly trim the set $\p_{nj}$ so that the resulting family will have all the desired properties of an $L^{\infty}$ forest of generation $n$.

In order to do so, we define the following three sets:
\begin{itemize}
\item $\r_{nj}^{1}$ - the family of tiles that are ``far away" from $\p_{nj}^{max}$:
\beq \label{r1}
\r_{nj}^{1}:=\left\{P\in\p_{nj}\:|\:
\forall\:\:P^{l}\:\:\Rightarrow\:\:\frac{3}{2}P\nleqslant
P^{l}\right\}\,.
\eeq
\item $\r_{nj}^{2}$ - the family of maximal tiles and of their neighbors:
\beq \label{r2}
\r_{nj}^{2}:=\left\{P\in\p_{nj}\:|\:\exists\:l\:st\:|I_P|=|I_{P^l}|\:,\:\frac{3}{2}P\leq P^{l}\right\}\,.
\eeq
\item $\r_{nj}^{3}$ - the family of minimal tiles:
\beq \label{r3}
\r_{nj}^3:=\left\{P\in\p_{nj}\:\big|\:P\:
\textrm{minimal}\right\}\:.
\eeq
\end{itemize}

With this, we define:
\beq \label{rest}
\r_{nj}:=\r_{nj}^{1}\cup\r_{nj}^{2}\cup\r_{nj}^{3}\,.
\eeq

\begin{claim}\label{CL1}
The set $\r_{nj}$ is a negligible family of tiles.
\end{claim}

Indeed, let us justify our claim as follows:

\begin{itemize}
\item \textbf{for the set $\r_{nj}^{1}$} - we proceed by contradiction: assume that there exist  $P_1,\,P_2\in \r_{nj}^{1}$ such that $P_1\lneq P_2$. Applying now $(B)$ we have that there exist $P^{l_1},\,P^{l_2}\in \p_{nj}^{max}$ such that $12P_j\trianglelefteq 4P^{l_j}$ with $j\in\{1,2\}$. Using now $(D)$ we must have $\frac{3}{2}P_1 \trianglelefteq 12 P_2  \trianglelefteq 4P^{l_2}$ contradicting thus the assumption that $P_1\in \r_{nj}^{1}$.

\item \textbf{for the set $\r_{nj}^{2}$} - again we proceed by contradiction: assume that there exist  $P_1,\,P_2\in \r_{nj}^{2}$ such that $P_1\lneq P_2$ and hence there exist $P^{l_1},\,P^{l_2}\in \p_{nj}^{max}$ such that $|I_{P_1}|=|I_{P^{l_1}}|< |I_{P^{l_2}}|=|I_{P_2}|$ and $\frac{3}{2}P_j\leq P^{l_j}$. Applying $(D)$ we have that $4 P^{l_1}\leq 4P^{l_2}$ contradicting the maximality assumption.

\item \textbf{for the set $\r_{nj}^{3}$}: this family is negligible from the definition of what means a minimal family of tiles - see Observation \ref{ordmaxrel}.
\end{itemize}

$\newline$
\noindent \textbf{Step 3.} \textsf{Verifying that the set $\S_{nj}:=\p_{nj}\setminus\r_{nj}$ is an $L^{\infty}-$forest.}
$\newline$

For the remaining set $\S_{nj}$, we proceed as follows:

\begin{enumerate}
\item  Set $S_m=\left\{P\in\S_{nj}\:|\:\frac{3}{2}P\leq P^{m}\right\}$. In what follows we only consider those sets $S_m$ which are non-empty. Without loss of generality
we may suppose $\S_{nj}=\bigcup_{m=1}^{s}S_m$;
\item  Introduce the ``clustering" relation among the sets $\left\{S_m\right\}_m$:
$$S_m\propto S_l$$
iff $\exists\: P_1\in S_m$ and
$\exists\: P_2\in S_l$ such that $2P_1\leq 2P^l$ or $2P_2\leq 2P^m$;
\item Define a second relation on $\left\{S_m\right\}_m$ given by:
$$S_m\backsim S_l$$
if $4P^m\leq4P^l$ or equivalently $4P^l\leq4P^m$.
\item  Deduce that $S_m\propto S_l$ implies $S_m\backsim S_l$ and making use of property (C) conclude
that $``\backsim"$ is an equivalence relation.
\item  Let $\hat{m}:=\left\{l\:|\:S_l\backsim S_m \right\}$; then
the cardinality of $\hat{m}$ is at most $c(d)$, and for
$$\hat{S}_m:=\bigcup_{m'\in\hat{m}}S_{m'}\:,$$
one has that $\hat{S}_m$ is a tree having as a top any $P^l$ with $l\in\hat{m}$.
\end{enumerate}

Let us justify (1)-(5). Relations $(1), (2)$ and $(3)$ are simply definitions and thus we only need to verify $(4)$ and $(5)$.

We start with item $(4)$. Assume that $S_m\propto S_l$ with $m\not=l$ and thus we have wlog that $\exists\: P_1\in S_m$ such that $2P_1\leq 2P^l$. Notice first that $|I_{P_1}|<|I_{P^l}|$ as otherwise we must have $P_1\in \r_{nj}^{2}$  which is not allowed. Thus $4P_1\trianglelefteq\:4P^{l}$ and $4P_1\trianglelefteq\:4P^{m}$ and hence from $(C)$ we conclude that $S_m\backsim S_l$.
Next we need to show that $``\backsim"$ is an equivalence relation. The only nontrivial part is to check transitivity. Assume thus that $S_m\backsim S_l$ and $S_l\backsim S_r$. Since $S_{l}\not=\emptyset$ we have that there exists $P$ with $|I_P|<|I_{P^l}|$ and $\frac{3}{2}P\leq P^{l}$. Since from $(A)$ we must have $|I_{P^m}|=|I_{P^l}|=|I_{P^r}|$ deduce that
$4P\trianglelefteq\:4P^{m}$ and $4P\trianglelefteq\:4P^{r}$. Thus, from $(C)$ we further have $4P^{m}\leq 4P^{r}$  which implies $S_m\backsim S_r$.

We pass now to proving item $(5)$. The fact that any equivalence class $\hat{m}$ has at most $c(d)$ elements is a direct consequence of item $(3)$. We will now focus on proving that $\hat{S}_m$ is a tree by verifying all the three items in Definition \ref{tree}:
\begin{itemize}
\item $\forall\:\:P\in\hat{S}_m\:\:\:\Rightarrow\:\:\:\:\frac{3}{2}P\leq
10 P_0$ for any $P_0\in\{P^l\}_{l\in\hat{m}}$.

This is a direct consequence of the fact that $\forall\:\:P\in\hat{S}_m$ there exists $P^l$ with $l\in\hat{m}$ such that
$\frac{3}{2} P\leq P^l$. However for any other $P^r$ with $r\in\hat{m}$ we have that $4P^r \leq 4 P^l$ and $4P^l \leq 4 P^r$ and hence $\frac{3}{2} P\leq 10 P^r$.

\item if $P\in\hat{S}_m$ and $P'\in N(P)$ such
that $\frac{4}{3}P'\leq 10P_{0}$ then $P'\in \hat{S}_m$.

Let us take $P\in\hat{S}_m$. Our first task is to first show a milder fact: any $P'\in N(P)$ is available in the family $\S_{nj}$. Recall now two important facts: $\S_{nj}\subset\p_{nj}$ with $\p_{nj}$ defined in \eqref{pjsn} and $\S_{nj}$ does not contain the minimal elements in $\p_{nj}$ since in particular $\S_{nj}\subseteq\p_{nj}\setminus\r_{nj}^3$. Now since $P\in\hat{S}_m$ then
\beq \label{condx}
\eeq
\begin{itemize}
\item $\exists\:l\in\hat{m}$ such that $\frac{3}{2}P\lneq P^l$;

\item $\exists\:P^{min}\in \p_{nj}\setminus\S_{nj}$ and a chain $\{P_i\}_{i=1}^{M}\subseteq\p_{nj}$ with $M\in\N$, $M\geq 2$ such that $P_M\lneq P_{M-1}\ldots\lneq P_{1}$ and $P_M=P^{min}$ and $P_1=P$.
\end{itemize}
This immediately implies the following key relation:
\beq \label{condx1}
4 P^{min} \vartriangleleft4 P' \vartriangleleft 4 P^l\:.
\eeq
From Definition \ref{mass} and relation \eqref{pjsn} we conclude that
$$\forall\:P'\in N(P)\:\:\Rightarrow\:\:P'\in\p_{nj}\;.$$

Now \eqref{r3} immediately implies that if $P'\in N(P)$ such that $\frac{4}{3}P'\lneq 10P^l$ then  $\frac{3}{2}P'\lneq P^l$ and hence $P'\in \S_{nj}$ and moreover $P'\in \hat{S}_m$.

\item if $P_1,\:P_2\in\hat{S}_m$ and $P_1\leq P \leq P_2$ then
$P\in\hat{S}_m$

We may assume wlog that $P_1\lneq P \lneq P_2$. Notice that in this case we immediately have that
$\frac{3}{2} P_1 \vartriangleleft \frac{3}{2} P \vartriangleleft\frac{3}{2} P_2$ which immediately implies both
that $P\in\p_{nj}$ and that $\frac{3}{2} P \vartriangleleft P^m$ and hence $P\in \hat{S}_m$.
\end{itemize}
$\newline$

This proves our item $(5)$ saying that $\hat{S}_m$ is a tree with top $P_0\in\{P^l\}_{l\in\hat{m}}$.

Finally, since for any two distinct $\hat{S}_m$ and  $\hat{S}_l$ we have that taking any correspondent $S_{m'}\in \hat{S}_m$ and $S_{l'}\in \hat{S}_l$  the relation $S_{m'}\propto S_{l'}$ does \textit{not} hold
we conclude that the set
\beq \label{for}
\S_{nj}=\bigcup_{m}\hat{S}_m\,,
\eeq
is an $L^{\infty}$-forest as in Definition \ref{infforest}.

\section{\bf The proofs of Propositions 1 and 2}\label{prop12}

In this section we will analyze the $L^p$ boundedness behavior of the operator associated with a generic (sparse) forest.

\begin{obs}\label{envir}
In what follows $\p$ is a generic (sparse) forest that should be thought of as
\beq\label{snjun}
\bigcup_{k,j\in\N}\bigcup_{j_1,\ldots j_{n-1}\in\N} \S_{nj}[A_{n}^{k}[A_{n-1}^{j_{n-1}}\diamond A_{n-1}^{j_{n-1}+1}]\ldots[A_{1}^{j_1}\diamond  A_{1}^{j_1+1}]],\,
\eeq
with each
$$\S_{nj}[A_{n}^{k}[A_{n-1}^{j_{n-1}}\diamond A_{n-1}^{j_{n-1}+1}]\ldots[A_{1}^{j_1}\diamond  A_{1}^{j_1+1}]]\subset \p_{n}[A_{n}^{k}[A_{n-1}^{j_{n-1}}\diamond A_{n-1}^{j_{n-1}+1}]\ldots[A_{1}^{j_1}\diamond  A_{1}^{j_1+1}]]$$ constructed in a similar fashion as $S_{nj}$ in \eqref{for}.
\end{obs}

\subsection{\bf Proof of Proposition 1}
$\newline$
We begin by restating the result that we need to prove:
$\newline$

\noindent\textbf{Proposition 1.} \textit{Let $\p\subseteq\P_n$ be a sparse forest. Then there exists $\eta\in(0,1/2)$, depending only on the degree $d$,
such that for $1<p<\infty$ we have
$$\left\|T^{\p}\right\|_{p}\lesssim_{p,d} 2^{-n\,\eta\,(1-\frac{1}{p^*})}\:.$$}

\subsubsection{\bf The $L^2$ bound}
$\newline$

Throughout the remaining reasonings we will assume wlog that
$$n\geq c(d)\geq d^d\,.$$
Assume $P=[\vec{\a},I]$ and $P'=[\vec{\a}',I']$ with $|I|\leq |I'|$. Inspecting now the proof of \eqref{v17} in Lemma \ref{interact} we realize that the following holds
\beq\label{estimlm0}
\left|T_{P^\prime}T^{*}_{P}f(x)\right|\lesssim{\left\lceil {\Delta}(P,P^\prime)\right\rceil}^{1/d}
\frac{\int_{E(P)}\left|f\right|}{\left|I^\prime\right|}\chi_{E(P^\prime)}(x)\:.
\eeq
Now, proceeding as in the corresponding proof of Proposition 1 in \cite{q}, we have
$$\int_{\TT}\left|\left ({T^{\p}}\right )^{*}f(x)\right|^2dx\lesssim\left|\sum_{P^{\prime}\in\p\atop{P^\prime=[\vec{\a}^\prime,I^\prime]}}
\int_{\TT}f(x)\left\{\sum_{P=[\vec{\a},I]\in\:\p\atop{\left|I\right|\leq\left|I^\prime\right|}}
\overline{T_{P^\prime }T^{*}_{P}f}(x)\right\}dx\:\right|$$
$$\lesssim\sum_{P^\prime\in\p}\int_{E(P^\prime)}|f|\left\{\sum_{P\in a(P^\prime)}{\left\lceil {\Delta}(P,P^\prime)\right\rceil}^{1/d}\frac{\int_{E(P)}\left|f\right|}{\left|I^\prime\right|}\right\}$$
$$+\:\sum_{P^\prime\in\p}\int_{E(P^\prime)}|f|\left\{\sum_{P\in b(P^\prime)}{\left\lceil {\Delta}(P,P^\prime)\right\rceil}^{1/d}\frac{\int_{E(P)}\left|f\right|}{\left|I^\prime\right|}\right\}=^{def}
A\:+\:B$$ where here we have used the following notations:
$$a(P^\prime)=\left\{P=[\vec{\a},I]\in\p\,,\:\:|I|\leq |I^\prime|\:\&\:I^{*}\cap {I'}^{*}\not=\emptyset\:\:|\:\:\Delta(P,P^\prime)\leq 2^{n\ep} \right\}\:,$$
$$b(P^\prime)=\left\{P=[\vec{\a},I]\in\p\,,\:\:|I|\leq |I^\prime|\:\&\:I^{*}\cap {I'}^{*}\not=\emptyset\:\:|\:\:\Delta(P,P^\prime)\geq 2^{n\ep} \right\}\:.$$
with $\ep\in(0,1)$ small enough (\textit{e.g.} $\ep=\frac{1}{100\,(N+d)}$ with $N$ defined in \eqref{v1}).

Further, we have
$$ A\lesssim\sum_{P^\prime\in\p}\int_{E(P^\prime)}\left|f(x)\right|\left\{\frac{1}{\left|I^\prime\right|}\sum_{P\in a(P^\prime)}\int_{E(P)}\left|f\right|\right\}dx=\int |f|\, V_a(|f|)\:, $$
where by definition
\beq\label{defva}
V_a(f):=\sum_{P'=[\vec{\a}',I']\in\p}\frac{\chi_{E(P')}}{|I'|}\,\sum_{P\in a(P^\prime)}\int_{E(P)}f\:.
\eeq
Similarly, using the definition of $b(P^\prime)$ we deduce
$$ B\lesssim\sum_{P^\prime=[\vec{\a}',I']\in\p}\int_{E(P^\prime)}\left|f(x)\right|\left\{\frac{2^{-n\,\frac{\ep}{d}}}{\left|I^\prime\right|}
\sum_{P\in b(P^\prime)}\int_{E(P)}\left|f\right|\right\}dx=2^{-n\,\frac{\ep}{d}}\,\int |f|\, V_b(|f|)\:, $$
where by definition
\beq\label{defvb}
V_b(f):=\sum_{P'=[\vec{\a}',I']\in\p}\frac{\chi_{E(P')}}{|I'|}\,\sum_{P\in b(P^\prime)}\int_{E(P)}f\:.
\eeq

We will now focuss on providing $L^2$-bounds on $V_a(f)$.

Fix $1<r<2$ and let $r'$ be the H\"older conjugate of $r$. Suppose wlog that $f\geq0$. Then
$$V_a(f)\leq\sum_{P'=[\vec{\a}',I']\in\p}\chi_{E(P')}\left(\frac{\int_{\tilde{I'}}f^r}{|I'|}\right)
^{\frac{1}{r}}\, \frac{\| \sum_{P\in a(P^\prime)}\chi_{E(P)} \|_{r'}}{|I'|^{\frac{1}{r'}}}\:.$$

The first key observation derived from the structure of the set $\p$ and the definition of $a(P)$ is

\begin{claim}\label{Cm}
The following Carleson measure type condition holds:
\beq\label{ap}
\| \sum_{P\in a(P')}\chi_{E(P)} \|_{r'}\lesssim_{r} 2^{-\frac{n}{r'}(1-10d\ep-10 N \ep)}\,|I'|^{\frac{1}{r'}}\:.
\eeq
Here $N$ stands for the parameter used for the mass definition in \eqref{v1}.
\end{claim}

As a consequence of \eqref{for2} in Definition \ref{forest}, it is enough to show \eqref{ap} for $\p$ a sparse $L^{\infty}$-forest.
\medskip

\noindent\textbf{Step 1} If  $\r\subset\p$ is such that $\r$ is a collection of {\it incomparable} tiles then the restriction of \eqref{ap} to $\r$ holds, that is
\beq\label{apinc}
\| \sum_{{P\in a(P')}\atop{P\in\r}}\chi_{E(P)} \|_{r'}\lesssim_{r} 2^{-\frac{n}{r'}(1-5d\ep-5 N \ep)}\,|I'|^{\frac{1}{r'}}\:.
\eeq
$\newline$

In order to show \eqref{apinc}, we first claim that
\beq\label{apinco}
\| \sum_{{P\in a(P')}\atop{P\in\r}}\chi_{E(P)} \|_{1}\lesssim 2^{-n(1-5d\ep-5 N \ep)}\,|I'|\:.
\eeq
Indeed, to see this we define
$$\I_{\r}(P'):=\{I\,|\,\exists\:P=[\vec{\a},I]\in a(P')\cap \r\}\,,$$
and assume wlog that $\I_{\r}(P')\not=\emptyset$.

Further, let $\I_{\r,min}(P')$ be the set of minimal intervals relative to inclusion belonging to the set $\I_{\r}(P')$.

Let
$$\breve{\I}_{\r}(P'):=\left\{I\subset 50 I'\,\big|\, \begin{array}{cc}
\textrm{Only one of the left or right children}\\
\textrm{of $I$ contains an element of} \:\I_{\r,min}(P')
\end{array}\right\}\,\cup\,\I_{\r,min}(P')\;.$$
Also, set
\beq\label{basis}
Basis(\p):=\bigcup_{P=[\vec{a},I]\in\p} I\:,
\eeq
and let
\beq\label{basisa}
\breve{a}(P'):=\left\{\begin{array}{cc} P=[\vec{a},I]\in\P\\I\subset Basis(\p)\end{array}\,\big|\,I\in \breve{\I}_{\r}(P')\:\:
\textrm{and}\:\:\Delta(P,P^\prime)\leq 2^{n\ep}\right\}\;.
\eeq

Using now the smoothening effect encoded in the mass definition \eqref{v1}, we deduce that
\beq\label{L1inc}
\eeq
$$\| \sum_{{P\in a(P')}\atop{P\in\r}}\chi_{E(P)} \|_{1}\leq \| \sum_{P\in \breve{a}(P')}\chi_{E(P)} \|_{1} $$
$$\lesssim 2^{-n}\,2^{5\ep n N}\,2^{5\ep n d} \sum _{I\in \breve{\I}_{\r}(P')}|I|\lesssim 2^{-n (1-5\ep N- 5\ep d)}\,|I'|\:.$$

The $L^{\infty}$ bound follows trivially since $\r$ is an incomparable family of tiles:
\beq\label{Linfinc}
 \| \sum_{P\in\r}\chi_{E(P)} \|_{\infty}\leq 1\:.
\eeq

By interpolating between \eqref{L1inc} and  \eqref{Linfinc} we deduce that \eqref{apinc} holds.

$\newline$
\noindent\textbf{Step 2} The general sparse forest case.
$\newline$

By Definition \ref{infforest}, we have that $\p\cap  a(P')=\bigcup_{j}\p_j$ with $\{\p_j\}_j$ sparse separated trees.

 Further set $\textrm{top}\: \p_j= P_j$ and let
$$\p_j^{1}=\{P\in\p_j\,|\,\textrm{there is no chain}\:\:P<P^{1}<\ldots<P^{n}=P_j\:\:\textrm{s.t.}\:\:P^{k}\in\p_j\}\:,$$
and $$\p_j^{2}:=\p_j\setminus\p_j^{1}\:.$$
In the above setting, by appealing to maximal chain decompositions, we notice that $\p\cap  a(P')$ can be written as
\beq\label{decsp}
\left(\bigcup_{l=1}^{n} \A_l\right)\cup\left(\bigcup_{j}\p_j^{2}\right)\:,
\eeq
such that
\begin{itemize}
\item each $\A_l$ is a set of incomparable tiles;

\item the second component verifies
\beq\label{disj}
\sum_{j\,,\:\p_j^{2}\not=\emptyset}\:\chi_{I_{P_j}}\leq 1\,.
\eeq
\end{itemize}
Indeed to see this we notice that if $P_i$ and $P_j$ are the tops of two separated trees such that $\p_i^{2},\,\p_j^{2}\not=\emptyset $ then either $I_{P_i}\cap I_{P_j}=\emptyset$ or we must have ${\Delta}(P_i,P_j)\gtrsim 2^{n}\max\{|\o_{P_i}|,\,|\o_{P_j}|\}$. However, only the first scenario is possible since condition $P_i,\,P_j\in a(P')$
requires ${\Delta}(P_i,P_j)\lesssim_{d} 2^{\ep n}\max\{|\o_{P_i}|,\,|\o_{P_j}|\}$.

 Finally, from Step 1, we know that \eqref{apinc} holds for each $\r=\A_l$, while from the fact that each $\p_j$ is a sparse tree we deduce that
 \beq\label{cmc1}
 \| \sum_{P\in \p_j^{2}} \chi_{E(P)}\|_{r'}\lesssim_{r} 2^{-n\,\frac{1}{r'}}\,|I_{P_j}|^{\frac{1}{r'}}\:.
 \eeq

Combining now Step 1 with  \eqref{disj} and \eqref{cmc1}, we conclude that Claim \ref{Cm} is true.

Now, in order to control the term $A$ it remains to show the following

\begin{claim}\label{contBMO}
With the previous notations, defining
\beq\label{vv}
\V f:=\sum_{P=[\vec{\a},I]\in\p}\chi_{E(P)}\left(\frac{\int_{\tilde{I}}f^r}{|I|}\right)^{\frac{1}{r}}\:,
\eeq
we have
\beq\label{v}
\|\V f\|_2\lesssim_{r} \|\sum_{P\in\p}\chi_{E(P)}\|_{BMO_C} \|f\|_2\lesssim \|f\|_2\:.
\eeq
\end{claim}

Set now $\I:=\{I\,|\,\exists\:P=[\vec{\a},I]\in\p\}$ and $E(I):=\bigcup_{{P=[\vec{\a},I_P]\in\p}\atop{I_P=I}} E(P)$.
Rewrite $\V$ as follows:
$$\V f=\sum_{I\in\I}\chi_{E(I)}\left(\frac{\int_{\tilde{I}}f^r}{|I|}\right)^{\frac{1}{r}}\:.$$
Denote with $\I_m:=\{I\in\I\,|\,\frac{\int_{\tilde{I}}f^r}{|I|}\approx 2^{m}\}$ and notice that $\I=\bigcup_{m\in\Z} \I_m$.
Also denote with $\I_m^{max}$ the set of maximal intervals (with respect of inclusion) in $\I_m$.
Assume wlog that $\int_{\TT} f^{r} \approx 2^{m_{0}}$ for some $m_{0}\in \Z$.
Now, for each $m\geq m_{0}$, notice then that $\I_m^{max}$ consists of pairwise disjoint intervals.

Then we have
$$\V f=\sum_{m\in\Z}\sum_{I\in\I_m}\chi_{E(I)}\left(\frac{\int_{\tilde{I}}f^r}{|I|}\right)^{\frac{1}{r}}
\lesssim (\int_{\TT} f^{r})^{\frac{1}{r}}\,+ \,\sum_{m\geq m_{0}}\sum_{J\in\I_m^{max}}\sum_{{I\subseteq J}\atop{I\in\I_m}}2^{\frac{m}{r}}\,\chi_{E(I)}\:,$$
and thus, ignoring the $L^r$ norm of $f$, one has
$$\|\V f\|_2^2\approx \sum_{m,m'\geq m_{0}}\sum_{{J\in\I_m^{max}}\atop{J'\in\I_{m'}^{max}}}2^{\frac{m+m'}{r}}
\int (\sum_{{I\subseteq J}\atop{I\in\I_m}}\chi_{E(I)})(\sum_{{I'\subseteq J'}\atop{I'\in\I_{m'}}}\chi_{E(I')})$$
$$\approx\sum_{m\geq m_{0}}\sum_{m'\geq m}\sum_{{J\in\I_m^{max}}}\sum_{{J'\subseteq J}\atop{{J'\in\I_{m'}^{max}}}}2^{\frac{m+m'}{r}}
\int (\sum_{{I\subseteq J}\atop{I\in\I_m}}\chi_{E(I)})(\sum_{{I'\subseteq J'}\atop{I'\in\I_{m'}}}\chi_{E(I')})$$
Let $1\leq q<\infty$ and $J\subseteq [0,1]$ fixed.

Applying John-Nirenberg to
\beq\label{Carlbmo}
\|\sum_{I\in\I}\chi_{E(I)}\|_{BMO_D}\lesssim 1\,,
\eeq
we deduce the Carleson packing condition
\beq\label{Carl}
\| \sum_{{I\subseteq J}\atop{I,\,J\in\I}}\chi_{E(I)} \|_{q}^{q}\lesssim_{q} |J|\:.
\eeq
 Now, from \eqref{Carl} and Cauchy-Schwarz, for $1<p<r<2$, we further have
$$\|\V f\|_2^2\lesssim \sum_{m}\sum_{m'\geq m}2^{\frac{m+m'}{r}}\sum_{{J\in\I_m^{max}}}\| \sum_{I\subseteq J}\chi_{E(I)} \|_{p'}
\|\sum_{{J'\subseteq J}\atop{{J'\in\I_{m'}^{max}}}}\sum_{I'\subseteq J'}\chi_{E(I')} \|_{p}$$
$$\lesssim \sum_{m}\sum_{{J\in\I_m^{max}}}2^{\frac{m}{r}}\,|J|^{\frac{1}{p'}}\sum_{m'\geq m}2^{\frac{m'}{r}}\,
(\sum_{{J'\subseteq J}\atop{{J'\in\I_{m'}^{max}}}}|J'|)^{\frac{1}{p}}$$
$$\lesssim \sum_{m}\sum_{{J\in\I_m^{max}}}2^{\frac{m}{r}}\,|J|^{\frac{1}{p'}}
\sum_{m'\geq m}2^{\frac{m'}{r}}\,2^{-\frac{m'}{p}}\,(\int_{\tilde{J}}f^r)^{\frac{1}{p}}$$
$$\lesssim \sum_{m}\sum_{{J\in\I_m^{max}}}2^{\frac{2m}{r}}\,|J|\lesssim
\sum_{m}2^{\frac{2m}{r}}\,2^{-m}\,\int_{(M_r f)^r \gtrsim 2^m}(M_r f)^r$$
$$\lesssim \int (M_r f)^2\lesssim_{r}\int f^2\:,$$
where here we denoted $M_r f(x):=\big(\sup_{x\in I}\frac{\int_{I^{*}} |f|^r}{|I|}\big)^{\frac{1}{r}}\:.$

This ends the proof of Claim \ref{vv}.

Thus, combining now  \eqref{ap} and \eqref{v}, for an appropriate choice of $\ep$, we conclude that
$$A \lesssim_{r} 2^{-\frac{n}{r'}(1-10 d\ep-10 N\ep)}\,\|f\|_2^2\lesssim 2^{-\frac{n}{2 r'}}\,\|f\|_2^2 \:.$$

The $B$ term can be similarly treated if one replaces \eqref{ap} with just
\beq\label{bp}
\| \sum_{P\in b(P')}\chi_{E(P)} \|_{r'}\lesssim_{r} \,|I'|^{\frac{1}{r'}}\:,
\eeq
thus obtaining
$$B \lesssim 2^{-n\,\frac{\ep}{d}}\,\|f\|_2^2\:.$$

Now, properly choosing $r$ and $\ep$, we conclude that there exists $\eta=\eta(d)\in (0,1)$ such that
\beq\label{Top2}
\left\|T^{\p} f\right\|_{2}\lesssim 2^{-n\,\frac{\eta}{2}}\,\|f\|_2\:.
\eeq

This ends our proof.

\subsubsection{\bf The $L^p$ bound}

Suppose first that $2\leq p<\infty$. For any $f\in L^1(\TT)$ we define the operator
\beq\label{Lop}
L^{\p}f(x):=\sum_{P=[\vec{\a},I_P]\in\p}\frac{\int_{I_{P^{*}}}f}{|I_P|}\chi_{E(P)}\:.
\eeq

Now, on the one hand, repeating the reasonings from the case $p=2$, one has
$$\| L^{\p}f\|_2\lesssim\|f\|_2\;.$$

On the other hand
$$\left\|L^{\p}\right\|_{\infty\rightarrow BMO_D}\lesssim 1\:.$$

Interpolating now between $L^2\rightarrow L^2$ and $L^{\infty}\rightarrow BMO_D$ we obtain that
\beq\label{Lopbound}
\|L^{\p}\|_p\lesssim p\:.
\eeq

Consequently, based on the straightforward relation
$$|T^{\p}f|\lesssim L^{\p}|f|\:,$$
we also get that for any $2\leq p<\infty$ one has
\beq\label{Topbound}
\|T^{\p}\|_p\lesssim p\:.
\eeq
Interpolating now between \eqref{Top2} and \eqref{Topbound} one obtains the desired conclusion (possibly by changing the exponent $\eta$ with a small factor).

For the case $1<p<2$ we need to focus on the behavior of
${T^{\p}}^{*}$.

Indeed, on the one hand we know that
$$\left\|{T^{\p}}^{*}\right\|_{2\rightarrow 2} = \left\|T^{\p}\right\|_{2\rightarrow
2}\lesssim 2^{-n\,\frac{\eta}{2}}\:.$$

On the other hand, for $f\in L^{\infty}$ we have
$$\|{L^{\p}}^{*}f\|_{BMO_D}=\|\sum_{P=[\vec{\a},I_P]\in\p}\frac{\int_{E(P)}f}{|I_P|}\chi_{I_{P^*}}\|_{BMO_D}\lesssim
\left\|f\right\|_{\infty}\:.$$

Thus, as before, for any $2\leq q=p'<\infty$ one has
$$\|{T^{\p}}^{*}f\|_{q}\lesssim \|{L^{\p}}^{*}|f|\|_{q}\lesssim q\:.$$

The claim now follows by interpolation.\footnote{We use here the
fact that $\left\|{T^{\p}}^{*}\right\|_{p'\rightarrow
p'}=\left\|{T^{\p}}\right\|_{p\rightarrow p}$.}

\begin{flushright}
$\Box$
\end{flushright}

\subsection{\bf Preparatives for the proof of Proposition 2}
$\newline$

As the name suggest, this section is meant for ``preparing the ground"
for the proof of Proposition 2. Most of the results presented here, have
a direct analogue in either \cite{f} or \cite{q}, and thus, we will not
insist much on their proofs but only treat the sensitive points that are different.

\subsubsection{\bf $L^2-$results. Main Lemma}
$\newline$

We start with the following

\begin{l1}\label{tr2}[$L^2$-\textsf{uniform mass tree estimate}]

Fix $\delta\in(0,1]$ and let $\p\subseteq\P$ be a tree with spacial support $I_0$ such that
\beq\label{masscon}
A_{\P,I_0}(P)<\delta\:\:\:\:\forall\:\:\:P\in\p\:.
\eeq
Then
\beq\label{treee}
\left\|T^\p\right\|_2\lesssim_{d}\delta^{\frac{1}{2}}\:. \eeq
\end{l1}

\begin{proof}
For more details please see Lemmas \ref{trp} and \ref{treecutp} in the next section.
\end{proof}

\begin{d0}\label{sep}[\textsf{Separated trees}]

Fix a number $\d\in(0,1]$. Let $\p_1$ and $\p_2$ be two trees with
tops $P_1=[\vec{\a}_1,I_1]$ and respectively
$P_2=[\vec{\a}_2,I_2]$. We say that $\p_1$ and
$\p_2$ are $\d^{-1}$-\textbf{separated} if either $I_1\cap
I_2=\es$ or else
\begin{itemize}
\item  $P=[\vec{\a},I]\in\p_1\:\:\&\:\:I\subseteq
I_2\:\:\:\:\Rightarrow\:\:\:\left\lceil
\Delta(P,P_2)\right\rceil<\d\:,$
\item  $P=[\vec{\a},I]\in\p_2\:\:\&\:\:I\subseteq
I_1\:\:\:\:\Rightarrow\:\:\:\left\lceil
\Delta(P,P_1)\right\rceil<\d\:.$
\end{itemize}
\end{d0}

\begin{d0}\label{sepcr}[\textsf{Separation and critical sets}]

Fix $\d\in(0,1)$ small and $\ep_0\in(0,1)$. Let $\p_1$ and $\p_2$ be two $\d^{-1}$-separated trees as above. Also let $q_j$ be the central polynomial of
$P_j$, $j\in\{1,2\}$, and $q_{1,2}=q_1-q_2$ the $(P_1,P_2)-$interaction polynomial. Recalling the construction in Lemma \ref{CZdec}, we define

\begin{itemize}
\item $I[s]$ - the {\bf separation set} of $\p_1$ and $\p_2$ by
\beq\label{seps}
I[s]:=J_s\left(q_{1,2},\,c_0(d)\,\d^{-1}\right)\:,
\eeq
for $J:=\tilde{I}_1\cap\tilde{I}_2$ and $c_0(d)>0$ properly\footnote{See Observation \ref{sepcrst} below.}chosen.

\item $I[c]$ - the {\bf ($\ep_0$-)critical
intersection set}  by
\beq\label{crits}
I[c]:=J_s\left(q_{1,2},\,c_0(d)\,\d^{-\ep_0}\right)\:,
\eeq
for $J:=\tilde{I}_1\cap\tilde{I}_2$.
\end{itemize}
\end{d0}

\begin{obs}\label{sepcrst} It is important to notice the
following three properties of our above-defined sets; these
facilitate the adaptation of the reasonings involved in the proofs
of Lemmas \ref{sept} and \ref{rt} to those of the corresponding Lemmas 2 and 4 in \cite{q}:
\begin{enumerate}
\item in what follows we will choose $c_0(d)$ in \eqref{seps} such that if $\{I^j\}_{j=1}^l$ is the decomposition
of $I[s]$ analogue to \eqref{defjs} then for any $j\in\{1,\ldots , l\}$
\beq\label{ob1}
\forall\:\:P=[\vec{\a},I_P]\in\p_1\cup \p_2\:\:\:\textrm{if}
\:\:I^j\cap \tilde{I}_P\not=\emptyset\:\:\textrm{then}\:\:|I_P|>|I^j|\,.
\eeq
Deduce that in particular we must have that for any $j\in\{1,\ldots , l\}$
\beq\label{ob2}
\Delta_{q_{1,2}}(I^j)\gtrsim_{d} \d^{-1}\,.
\eeq

\item from Lemma \ref{CZdec} relations \eqref{parthigh}-\eqref{cont}, we further deduce that for any dyadic
 $J\subset \tilde{I}_1\cap\tilde{I}_2$ such that $I[s]\cap J=\emptyset$ we have (for $c(d)\leq d^d$)
\beq\label{ob3}
\inf_{x\in J}|q_{1,2}(x)|\leq \sup_{x\in J}|q_{1,2}(x)|\leq c(d)\inf_{x\in J}|q_{1,2}(x)|\:,
\eeq
and
\beq\label{ob4}
\Delta_{q_{1,2}}(J)\gtrsim_{d} \d^{-1}\,.
\eeq
Moreover, one has
\begin{quote}
$\:\:\forall\:P=[\vec{\a},I_P]\in\p_1$ such that $I[s]\cap
\tilde{I}_P=\emptyset$ and $I_P\subset I_2$ we have
$$\textrm{Graph}(q_2)\cap\left(c(d)\d^{-1}\right)\widehat{P}=\emptyset\:.$$
\end{quote}
Of course, the same is true for the symmetric relation, {\it i.e.}
replacing the index $1$ with $2$ and vice versa.

\item  again, based on Lemma \ref{CZdec}, we deduce
\beq\label{ob5}
\forall\:\:P=[\vec{\a},I_P]\in\p_1\cup \p_2\:\: \textrm{we have}\:\:|\tilde{I}_P\cap I[c]|\lesssim_{d}{\d}^{\frac{1-\ep_0}{d}}|I_P|\,.
\eeq
\end{enumerate}
\end{obs}

\begin{l1}\label{sept}[\textsf{Interaction of separated trees}]

Let $\left\{\p_j\right\}_{j\in\left\{1,2\right\}}$  be two $\d^{-1}$-separated
trees with tops $P_j=[\vec{\a}_j,I_0]$. Then,
for any $f,\:g\in L^{2}(\TT)$ and $n\in \N$, we have that
\beq\label{v21} \left|\left\langle
{T^{\p_1}}^*f,\,{T^{\p_2}}^*g\right\rangle\right|\lesssim_{n,d}{\d}^n\left\|f\right\|_{L^{2}(\tilde{I}_0)}\left\|g\right\|_{L^{2}(\tilde{I}_0)}+\left\|\chi_{I[c]}{T^{\p_1}}^*f\right\|_2\left\|\chi_{I[c]}{T^{\p_2}}^*g\right\|_2\:.
\eeq
\end{l1}
\begin{proof}
In what follows, for conciseness, we will only present the main steps and adaptations for our proof. For further details see the analogue proof of Lemma 2 in \cite{q} as well as that of Lemma 4 in \cite{f}.

\begin{itemize}
\item We start by noticing that definition \eqref{seps} applies in our context to $J=\tilde{I}_0$. Applying now the steps in Lemma \ref{CZdec} we write \eqref{parthigh} in our setting, that is
\beq\label{parthigh1}
J_{l}(q,\l)=\bigcup_{W\in CZ_{(q,\l)}(J)} W\,,
\eeq
for $J=\tilde{I}_0$, $q=q_{1,2}$ and $\l=c_0(d)\,\d^{-1}$ chosen as in \eqref{seps}.

Notice that with the above notations and conventions, we have
\beq\label{parthigh2}
J_{l}(q,\l)=\tilde{I}_0\setminus I[s]\:.
\eeq

\item For $j\in\{1,2\}$, define the following tile-sets:
\beq\label{1c}
\p_j(I[s]):=\{P=[\vec{\a},I_P]\in\p_j\,|\,100\tilde{I}_P\cap I[s]\not=\emptyset\}\:,
\eeq
and, for each $W\in J_{l}(q,\l)$
\beq\label{1w}
\p_j(W):=\{P=[\vec{\a},I_P]\in\p_j\setminus \p_j(I[s])\,|\,I_P\subset W\:\:\&\:\: |\tilde{I}_P|\leq \frac{1}{5}\,|W|\}\:.
\eeq
Notice that this decomposition defines for each $j\in\{1,2\}$ a partition
\beq\label{partj}
\p_j=\p_j(I[s])\cup \bigcup_{W\in CZ_{(q,\l)}(J)}\p_j(W)\:.
\eeq

\item With this done, for $j\in\{1,2\}$, we define corresponding labels $K_j$ that can take values $K_j=I[s]$ or $K_j=W$ with $W$ as above, and set
\beq\label{partj1}
T^{*}_{j,K_j}=\sum_{P\in \p_j(K_j)}T^{*}_{P}\,.
\eeq

Deduce that
\beq\label{partj3}
\left\langle{T^{\p_1}}^*,\:{T^{\p_2}}^*
\right\rangle=\sum_{K_j}\left\langle{T^{*}_{1,K_1}},\:{T^{*}_{2,K_2}}
\right\rangle\:.
\eeq

\item Define a real-valued function $\f\in C_{0}^{\infty}(\R)$ with the following properties:
\begin{itemize}
\item $ supp\:\f\subset\left\{\frac{1}{4}\leq |x|\leq \frac{1}{2}\right\}$
\item $ \f\:is\: even $
\item $ |\hat{\f}(\xi)-1|\lesssim_{n} |\xi|^n \:\:\:\:\forall\:|\xi|\leq 1\:\:and\:n\: big\:enough $
\item $ |\hat{\f}(\xi)|\lesssim_{n} |\xi|^{-n}\:\:\:\:\forall\:|\xi|\geq 1$
    \end{itemize}

Also for $W$ as before, we define
$$d_{j,W}:=\min\{|I_P|\,|\,P=[\vec{\a},I_P]\in \p_j(W)\}\:.$$

Now, for $j\in\left\{1,2\right\}$ and $W\in CZ_{(q,\l)}(J)$, let
\beq\label{fprop1}
\f_{j,l}(x)=(\d^{1/3}d_{j,W})^{-1}\f((\d^{1/3}d_{j,W})^{-1}x)\:.
\eeq
and define the operators
\beq\label{fprop2}
\tilde{\f}_{j,W}\::L^2(\R)\longrightarrow L^2(\R)\:\:\:\: by\:\:\tilde{\f}_{j,W}f=\f_{j,W}*f\,,
\eeq
and
\beq\label{fprop3}
\Phi_{j,W}\::L^2(\R)\longrightarrow L^2(\R)\:\:\:\:by\:\:\Phi_{j,W}=\left(\prod_{l=1}^d M_{l,a_l^j}\right)\tilde{\f}_{j,W} \left(\prod_{l=1}^d M^{*}_{l,a_l^j}\right)\:.
\eeq
where in the last line we simply assumed that  $q_j\in \Q_{d-1}$ is the central polynomial corresponding to $P_j=[\vec{\a}_j,I_0]$ and
\beq \label{shifsss}
Q_j(y)=\sum_{l=1}^{d}a_l^j\:y^l\:,
\eeq
is the unique polynomial in $\Q_{d}$ such that $Q_j(x)=\frac{d}{dx}\,q_j$ and $Q_j(0)=0$.

\item Following now similar reasoning with those in Lemma 2 in \cite{q}, for $j\in\left\{1,2\right\}$ and $W$ as before,
we decompose
\beq\label{dec}
T^{*}_{j,W}f=\Phi_{j,W}{T^{*}_{j,W}}f\:+\:\Omega_{j,W}f\,,
\eeq
and deduce that
\beq \label{om}
 \left\|\Omega_{j,W}\right\|_2\lesssim_{n}\d^{n}\,,
\eeq
and for any $n\in\N$
\beq \label{fi} \left|\left\langle\Phi_{1,W_1}{T^{*}_{1,W_1}}f
,\Phi_{2,W_2}{T^{*}_{2,W_2}}g\right\rangle\right|\lesssim_{n}\d^{n}\left\|f\right\|_2\left\|g\right\|_2\:.
\eeq
Remark that if $W_1,\,W_2$ are not the same or adjacent then the RHS of \eqref{fi} is in fact zero.

\item After some relatively involved computations one concludes that

\beq\label{nonc}
\sum_{{{K_1,K_2}\atop{K_1\not=I[s]}}\atop{\textrm{or}\:K_2\not=I[s]}}\left|\left\langle
{T^{*}_{1,K_1}}f,\,{T^{*}_{2,K_2}}g\right\rangle\right|\lesssim_{n,d}
{\d}^n\left\|f\right\|_{L^{2}(\tilde{I}_0)}\left\|g\right\|_{L^{2}(\tilde{I}_0)}\:,\eeq
and
\beq\label{crit}
\left|\left\langle{T^{*}_{1,I[s]}}f,\,{T^{*}_{2,I[s]}}g\right\rangle\right|\lesssim_{n,d}
{\d}^n\left\|f\right\|_{L^{2}(\tilde{I}_0)}\left\|g\right\|_{L^{2}(\tilde{I}_0)}+
\left\|\chi_{I[c]}{T^{\p_1}}^*f\right\|_2\left\|\chi_{I[c]}{T^{\p_2}}^*g\right\|_2\:.\eeq
finishing our proof.
\end{itemize}

\end{proof}

\begin{d0}\label{nor}[\textsf{Normal tree}]

A tree $\p$ with top $P_0=[\vec{\a}_0,I_0]$ is called \textbf{normal} if for any
$P=[\vec{\a},I]\in\p$ we have $100 I\cap (I_0)^c=\emptyset\;.$
\end{d0}

\begin{obs}\label{norob} Notice that if $\p$ is a normal tree as above
then $$supp\:{{T^{\p}}^{*}f}\subseteq I_0\:.$$
\end{obs}

\begin{d0}\label{row} [\textsf{Row}]

A \emph{row} is a collection $\p=\bigcup_{j\in \N}\p^{j}$ of normal
trees $\p^{j}$ with tops
$P^{j}_0=[\vec{\a}^{j}_0,I^j_0]$ such that the
$\left\{I^j_0\right\}$ are pairwise disjoint.
\end{d0}

\begin{l1}\label{rt}[\textsf{Row-tree interaction}]

Let $\p$ be a row as above, let $\p'$ be a tree with top $P'_0=[\vec{\a}'_0,I'_0]$ and suppose
that $\forall\: j\in \N$, $I_0^{j}\subseteq I_0'$ and $\p^{j},\:\p'$
are $\d^{-1}$separated trees; denote by $I^{j}[c]$ the critical intersection
set between each $\p^{j}$ and $\p'$.

Then for any $f,\:g\in L^{2}(\TT)$ and $n\in \N$ we have that
$$\left|\left\langle {T^{\p'}}^*f,{T^{\p}}^*g\right\rangle\right|
\lesssim_{n,d}{\d}^n\left\|f\right\|_{2}\left\|g\right\|_{2}+
\left\|\sum_{j}\chi_{I^{j}[c]}{T^{\p'}}^*f\right\|_2\left\|\sum_{j}\chi_{I^{j}[c]}{T^{\p^{j}}}^*g\right\|_2\:.$$
\end{l1}

\begin{proof}
Again this proof follows the corresponding proof in \cite{q} with many elements borrowed from the proof of Lemma \ref{sept} (see also Lemma 5 in \cite{f}). For concreteness, as before, we present only a brief outline of the proof.

\begin{itemize}
\item notice that it is enough to show that for each $j$, we have
\beq\label{crit}
\eeq
$$\left|\left\langle {T^{\p'}}^*f,{T^{\p^{j}}}^*g\right\rangle\right|\lesssim_{n}\d^n\left(\left\|M(Mf)\right\|_{L_2(I_0^j)}+\left\|M(M({T^{\p'}}^*f))\right\|_{L_2(I_0^j)}\right)\left\|g\right\|_{L_2(I_0^j)}$$
$$+\left\|\chi_{I^j[c]}{T^{\p'}}^*f\right\|_2\left\|\chi_{I^j[c]}{T^{\p^{j}}}^*g\right\|_2\:,$$
where as usual here $M$ stands for the Hardy-Littewood maximal function.

\item We repeat the steps from Lemma \ref{CZdec}, as follows:

Define
\beq\label{parthigh1}
J_{l}(q,\l)=\bigcup_{W\in CZ_{(q,\l)}(J)} W\,,
\eeq
for $J=I_0^{j}$, $q$ the $(P'_0,\,P^j_0)$-interaction polynomial and $\l=c_0(d)\,\d^{-1}$ chosen as in \eqref{seps}.

\item Define the following tile-sets:
\beq\label{1cr}
\p^j(I^j[s]):=\{P=[\vec{\a},I_P]\in\p^j\,|\,100\tilde{I}_P\cap I^j[s]\not=\emptyset\}\:,
\eeq
\beq\label{2cr}
\p'(I^j[s]):=\{P=[\vec{\a},I_P]\in\p'\,|\,100\tilde{I}_P\cap I^j[s]\not=\emptyset\}\:,
\eeq
and, for each $W\in J_{l}(q,\l)$
\beq\label{1wr}
\p^j(W):=\{P=[\vec{\a},I_P]\in\p^j\setminus \p^j(I[s])\,|\,I_P\subset W\:\:\&\:\: |\tilde{I}_P|\leq \frac{1}{5}\,|W|\}\:,
\eeq
\beq\label{1wr1}
\p'(W):=\{P=[\vec{\a},I_P]\in\p'\setminus \p'(I[s])\,|\,I_P\subset W\:\:\&\:\: |\tilde{I}_P|\leq \frac{1}{5}\,|W|\}\:.
\eeq
We thus get that
\beq\label{partj1r}
\p^j=\p^j(I[s])\cup \bigcup_{W\in CZ_{(q,\l)}(J)}\p^j(W)\:,
\eeq
and
\beq\label{partj2r}
\p'=\p'(I[s])\cup \bigcup_{W\in CZ_{(q,\l)}(J)}\p'(W)\:.
\eeq
\item For $K=I^j[s]$ or $K=W$ with $W$ as above, we let
\beq\label{partj1rr}
T^{*}_{j,K}=\sum_{P\in \p^j(K)}T^{*}_{P}\,,
\eeq
and
\beq\label{partj2rr}
T^{'*}_{K}=\sum_{P\in \p'(K)}T^{*}_{P}\,.
\eeq
\item further define
$$d_{W}:=\min\{|I_P|\,|\,P=[\vec{\a},I_P]\in \p^j(W)\cup \p'(W)\}\:,$$
and with $\f$ as in Lemma \ref{sept} we let
\beq\label{fprop1r}
\f(x)=(\d^{1/3}d_{W})^{-1}\f((\d^{1/3}d_{W})^{-1}x)\:.
\eeq
Next, define the operators
\beq\label{fprop2r}
\tilde{\f}_{W}\::L^2(\R)\longrightarrow L^2(\R)\:\:\:\: by\:\:\tilde{\f}_{W}f=\f_{W}*f\,,
\eeq
and
\beq\label{fprop3r}
\Phi_{j,W}\::L^2(\R)\longrightarrow L^2(\R)\:\:\:\:by\:\:\Phi_{j,W}=\left(\prod_{l=1}^d M_{l,a_l^j}\right)\tilde{\f}_{W} \left(\prod_{l=1}^d M^{*}_{l,a_l^j}\right)\:,
\eeq
and
\beq\label{fprop4r}
\Phi'_{W}\::L^2(\R)\longrightarrow L^2(\R)\:\:\:\:by\:\:\Phi'_{W}=\left(\prod_{l=1}^d M_{l,a'_l}\right)\tilde{\f}_{W} \left(\prod_{l=1}^d M^{*}_{l,a'_l}\right)\:.
\eeq
where, as usual, $Q_j(y)=\sum_{l=1}^{d}a_l^j\:y^l$ and $Q(y)=\sum_{l=1}^{d}a'_l\:y^l$ are the unique polynomials in $\Q_{d}$ with no constant term whose derivatives equal the central polynomials corresponding to $P^{j}_0$ and $P'_0$ respectively.

\item Finally, we end this decomposition by setting
\beq\label{decr}
T^{*}_{j,W}f=\Phi_{j,W}{T^{*}_{j,W}}f\:+\:\Omega_{j,W}f\,,
\eeq
and
\beq\label{decrr}
T^{'*}_{W}f=\Phi'_{W}{T^{'*}_{W}}f\:+\:\Omega'_{W}f\,.
\eeq
\item Withe the above notations and conventions, one notices that
\beq \label{rowdec}
\begin{array}{rl}
\left\langle {T^{\p'}}^*f,T^{*}_{j,W}g\right\rangle=\left\langle
\Phi'_{W}{T^{\p'}}^*f,\Phi_{j,W}T^{*}_{j,W}g\right\rangle+\\
\left\langle\Phi'_{W}{T^{\p'}}^*f,\Omega_{j,W}g
\right\rangle+\left\langle\Omega'_{W}f,T^{*}_{j,W}g\right\rangle=I+II+III\:.
\end{array}
\eeq
\item Using now similar reasonings with the ones in Lemma \ref{sept} we have:

- for the first term
\beq \label{term1}
|I|\lesssim_{n}\d^n\left\langle M\left({T^{\p'}}^*f\right),|T^{*}_{j,W}g|\right\rangle\lesssim
\d^n\left\|M({T^{\p'}}^*f)\right\|_{L^2(\frac{3}{2}W\cap
I_0^j)}\left\|g\right\|_{L^2(\frac{3}{2}W\cap I_0^j)}\,,
\eeq
where here we used
\beq \label{rowdec11}
\left|\Phi'_{W}\Phi_{W}h(x)\right|\lesssim_{n}\d^n \left(u_{W}*|h|\right)(x)\,,
\eeq
with $h\in L^1(\TT)$ and
$$u_{W}(x):=(\d^{1/3}d_{W})^{-1}\chi_{\left\{|t|\leq2\d^{1/3}d_{W}\right\}}(x)\:.$$

- for the second term
\beq \label{term2}
\begin{array}{rl}
|II|\lesssim\left\langle \r_{W}*\left\{\chi_{\frac{3}{2}W\cap I_0^j}M({T^{\p'}}^*f)\right\},|g|\right\rangle\\
\lesssim\d^n\left\|M({T^{\p'}}^*f)\right\|_{L^2(\frac{3}{2}W\cap
I_0^j)}\left\|Mg\right\|_{L^2(\frac{3}{2}W\cap I_0)}\:,
\end{array}
\eeq
where here we make use of the fact that
\beq \label{omr}
\left|\Omega_{j,W}^{*}h(x)\right|\lesssim\left(\r_{W}*|h|\right)(x)\:\:\:\textrm{and}\:\:\:\left\|\r_{W}\right\|_1\lesssim (\d)^{n/3}\,,\,,
\eeq
where
 $$\r_{W}(y):=\sum_{2^k\leq(d_{W})^{-1}}(\d^{1/3}2^{k}d_{W})^{n}2^{k}\chi_{[-2^{-k},2^{-k}]}(y)\:.$$

- for the third term, we notice that
\beq \label{omrst}
x\in \frac{3}{2}W\:\:\Rightarrow\:\:\left|{\Omega'_{W}}^{*}h(x)\right|\lesssim\left(\r_{W}*|h|\right)(x)\,,
\eeq
and as a consequence
\beq \label{term3}
\begin{array}{rl}
|III|\lesssim\left\langle |f|,\r_{W}*|{T_{j,W}}^*g|
\right\rangle\lesssim\d^n\left\|Mf\right\|_{L^2(\frac{3}{2}W\cap
I_0^j)}\left\|g\right\|_{L^2(\frac{3}{2}W\cap
I_0^j)}\:.
\end{array}
\eeq

\item Putting together \eqref{term1} - \eqref{term3}, we obtain
\beq \label{rowe0}
\begin{array}{cl}
\left|\left\langle {T^{\p'}}^*f,{T_{j,W}}^*g\right\rangle\right|\lesssim\\
\d^n\left(\left\|Mf\right\|_{L^2(\frac{3}{2}W\cap
I_0^j)}+\left\|M({T^{\p'}}^*f)\right\|_{L^2(\frac{3}{2}W\cap
I_0^j)}\right)\left\|Mg\right\|_{L^2(\frac{3}{2}W\cap
I_0^j)}\,,
\end{array}
\eeq
and after applying Cauchy-Schwarz we deduce
\beq \label{rowe}
\begin{array}{cl}
\left|\left\langle {T^{\p'}}^*f,{T^{\p^{j}}}^*g\right\rangle\right|\lesssim_{n}\\
\d^n\left(\left\|Mf\right\|_{L^2(I_0^j)}+\left\|M({T^{\p'}}^*f)\right\|_{L^2(
I_0^j)}\right)\left\|Mg\right\|_{L^2(I_0^j)}+ \left|\left\langle
{T^{\p'}}^*f,T^{*}_{j,I^j[s]}g\right\rangle\right|\:.
\end{array}
\eeq

\item It remains to estimate the last term in the right-hand side of \eqref{rowe}.

Recall the definition of $\p'(I^j[s])$ and further set
$$\p'(I_0^j):=\left\{P=[\vec{\a},I]\in\p'\setminus\p'(I^j[s])\,|\:|I|\leq \frac{|I_0^j|}{100}\right\}\,,$$
and
$$\grave{\p}:=\p'\setminus\left(\p'(I^j[s])\cup \p'(I_0^j)\right)\,.$$
Deduce thus that
\beq \label{decro}
{T^{\p'}}^*f={T^{\p'(I_0^j)}}^*f\:+\:{T^{\grave{\p}}}^*f\:+\:T^{'*}_{I^j[s]}\:.
\eeq

For the first term, from Lemma \ref{sept}, we deduce that
\beq \label{tr1}
\left|\left\langle {T^{\p'(I_0^j)}}^*f,T^{*}_{j,I^j[s]}g\right\rangle\right|\lesssim_{n}
{\d}^n\left\|f\right\|_{L^2(I_0^j)}\left\|g\right\|_{L^2(I_0^j)}\:.
\eeq
For the remaining terms, we follow similar reasonings with the ones above, and get
\beq \label{tr2}
\left|\left\langle {T^{\grave{\p}}}^*f,T^{*}_{j,I^j[s]}g\right\rangle\right|\lesssim_{n}{\d}^n\left(\left\|Mf\right\|_{L^2( I_0^j)}+\left\|M\left\{\M({T^{\p'}}^*f)\right\}\right\|_{L^2( I_0^j)}\right)\left\|Mg\right\|_{L^2(I_0^j)}
\eeq
and
\beq \label{tr3}
\eeq
$$\left|\left\langle T^{'*}_{I^j[s]},T^{*}_{j,I^j[s]}g\right\rangle\right|\lesssim_{n}{\d}^n\left(\left\|Mf\right\|_{L^2( I_0^j)}+\left\|M\left\{\M({T^{\p'}}^*f)\right\}\right\|_{L^2( I_0^j)}\right)\left\|Mg\right\|_{L^2(I_0^j)}$$
$$+\left\|\chi_{I^j[c]}{T^{\p'}}^*f\right\|_2\left\|\chi_{I^j[c]}{T^{\p}}^*g\right\|_2\:,$$
where
$$\M({T^{\p'}}^*f)=^{def}\sup_{m\in N}\left|\sum_{{P=[\vec{\a},I]\in\ \p'}\atop{|I|\geq 2^{-m}}}T_{P}^{*}f\right|\:.$$

\item So to summarize, we proved that
$$\left|\left\langle {T^{\p'}}^*f,{T^{\p^j}}^*g\right\rangle\right|\lesssim_{n}\d^n\left(\left\|Mf\right\|_{L^2( I_0^j)}+\left\|M\left\{\M({T^{\p'}}^*f)\right\}\right\|_{L^2( I_0^j)}\right)\left\|Mg\right\|_{L^2(I_0^j)}$$
$$+\left\|\chi_{I[c]^j}{T^{\p'}}^*f\right\|_2\left\|\chi_{I[c]^j}{T^{\p^j}}^*g\right\|_2\:.$$

Now the conclusion follows if we add the observation that
\beq
\label{treemax} \M({T^{\p'}}^*f)\leq Mf\:+\:M({T^{\p'}}^*f)\:.
\eeq
\end{itemize}

\end{proof}

\begin{l1}\label{ed} [$L^2$-\textsf{small adjoint tree support}]

Let $\p\:$ be a tree with top
$P_{0}=[\vec{\a}_{0},I_{0}]$; suppose also that
we have a set $A\subseteq \tilde{I_0}$ with the property that
 \beq
\label{masx} \:\exists\:\d\in(0,1)\:\:\:st\:\:\:\:\:\forall
\:P=[\vec{\a},I]\in\p\:\:\operatorname{we\:have\:}\:\:|I^{*}\cap A|\leq \d
|I|. \eeq

Then for any $f\in L^2(\TT)$ we have
\beq \label{cut}
\left\|\chi_{A}{T^{\p}}^{*}f\right\|_2\lesssim \d^{\frac{1}{2}}\left\|f\right\|_{2}.
\eeq
\end{l1}
\begin{proof}
This follows immediately from the corresponding proof of Lemma \ref{lptree} in the next section.
\end{proof}

We are now in the position to state the main result of this section.
$\newline$

\noindent\textbf{Main Lemma.} {\it Let $\p\subset \P_n$ be an $L^{\infty}-$forest of generation $n$.

\noindent Then there exists $\eta=\eta(d)\in (0,1)$ such that
$$\|T^{\p}f\|_2\lesssim 2^{-\frac{n}{2}\,\eta}\,\|f\|_2\:.$$

Moreover, if $\p$ is normal and $2^{100\,n\,d}$-separated,\footnote{As expected, an $L^{\infty}-$forest $\p$
is called normal if all the trees inside are normal; same principle applies for the $\d^{-1}-$separateness condition.} then
- decomposing $\p$ canonically into a union of rows $\{\r_j\}$ - one has
\beq \label{almostortho}
\|{T^{\p}}^{*}f\|_2^2\lesssim \sum_{j}\|{T^{\r_j}}^{*}f\|_2^2\,+\,2^{-5\,n}\,\|f\|_2^2\,,
\eeq
from which one deduces the improved bound
$$\|T^{\p}f\|_2\lesssim 2^{-\frac{n}{2}}\,\|f\|_2\:.$$}

\begin{o0}\label{Almoortog}
\textit{Relation \eqref{almostortho} should be regarded as a strong almost orthogonality relation arising from the good geometric properties imposed on $\p$ and which essentially states - up to a negligible term - that one expects
$$\|{T^{\p}}^{*}f\|_2^2\lesssim \sum_{j}\|{T^{\r_j}}^{*}f\|_2^2\:.$$}
\end{o0}

\begin{proof}

Recalling Definition \ref{infforest}, our hypothesis implies that we can decompose our forest $\p$ as
\beq \label{decp}
 \p=\bigcup_{j=1}^{c\,2^n} \r_j\,,
\eeq
 with each $\r_j$ being a maximal collection of spatially disjoint trees. Further, we decompose each $\r_j$ in a disjoint union of maximal trees $\{\T_{j,k}\}_k$. If $I_{jk}$ stands for the time-interval of the top of $\T_{j,k}$, define the boundary component
$$\T_{j,k}^{bd}:=\{P\in\T_{j,k}\,|\,100 I_P\cap (I_{jk})^c\not=\emptyset\}\,.$$

Also let $\tilde{\T}_{j,k}$ be the set of tiles defined inductively as follows: we erase the collection of minimal tiles in $\T_{j,k}$ and reapply this step to the newly defined $\T_{j,k}$ for $100\,n\,d$ times, moment at which we stop and collect all the successively erased tiles in the set $\tilde{\T}_{j,k}$.

Next, we notice that the set $\bigcup_{j,k} \T_{j,k}^{bd}$ is a sparse $L^{\infty}$-forest while the set $\bigcup_{j,k} \tilde{\T}_{j,k}$ can be decomposed in at most $c\,n$ negligible sets. Deduce thus that the family
$$\bigcup_{j,k} \T_{j,k}^{bd}\cup \tilde{\T}_{j,k}\,,$$
can be easily treated via Proposition 1. Thus, we can erase this set of tiles from our initial forest.

We are now left with the case $\p$ normal and $2^{100\,n\,d}$-separated $L^{\infty}-$forest. In this new context, we notice that for the updated family $\p$ the analogue of \eqref{decp} is precisely the row decomposition of $\p$.

The key argument now  is that the operators $\{T^{\r_j}\}_j$ are almost orthogonal. More precisely, for $k\not=j$, we claim that
\begin{itemize}
\item $\|{T^{\r_k}}^{*}\,T^{\r_j}\|_{2\mapsto 2}=0$;

\item $\|T^{\r_k}\,{T^{\r_j}}^{*}\|_{2\mapsto 2}\lesssim 2^{-10\,n}$.
\end{itemize}

The first item is a direct consequence of the pairwise disjointness of the sets $\{\textrm{supp}\,T^{\r_j}\}_j$. For the second item one needs to make use of the strong ($2^{100\,n\,d}$)-separateness hypothesis and successively apply Lemmas \ref{rt} and \ref{ed}. We leave these details to the reader.
\end{proof}

\subsubsection{\bf $L^p-$results}
$\newline$

In this section, we will discuss several $L^p$ versions of the results presented
in the previous section.

\begin{l1}\label{trp}[$L^p$-\textsf{uniform mass tree estimate}]

Fix $\delta\in(0,1]$ and let $\p\subseteq\P$ be a tree with spacial support $I_0$ such that
\beq\label{masscon}
A_{\P,I_0}(P)<\delta\:\:\:\:\forall\:\:\:P\in\p\:.
\eeq
Then, for $1<p<\infty$, we have
\beq\label{treep}
\left\|T^\p\right\|_p\lesssim_{p,d}\delta^{\frac{1}{p}}\:. \eeq
\end{l1}

\begin{obs}\label{treeimp}
The $L^p$-tree lemma above follows directly from Lemma \ref{treecutp} below by simply setting $A=[0,1]$. Indeed, with the notations in Lemma \ref{treecutp}, one has that \eqref{masscon} immediately implies
\beq \label{eqj1}
\frac{E_{A}^{\p}(J)}{|J|}\lesssim \delta \:\:\:\:\:\:\forall\:J\in CZ(\I_{\p}^{min})\,.
\eeq
\end{obs}

\begin{l1}\label{treecutp}[$L^p$-\textsf{small tree support}]

Let $\p\:$ be a tree. Define
\beq \label{mastree}
\I_{\p}^{min}:=\{I\,|\,\exists\:P=[\vec{\a},I]\in\p\:\:\:\textrm{minimal relative to}\:\leq\:\textrm{in}\:\p\}\:,
\eeq
and let $CZ(\I_{\p}^{min})$ be the Calderon-Zygmund decomposition of the interval $[0,1]$ relative to $\I_{\p}^{min}$.

Assume we are given $A\subseteq[0,1]$ measurable set and define
\beq \label{eqj}
E_{A}^{\p}(J):=\bigcup_{P\in\p} (A\cap E(P)\cap J)\:\:\:\:\:\:\forall\:J\in CZ(\I_{\p}^{min})\,.
\eeq
Then, for $1<p<\infty$, we have
\beq \label{mastree}
\left\|\chi_{A}\,T^{\p}f\right\|_p\lesssim_{p,d} \left(\sup_{J\in CZ(\I_{\p}^{min})}\,\frac{|E_{A}^{\p}(J)|}{|J|}\right)^{\frac{1}{p}}\,\left\|f\right\|_{p}.
\eeq
 \end{l1}

\begin{proof}
We start by setting the parameters of our tree, that is: we
fix the top $P_0=[\vec{\a_0},I_0]$, and the corresponding central polynomial $q_0\in \Q_{d-1}$. Once we specified
the polynomial $q_0$ there is a unique $Q_0\in \Q_{d}$ such that
 $Q_0(x)=\frac{d}{dx}\,q_0$ and $Q_0(0)=0$. Assume wlog that
\beq \label{shifs}
Q_0(y)=\sum_{j=1}^{d}a_j^0\:y^j\:.
\eeq
Then, proceeding as in \cite{q}, we define
\beq \label{shif}
\T^{\p}:=\left(\prod_{j=1}^d M^{*}_{j,a_j^0}\right)\, T^{\p}\,\left(\prod_{j=1}^d M_{j,a_j^0}\right)\,,
\eeq
and $g(x):=\prod_{j=1}^d M^{*}_{j,a_j^0} f(x)$ and notice that
$$\left\|\chi_{A}\,T^{\p}f\right\|_p=\left\|\chi_{A}\,\T^{\p}g\right\|_p\:.$$
Following similar reasonings with the ones in \cite{f}, for a fixed $x\in \TT$, we further define
$$k_0(x):=\inf\{k\in\N\,|\,\exists\:P\in\p\:\:\textrm{s.t.}\:|I_P|=2^{-k}\:\&\:\chi_{E(P)}(x)\not=0\}\,,$$
$$k_1(x):=\sup\{k\in\N\,|\,\exists\:P\in\p\:\:\textrm{s.t.}\:|I_P|=2^{-k}\:\&\:\chi_{E(P})(x)\not=0\}\,,$$
and  notice that
\beq \label{treehilb}
\T^{\p}g(x):=\sum_{k=k_0(x)}^{k_1(x)}\,\int_{\TT}\psi_{k}(y)\,e^{i\,(Q(x)-Q_0(x)+Q_0(x-y)-Q(x-y))}\,g(x-y)\,dy\:.
\eeq
Relation \eqref{treehilb} above is the key place where we use the convexity of the tree and thus the very reason for which we need to remove the possible boundary effect by requiring item 2 in Definition \ref{tree}.

With this, we notice that
$$|\chi_{A}(x)\,\T^{\p}g(x)|$$
$$\leq\chi_{A}(x)\,\sum_{k=k_0(x)}^{k_1(x)}\,\int_{\TT}\,|\psi_{k}(y)|\,|e^{i\,(Q(x)-Q_0(x)+Q_0(x-y)-Q(x-y))}-1|\,|g(x-y)|\,dy$$
$$+\chi_{A}(x)\,\left|\sum_{k=k_0(x)}^{k_1(x)}\,\int_{\TT}\psi_{k}(y)\,g(x-y)\,dy\right|=:\A(x)\,+\,\B(x)\:.$$
Let now
\beq \label{fmaxa}
M_{A}^{\p}g(x):=\left\{
                        \begin{array}{ll}
                        \sup_{I\supset J}\frac{1}{|I|}\int_{I}|g|, \  \mbox{if} \  x\in E_{A}^{\p}(J) \:\textrm{and}\:J\in CZ(\I_{\p}^{min})\\
                        0 \qquad, \  \mbox{otherwise}
                        \end{array} \right.
\:.\eeq

Now, applying some elementary reasonings, we further deduce:
\beq \label{keytr}
\A(x)\lesssim_d M_{A}^{\p} g(x)\:,
\eeq
and
\beq \label{keytr}
\B(x)\lesssim M_{A}^{\p}(R*g)(x)\:,
\eeq
where we set $R(y)=\sum_{k\in \N D}\psi_{k}(y)$ and assumed wlog, based on Observation \ref{redu},
that $\p\subset\bigcup_{k\in\N}\P_{k D}$.

It only remains to notice that
\beq \label{keymax}
\left\|M_{A}^{\p}g\right\|_p\lesssim \left(\sup_{J\in CZ(\I_{\p}^{min})}\,\frac{|E_{A}^{\p}(J)|}{|J|}\right)^{\frac{1}{p}}\,\left\|M g\right\|_p
\lesssim_{p}\left\|g\right\|_p\,,
\eeq
and
\beq \label{keyhilb}
\left\|R*g\right\|_p\lesssim_p\left\|g\right\|_p\,.
\eeq
Thus, we conclude that \eqref{mastree} holds.
\end{proof}

\begin{l1}\label{lptree}[$L^p$-\textsf{small adjoint tree support}]

Let $\p\:$ be a tree with spacial support $I_0$. Recalling \eqref{partint1}, we define
\beq \label{mastrees}
\I_{\p^{*}}^{min}:=\bigcup_{j=1}^{14}\{I_{P*}^{j}\,|\,\exists\:P=[\vec{\a},I_P]\in\p\:\:\textrm{minimal relative to}\:\leq\:\textrm{in}\:\p\}\:,
\eeq
and let $CZ(\I_{\p^*}^{min})$ be the Calderon-Zygmund decomposition of the interval $[0,1]$ relative to $\I_{\p^*}^{min}$.

Assume we are given $A\subseteq[0,1]$ measurable set and define
\beq \label{eqjs}
E_{A}^{\p^*}(J):=A\cap J\:\:\:\:\:\:\forall\:J\in CZ(\I_{\p^*}^{min})\,.
\eeq
Then, for $1<p<\infty$, we have
\beq \label{mastrees}
\left\|\chi_{A}{T^{\p}}^{*}f\right\|_p\lesssim_{p,d} \left(\sup_{J\in CZ(\I_{\p^*}^{min})}\,\frac{|E_{A}^{\p^*}(J)|}{|J|}\right)^{\frac{1}{p}}\,\sup_{P\in\p} A_{\P,I_0}(P)^{\frac{1}{p'}}\,\left\|f\right\|_{p}\:.
\eeq
\end{l1}
\begin{proof}

We first notice that by applying the same reasonings as in the previous proof, specifically relying on \eqref{shif}, one can assume wlog that our tree lives at frequency zero, or with other words that $Q_0\equiv 0$ in \eqref{shifs}.\footnote{Alternatively one can apply the reasonings in our proof to the operator $\chi_{A}{\T^{\p}}^{*}$ where here $\T^{\p}$ is defined by \eqref{shif}.}

Fix now $x\in \TT$ and notice that with the previous notations one has
\beq \label{tst}
\chi_{A}{T^{\p}}^{*}f(x)=\sum_{J\in CZ(\I_{\p^*}^{min})}\chi_{J\cap
A}\left\{\sum_{P\in\p}{T_{P}}^{*}f(x)\right\}
\eeq
Next, fixing $J\in CZ(\I_{\p^*}^{min})$ and assuming wlog that $x\in J$, we deduce
$$\left|{T^{\p}}^{*}f(x)-\frac{1}{|J|}\int_{J}{T^{\p}}^{*}f(s)ds\right|$$
\beq \label{ets}
=\left|\frac{1}{|J|}\int_{J}\left\{\sum_{{P\in\p}\atop{2^{-k}=|I_P|\geq
|J|}}\int_{\TT}\left[\f_k(x-y)-\f_k(s-y)\right]f(y)\chi_{E(P)}(y)dy\right\}ds\right|
\eeq
$$\lesssim\sum_{{P\in\p}\atop{2^{-k}=|I_P|\geq
|J|}}2^{k}\,|J|\,\frac{\int_{E(P)}|f|}{|I_P|}\,.$$
Thus, from \eqref{ets}, we deduce that
\beq \label{ets1}
\eeq
$$\sum_{J\in CZ(\I_{\p^*}^{min})}\chi_{J\cap A}\,\left|{T^{\p}}^{*}f(x)-\frac{1}{|J|}\int_{J}{T^{\p}}^{*}f(s)ds\right|$$
$$\lesssim \sum_{J\in CZ(\I_{\p^*}^{min})}\chi_{J\cap A}\sum_{{P\in\p}\atop{J\subset I_{P^*}}}\frac{|J|}{|I_P|}\,\frac{\int_{E(P)}|f|}{|I_P|}\:.$$

Denote with
\beq \label{maxpst}
M_{\p^*}f(x)=\sum_{J\in CZ(\I_{\p^*}^{min})}\chi_{J}(x)\sup_{J\subseteq I}\frac{1}{|I|}\int_{I}|f|(s)ds\:.
\eeq
Deduce from \eqref{ets1} and \eqref{maxpst} that
\beq \label{maxpst1}
\eeq
$$\|\chi_{A}{T^{\p}}^{*}f\|_{p}\lesssim \left(\sup_{J\in CZ(\I_{\p^*}^{min})}\,\frac{|E_{A}^{\p^*}(J)|}{|J|}\right)^{\frac{1}{p}}\,\|M_{\p^*} ({T^{\p}}^{*}f)\|_p$$
$$+\left(\sup_{J\in CZ(\I_{\p^*}^{min})}\,\frac{|E_{A}^{\p^*}(J)|}{|J|}\right)^{\frac{1}{p}}\,\left\|\sum_{J\in CZ(\I_{\p^*}^{min})}\chi_{J}\sum_{{P\in\p}\atop{J\subset I_{P^*}}}\,\frac{|J|}{|I_P|}\,\frac{\int_{E(P)}|f|}{|I_P|}\right\|_p\:.$$
Now  we notice that
\beq \label{l1}
\left\|\sum_{J\in CZ(\I_{\p^*}^{min})}\chi_{J}\sum_{{P\in\p}\atop{J\subset I_{P^*}}}\,\frac{|J|}{|I_P|}\,\frac{\int_{E(P)}|f|}{|I_P|}\right\|_1\lesssim \|f\|_1\,,
\eeq
and
\beq \label{l2}
\left\|\sum_{J\in CZ(\I_{\p^*}^{min})}\chi_{J}\sum_{{P\in\p}\atop{J\subset I_{P^*}}}\,\frac{|J|}{|I_P|}\,\frac{\int_{E(P)}|f|}{|I_P|}\right\|_{\infty}\lesssim \sup_{P\in\p} A_{0}(P)\,\|f\|_{\infty}\,.
\eeq
Interpolating now between \eqref{l1} and \eqref{l2} we get that for $1<p<\infty$
\beq \label{l3}
\left\|\sum_{J\in CZ(\I_{\p^*}^{min})}\chi_{J}\sum_{{P\in\p}\atop{J\subset I_{P^*}}}\,\frac{|J|}{|I_P|}\,\frac{\int_{E(P)}|f|}{|I_P|}\right\|_{p}\lesssim \left(\sup_{P\in\p} A_{0}(P)\right)^{\frac{1}{p'}}\,\|f\|_{p}\,.
\eeq
Also, we trivially have that
\beq \label{l5}
\|M_{\p^*} ({T^{\p}}^{*}f)\|_p\lesssim \|M({T^{\p}}^{*}f)\|_p\lesssim_{p} \|{T^{\p}}^{*}f\|_p\lesssim_{p,d} \sup_{P\in\p} A_{\P,I_0}(P)^{\frac{1}{p'}}\,\left\|f\right\|_{p}\:.
\eeq
Combining now \eqref{maxpst1}, \eqref{l3} and \eqref{l5} we conclude that \eqref{mastrees} holds.
\end{proof}

\subsection{\bf Proof of Proposition 2}
$\newline$

We start by restating the result that we need to prove\footnote{Excepting the display in the main statement, for notational simplicity, throughout the remaining part of the paper we will drop the $d-$dependence in all the inequalities signs.}:
$\newline$

\noindent\textbf{Proposition 2.} \textit{Let $\p\subseteq\P_n$ be a forest. Then there exists $\eta\in(0,1/2)$, depending only on the degree $d$,
such that for $1<p<\infty$ we have
$$\left\|T^{\p}\right\|_{p}\lesssim_{p,d} 2^{-n\,\eta\,(1-\frac{1}{p^*})}\:.$$}

\subsubsection{\bf The $L^2$ bound}
$\newline$

Recalling the setting described in Observation \ref{envir} as well as Definition \ref{forest} and appealing to a pigeonhole principle,  from now we can assume wlog that
\begin{itemize}
\item the family $\p$ can be written as
\beq\label{unionp}
\p=\bigcup_{k\geq 0}\p_n^{k}\,,
\eeq
is a (BMO)-forest of generation $n$ such that for each $k\geq 0$
\beq\label{unionp1}
\p_n^{k}\subset\P_n\:\:\textrm{is an}\:L^{\infty}-\textrm{forest of generation}\:n\,;
\eeq
\item the spacial support of the tiles in $\p_n^{k}$ is contained in $A_n^k$ a set that can be represented
as a finite union of maximal (disjoint) dyadic intervals;

\item there exists $c>0$ such that for each $k\in\N$ one has
\beq\label{supct}
A_n^k\prec_{n\,c} A_n^{k+1}\,;
\eeq

\item if $P=[\vec{\a}, I_P]\in\p_n^k$ then
\beq\label{unionp2}
\eeq
\begin{itemize}
\item $I_P\subseteq A_n^k$;

\item $I_P\nsubseteq A_n^{k+1}$;

\item $A_{\P, A_n^{k}}(P)\in (2^{-n},\,2^{-n+1}]$.
\end{itemize}
\end{itemize}

\begin{obs}\label{redf}
Now following similar reasonings with the ones described in Section \ref{redmp} and defining
\beq \label{cnkad}
 \check{\C}_n^k:=\left\{P\in\p_n^k\:|\:\operatorname{there\:are\: no\:chains}\:P\lneq P_{1}\lneq\ldots\lneq
P_{n}\:\&\:\left\{P_j\right\}_{j=1}^{n}\subseteq\p_n^k\:\right\}\:,
\eeq
we have that
\begin{itemize}
\item the set  $\check{\C}_n:=\bigcup_{k\geq 0}\check{\C}_n^k$ can be decomposed in a union of at most $n$ sparse forests; applying Proposition 1 to each of the resulting sparse forests we have that the associated operator $T^{\check{\C}_n}$ is under control.

\item erasing from each $\p_n^k$ the corresponding set $\check{\C}_{100 n d}^k$ one has that
 $$\p:=\bigcup_{k\geq 0}\p_n^k\,,$$
 is a BMO-forest of generation $n$ such that each $\p_n^k$ is an $L^\infty$-forest having the property that any two trees inside $\p_n^k$ are $2^{100\,n\,d}-$separated.
\end{itemize}
\end{obs}
Let now
\beq \label{anks}
\check{\p}_n^{k}:=\left\{\begin{array}{cc} P=[\vec{\a},I]\in \p_n^{k}\\
(\textrm{hence}\:I\subseteq A_n^{k}) \end{array}\:\big|\: \begin{array}{cc}
\:\textrm{if}\:J\subseteq A_n^{k+1} \:\:\textrm{s.t.}\:\:20 I\cap J^c\not=\emptyset\:\\
 \textrm{then}\: |I|\geq |J|
\end{array}\right\}\:,
\eeq
and define
\beq \label{pnbd}
\p_{n,bd}^{k,e}:=\p_n^{k}\setminus \check{\p}_n^{k}\:.
\eeq
Moreover, letting $\p_{n}^{k,max}$ be the set maximal set of tiles in $\p_n^{k}$, we set
\beq \label{pnbdi}
\p_{n,bd}^{k,i}:=\{P\in\A_n^{k}\,|\, \:\exists\:P_{kj}\in\p_{n}^{k,max}\:\textrm{s.t.}\:P\leq P_{kj}\:\:\&\:\: 20 I_{P}\cap (I_{P_{kj}})^c\not=\emptyset\:\}\:.
\eeq
Then, for each $\p_n^{k}$, we define its {\it boundary} forest component as
\beq \label{pnbdks}
\p_{n,bd}^{k}=\p_{n,bd}^{k,i}\cup\p_{n,bd}^{k,e}\:.
\eeq
The {\it normal} forest component is defined as
\beq \label{pnmks}
\p_{n,nm}^{k}:=\p_n^{k}\setminus \p_{n,bd}^{k}\:.
\eeq

Finally, we set
\beq \label{bdk}
\p_{bd}:=\bigcup_{k\in\N}\p_{n,bd}^{k}\,,
\eeq
and
\beq \label{nmk}
\p_{nm}:=\bigcup_{k\in\N}\p_{n,nm}^{k}\,.
\eeq

Now, here is our plan:
\begin{itemize}
\item for estimating the $L^2$-bound of the operator $T^{\p_{nm}}$ we will show that the family $\{T^{\p_{n,nm}^{k}}\}_k$ consists of almost orthogonal operators;

\item for treating the operator $T^{\p_{bd}}$ one simply notices that $\p_{bd}$ is a sparse forest and hence falls under the hypothesis of Proposition 1.
\end{itemize}

\begin{claim}\label{L2forest}
With the above notations, for $\eta=\eta(d)\in (0,1)$, one has
\beq\label{tpl2}
\left\|T^{\p_{nm}}\right\|_{2}\lesssim 2^{-\frac{n}{2}\,\eta}\:.
\eeq
\end{claim}

 In order to prove the above claim, using $TT^{*}-$method, it is enough to show that for some $c>0$ and $|k-k'|>10(1+c^{-1})$ we have
\beq\label{TT*n}
\|T^{{\p_{n,nm}^{k}}}\,{T^{\p_{n,nm}^{k'}}}^{*}\|_2 \lesssim e^{-c\:|k-k'|\:n}\:,
\eeq
\beq\label{T*Tn}
\|{T^{{\p_{n,nm}^{k}}}}^{*}\,T^{\p_{n,nm}^{k'}}\|_2 \lesssim e^{-c\:|k-k'|\:n}\:.
\eeq
Indeed, \eqref{tpl2} will then be easily derived, since for any $k\in\N$
\beq \label{fforest}
\|T^{{\p_{n,nm}^{k}}}f\|_2\lesssim 2^{-\frac{n}{2}\,\eta}\,\left\|f\right\|_{2}.
\eeq
Notice that \eqref{fforest} is a direct consequence of the Main Lemma and the second item in Observation \ref{redf}.

With this being said, let us start by proving \eqref{TT*n}.

Without loss of generality, we can suppose that $k'>k+10(1+c^{-1})$. Applying Cauchy-Schwarz we have
$$\left|\left\langle {T^{{\p_{n,nm}^{k}}}}^{*}f,\,{T^{\p_{n,nm}^{k'}}}^{*}g\right\rangle\right|\leq \|\chi_{A_n^{k'}}{T^{{\p_{n,nm}^{k}}}}^{*}f\|_2 \:\|{T^{{\p_{n,nm}^{k'}}}}^{*}g\|_2\:.$$
Here we have used that $\p_{n,nm}^{k'}$ is normal and thus $\textrm{supp } {T^{{\p_{n,nm}^{k'}}}}^{*}\subseteq A_n^{k'}$.

Next, from the way in which we have constructed $\p_{n,nm}^{k}$, we have that $$\forall\:P\in\p_{n,nm}^{k}\:\textrm{s.t.}\:I_{P^{*}}\cap A_n^{k+1}\not=\emptyset\:\Rightarrow\:I_{P^{*}}\nsubseteq A_n^{k+1}\:.$$

Thus, for any $P\in\p_{n,nm}^{k}$, we either have $I_{P^{*}}\cap A_n^{k+1}=\emptyset$ or the following relation holds:
\beq\label{densityofAIP}
\frac{|I_{P^{*}}\cap A_n^{k'}|}{|I_{P^{*}}|}\leq \frac{|I_{P^{*}}\cap A_n^{k'}|}{|I_{P^{*}}\cap A_n^{k+1}|}\lesssim e^{-c\:|k-k'|\:n}\:.
\eeq

Reaching this point, we remember that $\p_{n,nm}^{k}$ is an $L^{\infty}-$forest of $n^{th}$ generation and hence
\beq\label{rowdec}
\p_{n,nm}^{k}=\bigcup_{j=1}^{c\,2^{n}}\r_j^{k}\:,
\eeq
with each $\r_j^{k}$ a row.

Then, using \eqref{densityofAIP}, and applying Lemma \ref{ed} for $A:=A_n^{k'}$, we obtain
\beq \label{rowest}
\left\|\chi_{A_n^{k'}}\,{T^{\p_{n,nm}^{k}}}^{*}f\right\|_2\lesssim
\sum_{j=1}^{c\,2^{n}} \left\|\chi_{A_n^{k'}}\,{T^{\r_j^{k}}}^{*}f\right\|_2\lesssim e^{-c\:|k-k'|\:n}\,\left\|f\right\|_{2}\:,
\eeq
which proves \eqref{TT*n}.

We will now move on to the proof of \eqref{T*Tn}.

As before, we can start by first applying Cauchy-Schwartz
$$\left|\left\langle T^{{\p_{n,nm}^{k}}}f,\,T^{\p_{n,nm}^{k'}}g\right\rangle\right|\leq \|\chi_{A_n^{k'}}\,T^{{\p_{n,nm}^{k}}}f\|_2 \:\|T^{{\p_{n,nm}^{k'}}}g\|_2\:.$$

Based on \eqref{rowdec} and the fact that the operators $\{T^{\r_j^{k}}\}_j$ have disjoint supports, we have
$$\|\chi_{A_n^{k'}}\,T^{{\p_{n,nm}^{k}}}f\|_2 ^2=\sum_{j} \|\chi_{A_n^{k'}}\,T^{\r_j^{k}}f\|_2 ^2\lesssim 2^{n} \sup_{j}\|\chi_{A_n^{k'}}\,T^{\r_j^{k}}f\|_2 ^2\:.$$

Now, applying Lemma \ref{treecutp} to our row $\r_j^{k}$ (with the obvious adaptation of the partition $CZ(\I_{\p^*}^{min})$ there to our new context - call this new partition $\J_{\r_j^{k}}$) we have
\beq \label{massrow}
\|\chi_{A_n^{k'}}\,T^{\r_j^{k}}f\|_2 \lesssim 2^{-\frac{n}{2}}\,\left(\sup_{J\in\J_{\r_j^{k}}}\,\frac{|E_{A_n^{k'}}(J)|}{|J|}\right)^{\frac{1}{2}}\,\left\|f\right\|_{2}.
\eeq
Here, it is easy to remark that, from the construction of $\p_{n,nm}^{k}$, we have
$$\sup_{J\in\J_{\r_j^{k}}}\,\frac{|E_{A_n^{k'}}(J)|}{|J|}\lesssim e^{-c\:|k-k'|\:n}\,.$$
Thus, combining this last observation with \eqref{massrow}, we deduce
$$\|\chi_{A_n^{k'}}\,T^{{\p_{n,nm}^{k}}}f\|_2 \lesssim e^{-c\:|k-k'|\:n}\,\left\|f\right\|_{2}\:, $$
which together with \eqref{fforest} implies \eqref{T*Tn}.

\subsubsection{\bf The $L^p$ bound}\label{Lpb}
$\newline$

In this section, based on assumptions \eqref{unionp}-\eqref{unionp2} and Observation \ref{redf}, we will show that
\beq \label{lpnforestim}
\left\|T^{\p_{nm}}\right\|_{p}\lesssim_{p} 2^{-n\,\eta\,(1-\frac{1}{p^*})}\:.
\eeq

Our proof will be split in two cases:

$\newline$
\noindent \textbf{Case 1.} \textsf{Assume $1<p<2$.}
$\newline$

In this situation we notice that $p^*=p$ and thus \eqref{lpnforestim} is equivalent with
\beq \label{lpnforestimstar}
\left\|{T^{\p_{nm}}}^{*}\right\|_{p'}\lesssim_{p'} 2^{-\frac{n\eta}{p'}}\:.
\eeq

Firstly we notice - based on elementary interpolation techniques - that it is enough to prove \eqref{lpnforestimstar} only for $p'\in 2\N$ with $p'\geq 2$.

In this context, at the heuristic level, our goal is to show that:
\beq \label{Lpineq}
\|\sum_{k} {T^{{\p_{n,nm}^{k}}}}^{*}f\|_{p'}^{p'}\lesssim_{p'} \sum_{k}\|{T^{{\p_{n,nm}^{k}}}}^{*}f\|_{p'}^{p'}\,+\textrm{Error}\:,
\eeq
where the  "Error" term above is appropriately small and will be made precise in what follows.

Then, we notice that (up to conjugation), we have

$$\|\sum_{k} {T^{{\p_{n,nm}^{k}}}}^{*}f\|_{p'}^{p'}\approx_{p'}\sum_{{(k_1,\ldots , k_{p'}),\:(r_1,\ldots , r_{p'})\in \N^{p'}}\atop{r_1+\ldots +r_{p'}=p'}}
\int\large({T^{{\p_{n,nm}^{k_1}}}}^{*}f\large)^{r_1}\ldots\large({T^{{\p_{n,nm}^{ k_{p'}}}}}^{*}f\large)^{r_{p'}}$$
and after applying the H\"older and Jensen inequalities we further have
$$\|\sum_{k} {T^{{\p_{n,nm}^{k}}}}^{*}f\|_{p'}^{p'}\lesssim_{p'}$$
$$\sum_{{(k_1,\ldots , k_{p'}),\:(r_1,\ldots , r_{p'})\in \N^{p'}}\atop{r_1+\ldots+r_{p'}=p'}}\left(\int_{\bigcap_{j=1}^{p'} A_{n}^{k_j}}\,\left|{T^{{\p_{n,nm}^{k_1}}}}^{*}f\right|^{p'}\right)^{\frac{r_1}{p'}}\ldots
\left(\int_{\bigcap_{j=1}^{p'} A_{n}^{k_j}}\,\left|{T^{{\p_{n,nm}^{k_{p'}}}}}^{*}f\right|^{p'}\right)^{\frac{r_{p'}}{p'}}$$
$$\lesssim_{p'}\sum_{k}\sum_{m\in\N} |m+1|^{100p'}\,\int_{A_{n}^{k+m}}\,\left|{T^{{\p_{n,nm}^{k}}}}^{*}f\right|^{p'}\:.$$

Thus, we have just proved that for $p'\in 2\N$, with $p'>1$, we have that
\beq \label{L2pineq}
\|\sum_{k} {T^{{\p_{n,nm}^{k}}}}^{*}f\|_{p'}^{p'}\lesssim_{p'}\sum_{k}\sum_{m\in\N} |m+1|^{100p'}\,\int_{A_{n}^{k+m}}\,\left|{T^{{\p_{n,nm}^{k}}}}^{*}f\right|^{p'}\:.
\eeq
which trivially translates into
 \beq \label{Lppreineq}
\|\sum_{k} {T^{{\p_{n,nm}^{k}}}}^{*}f\|_{p'}^{p'}\lesssim_{p'}\sum_{k}\|{T^{{\p_{n,nm}^{k}}}}^{*}f\|_{p'}^{p'}\,+
\,\sum_{{m\geq10\,p'}\atop{m\in\N}} \sum_{k} m^{100p'}\,\int_{A_{n}^{k+m}}\,\left|{T^{{\p_{n,nm}^{k}}}}^{*}f\right|^{p'}\,.
\eeq
Notice that \eqref{Lppreineq} is the precise formulation of the heuristic described in relation \eqref{Lpineq}.

Next step will be to treat the main term
\beq \label{mtrm}
A=\sum_{k}\|{T^{{\p_{n,nm}^{k}}}}^{*}f\|_{p'}^{p'}\,.
\eeq

We first prove that it is enough to show that \eqref{lpnforestimstar} holds for $\p_{n,nm}^{k}$ (uniformly in $k$), that is
\beq \label{lpnforestimstarpar}
\left\|{T^{\p_{nm}^k}}^{*}f\right\|_{p'}\lesssim_{p'} 2^{-\frac{n\eta}{p'}}\,\|f\|_{p'}\:.
\eeq
Indeed, assume for the moment that \eqref{lpnforestimstarpar} holds.

Then, we first split the input of ${T^{{\p_{n,nm}^{k}}}}^{*}$ in disjoint sets $\{\chi_{A_n^{k+l}\setminus A_n^{k+l+1}}\}_{l\in\N}$ and notice that based on \eqref{lpnforestimstarpar}, for any $l\in\N$ one has
\beq \label{lpnforestimstarpar1}
\left\|{T^{\p_{nm}^k}}^{*} \chi_{A_n^{k+l}\setminus A_n^{k+l+1}}f\right\|_{p'}\lesssim_{p'} 2^{-\frac{n\eta}{p'}}\,\|\chi_{A_n^{k+l}\setminus A_n^{k+l+1}} f\|_{p'}\:.
\eeq
For $l\geq 2$ however, we can do better. To see this, we first apply standard H\"older inequality relative to the row decomposition of $\p_{n,nm}^{k}$:
\beq \label{kforestpps}
\|{T^{{\p_{n,nm}^{k}}}}^{*}f\|_{p'}\lesssim_{p'} (2^n)^{\frac{1}{p}}\,\left\{\sum_{j=1}^{c2^n}\|{T^{\r_{j}^{ k}}}^{*}f\|_{p'}^{p'}\right\}^{\frac{1}{p'}}.
\eeq

Now, from \eqref{densityofAIP} and Lemma \ref{treecutp}, we deduce\footnote{As in the $L^2-$case, if necessary, we can replace the original decompositions $\{\p_{n,nm}^{k}\}_{k\in\N}$, $\{A_n^{k}\}_{k\in\N}$ by the sparser correspondents $\{\p_{n,nm}^{k}\}_{k\in L\N}$, $\{A_n^{k}\}_{k\in L\N}$ where here $L\in\N$ with $L>10(1+c^{-1})$.}
\beq \label{prowpart}
\|{T^{\r_{j}^{k}}}^{*} (\chi_{A_n^{k+l}\setminus A_n^{k+l+1}}\,\cdot)\|_{p'}=\|\chi_{A_n^{k+l}\setminus A_n^{k+l+1}}\,T^{\r_{j}^{k}}\|_{p}
\lesssim_p \min (2^{-\frac{l\,n}{p}},\:2^{-\frac{n}{p}})\:.
\eeq

Denoting now with $E_{j}^k:=\bigcup_{P\in \r_{j}^{k}} E(P)$ and using the fact that $\{E_{j}^k\}_{j}$ are pairwise disjoint, we have from \eqref{kforestpps} and \eqref{prowpart} that
\beq \label{kforestpps1}
\|{T^{{\p_{n,nm}^{k}}}}^{*} \chi_{A_n^{k+l}\setminus A_n^{k+l+1}}f\|_{p'}\lesssim_{p'} \min\{1, 2^{-\frac{n(l-1)}{p}}\}\,\|\chi_{A_n^{k+l}\setminus A_n^{k+l+1}}f\|_{p'}\,.
\eeq

Deduce from \eqref{lpnforestimstarpar1}, \eqref{kforestpps1} and H\"older's inequality that
\beq \label{revLpineqsplit1}
\|{T^{{\p_{n,nm}^{k}}}}^{*}f\|_{p'}^{p'}\lesssim_{p'}
\sum_{l\in\N} (l+1)^{p'}\,2^{-n\eta}\,\min\{1, 2^{-\frac{n(l-1-\eta)p'}{p}}\}\,\|\chi_{A_n^{k+l}\setminus A_n^{k+l+1}}f\|_{p'}^{p'}\:.
\eeq
Replacing now \eqref{revLpineqsplit1} in \eqref{mtrm} and summing over $k$ we conclude that
\beq \label{Aterm}
A\lesssim_{p'} 2^{-n \eta}\,\|f\|_{p'}^{p'}\:.
\eeq

Returning now to the proof of \eqref{lpnforestimstarpar} the simplest approach is provided by the following short argument - which holds uniformly in $k$:
\begin{itemize}
\item for the case $p=p'=2$ we already know that \eqref{lpnforestimstarpar} holds from \eqref{fforest} (or equivalently from the Main Lemma);

\item for $p'\in 2\N\setminus\{0\}$, applying similar steps with those in \eqref{kforestpps}-\eqref{revLpineqsplit1} but with no extra-assumption \eqref{lpnforestimstarpar} we notice the trivial bound:
\beq \label{revLpineqsplit11}
\|{T^{{\p_{n,nm}^{k}}}}^{*}f\|_{p'}\lesssim_{p'}\|f\|_{p'}\:.
\eeq
\end{itemize}
Conclude from the above using standard interpolation that \eqref{lpnforestimstarpar} holds.

A different, more involved, but direct approach (\textit{i.e.} not appealing to formal interpolation), was present in an earlier version of this paper and was based on the following \emph{heuristic} hinted by an \textit{informal} interpolation argument:

Recalling \eqref{almostortho} in Main Lemma, we know - ignoring the error term - that
\beq \label{kforest1}
\|{T^{{\p_{n,nm}^{k}}}}^{*}f\|_{2}\lesssim \left(\sum_{j=1}^{c2^n}\|{T^{\r_{j}^{k}}}^{*}f\|_{2}^2\right)^{\frac{1}{2}}\,.
\eeq
We also trivially have
$$\|{T^{{\p_{n,nm}^{k}}}}^{*}f\|_{\infty}\lesssim \sum_{j=1}^{c2^n}\|{T^{\r_{j}^{k}}}^{*}f\|_{\infty}\,.$$

Thus, at the heuristic level, we expect for any $2\leq p'<\infty$ to have
\beq \label{kforestpp}
\|{T^{{\p_{n,nm}^{k}}}}^{*}f\|_{p'}\lesssim_{p'} \left\{\sum_{j=1}^{c2^n}\|{T^{\r_{j}^{ k}}}^{*}f\|_{p'}^{p}\right\}^{\frac{1}{p}}
\lesssim (2^n)^{\frac{1}{p}-\frac{1}{p'}}\,\left\{\sum_{j=1}^{c2^n}\|{T^{\r_{j}^{ k}}}^{*}f\|_{p'}^{p'}\right\}^{\frac{1}{p'}}.
\eeq

The key message here is that one is able to \emph{decouple} the information carried by the rows of a forest with a gain of $(2^n)^{-\frac{1}{p'}}$ over the trivial H\"older bound. Notice at this point that \textit{any} gain over the trivial bound $(2^n)^{\frac{1}{p}}$ would be enough for our claim \eqref{lpnforestimstarpar}.

The precise form of this decoupling argument is given by the following:

\begin{o0}\label{Almoortogpcas}
\textit{Let $p'\in 2\N$, $p'\geq 2$ and $\p_{n,nm}^{k}$ be an $L^{\infty}$-forest of generation $n$ whose standard decomposition into rows is given by $\{\r_{j}^{ k}\}_{j=1}^{c 2^n}$. Assume that any two distinct trees within this row decomposition are normal and $2^{100\,n\,d\,p'}$-separated. Then, there exists $\eta\in (0,1)$ such that the following holds:
\beq \label{pdecoupling}
\|{T^{{\p_{n,nm}^{k}}}}^{*}f\|_{p'}\lesssim_{p'} (2^n)^{\frac{1}{p}-\frac{\eta}{p'}}\,\left\{\sum_{j=1}^{c2^n}\|{T^{\r_{j}^{ k}}}^{*}f\|_{p'}^{p'}\right\}^{\frac{1}{p'}}\,+\,2^{-\frac{10 n}{p'}}\,\|f\|_{p'}\,.
\eeq
The proof of this statement relies in a key fashion on the separateness assumption of the trees further reflected into the time-frequency localization properties of each of the maximal trees belonging to the forest. The beauty of this approach is that it provides the desired $L^{p'}-$decay (for $p'\in 2\N$, $p'\geq 2$) in a \emph{direct} fashion with \emph{no actual usage} of interpolation methods. Due to space limitation and much easier alternative provided above we choose not to present here a proof of this statement.}
\end{o0}

We pass now to the error term
$$B:=\sum_{m\geq 10\,p'} \sum_{k} m^{100p'}\,\int_{A_{n}^{k+m}}\,\left|{T^{{\p_{n,nm}^{k}}}}^{*}f\right|^{p'}\,.$$

We first notice that
\beq \label{holdlprow}
\int_{A_{n}^{k+m}}\,\left|{T^{{\p_{n,nm}^{k}}}}^{*}f\right|^{p'}\lesssim
 (2^n)^{p'}\,\sum_{j}\int_{A_{n}^{k+m}}\,\left|{T^{\r_{j}^{k}}}^{*}f\right|^{p'}\:.
\eeq

Now, based on \eqref{densityofAIP} and Lemma \ref{lptree}, we deduce that for each $j$ we have
\beq \label{rowakm}
\int_{A_{n}^{k+m}}\,\left|{T^{\r_{j}^{k}}}^{*}f\right|^{p'}\lesssim 2^{-m\,n}\,\|f\|_{p'}^{p'}\:.
\eeq

Combining \eqref{prowpart} with \eqref{rowakm} we further have
\beq \label{rowakm1}
\int_{A_{n}^{k+m}}\,\left|{T^{\r_{j}^{k}}}^{*}\chi_{A_n^{k+l}\setminus A_n^{k+l+1}}\,f\right|^{p'}
\lesssim_{p'} 2^{-\frac{p'(l+1)n}{4p}}\,2^{-\frac{m\,n}{2}}\,\|\chi_{A_n^{k+l}\setminus A_n^{k+l+1}}\,f\|_{p'}^{p'}\:.
\eeq

Next, proceeding in a similar fashion with \eqref{revLpineqsplit1}, we have

\beq \label{rowjcut}
\int_{A_{n}^{k+m}}\,\left|{T^{\r_{j}^{k}}}^{*}\,f\right|^{p'}
\lesssim_{p'} 2^{-\frac{m\,n}{2}}\,\sum_{l\in\N} (l+1)^{p'}\,2^{-\frac{l\,n\,p'}{4p}}\, \|\chi_{A_n^{k+l}\setminus A_n^{k+l+1}}\,\chi_{E_{j}^k}f\|_{p'}^{p'}\:.
\eeq

Putting together \eqref{holdlprow} and \eqref{rowjcut} we deduce that
$$B\lesssim \sum_{m\geq 10\,p'} m^{100p'}\,2^{-\frac{m\,n}{2}}\,(n\,2^{n})^{p'}\sum_k\sum_{l\in\N} (l+1)^{p'}\,2^{-\frac{l\,n\,p'}{4p}}\, \|\chi_{A_n^{k+l}\setminus A_n^{k+l+1}}\,f\|_{p'}^{p'}\:,$$
and hence
\beq \label{Bestim}
B\lesssim_{p'} 2^{-n}\,\|f\|_{p'}^{p'}\:.
\eeq

Finally, from \eqref{Aterm} and  \eqref{Bestim}, we conclude that \eqref{lpnforestim} holds.

$\newline$
\noindent \textbf{Case 2.} \textsf{Assume $2<p<\infty$.}
$\newline$

In this situation we have that $p^*=p'$ and hence \eqref{lpnforestim} is equivalent with
\beq \label{lpnforestimstar2}
\left\|T^{\p_{nm}}\right\|_{p}\lesssim_{p} 2^{-\frac{n}{p}}\:.
\eeq

Once at this point, we notice that that we can follow line by line the same arguments as in Case 1 by just dropping the adjoint symbol in the corresponding proof. The key aspect that allows us to work with only this simple modification is that all the $L^{\infty}-$forests appearing in the reasonings from Case 1 consist of \emph{normal} trees and hence $\textrm{supp}\,T^{{\p_{n,nm}^{k}}}\subseteq A_n^k$.

\begin{flushright}
$\Box$
\end{flushright}

\section{\bf Remarks}\label{remark}

In this section we will discuss some applications and consequences of the discretization procedure presented in Section 5.1.

1) The first remark is a consequence of a fruitful conversation that I had with C. Thiele and M. Bateman,
and refers to a vector-valued variant of the Carleson Theorem. More precisely, using the above discretization procedure
(and thus eliminating the exceptional sets), we devised an alternative proof that for any
$1<p,q<\infty$ one has
\footnote{Here we use the notations from Section 1.}
\beq \label{vectorval}
\left\|(\sum_{k} |C f_k|^{q})^{\frac{1}{q}}\right\|_p\lesssim_{p,q} \left\|(\sum_{k} |f_k|^{q})^{\frac{1}{q}}\right\|_p\:,
\eeq
an inequality that had been proven in \cite{GMF} using weighted
and extrapolation theory. Nevertheless,
as a consequence of the Theorem presented in this paper, one has that \eqref{vectorval}
holds with $C=C_{1,1}$ replaced by $C_{d,1}$.
\medskip

2) As mentioned in the introduction, our discretization procedure was designed for obtaining the following informal principle:
\medskip

If $\p=\bigcup_{k}\p_k \subseteq \P_n$ is a collection of separated trees, and  $\cN_{\p}$ stands for the usual counting function associated with $\p$, then
\beq \label{infp}
\|\sum_{k} {T^{\p_k}}^{*}\,f\|_2\lesssim \log (10+ \|\cN_{\p}\|_{BMO_{C}})\,\left(\sum_k \|{T^{\p_k}}^{*}\,f\|_2^2\right)^{\frac{1}{2}}\:.
\eeq

In both \cite{dtt} and \cite{f}, this principle was only present in a weaker form with $\|\cN_{\p}\|_{BMO_{C}}$ replaced by $\|\cN_{\p}\|_{L^{\infty}}$ in \eqref{infp}, thus causing some intricate technicalities in order to treat the so-called ``exceptional sets" on which the quantity $\|\cN_{\p}\|_{L^{\infty}}$ is too large. Through our new approach having as a consequence the rigorous proof of \eqref{infp} we are now able to completely discard the analysis of the exceptional sets,
answering thus to a question raised by C. Fefferman (see the Remarks in \cite{f}).
\medskip

3) In an earlier version of the paper we presented a slightly different partition of the collection of tiles $\P$, that has its own merits and that we choose to present very briefly below as an alternative:

Set $\P_{0}=\emptyset$ and suppose that we have defined (for some $n\geq 1$) the sets $\{\P_k\}_{k<n}$. We describe now the algorithm for constructing the set $\P_n$.

First step consists of selecting the family $\p_{n}^{0,max}$ of the maximal tiles $P\in\P\setminus\bigcup_{k<n}\P_k$ with $\frac{|E(P)|}{|I_P|}\geq 2^{-n}$. After that, we collect the time-intervals of these maximal tiles into the set $\I_{n}^{0}$ and form with them the counting function
$$\cN_{n}^{0}:=\sum_{I\in\I_{n}^{0}} \chi_{I}\:.$$

Next, by using John-Nirenberg inequality we remark that the set

$$A_n^{1}:=\{x\in [0,1]\,|\, \sum_{I\in\I_{n}^{0}}\chi_{I}(x)> c\,n\,\|\cN_n^{0}\|_{BMO_{C}}\}\:,$$
has the measure  $|A_n^{1}|\leq e^{-100\,n }$.

Further, we construct $\p_{n}^{1,max}$ to be the collection of maximal tiles $P\in\P\setminus\bigcup_{k<n}\P_k$ with $\frac{|E(P)|}{|I_P|}\geq 2^{-n}$ and $I_P\subseteq A_n^{1}$. Also, as before, define  $\I_{n}^{1}:=\{I\,|\,P=[\vec{\a},I]\in\p_{n}^{1,max}\}$, the counting function $\cN_{n}^{1}:=\sum_{I\in\I_{n}^{1}}\chi_{I}$ and the exceptional set
$$A_n^{2}:=\{x\in [0,1]\,|\, \sum_{I\in\I_{n}^{1}}\chi_{I}(x)> c\,n\,\|\cN_n^{1}\|_{BMO_{C}}\}\:.$$

Proceeding by induction, at the end of the day, we will have constructed the collection of sets of maximal tiles $\{\p_{n}^{k,max}\}_k$, the collection of sets representing the time-intervals - $\{\I_{n}^{k}\}_k$, the collection of counting functions $\{\cN_{n}^{k}\}_k$ and finally the level sets $\{A_n^{k}\}_k$.

Reaching this point, we state the following important consequences of our construction:
\beq\label{excset}
|A_n^{k}|\leq e^{- 100\,|k-l|\:n}\,|A_n^{l}|\:,
\eeq
\beq\label{countingset}
\sup_{k}\|\cN_n^{k}\|_{BMO_{C}}\leq 2^{n}\:\:\textrm{and}\:\:\sup_{k}\|\cN_n^{k}\|_{L^{\infty}(A_n^{k}\setminus A_n^{k+1})}\lesssim n\,2^{n}\:.
\eeq

Moreover, if we set the counting function $$\cN_{n}:=\sum_{k} \cN_{n}^{k} \:,$$ we also have
\beq\label{globalcountingset}
\|\cN_n\|_{BMO_{C}}\lesssim 2^{n}\:.
\eeq

With this done, for each $k\in\N$, we define
\beq\label{pk}
\p_n^k:=\left\{P=[\vec{\a},I]\,|\,\begin{array}{rl} I\subseteq A_n^{k},\:I\nsubseteq A_n^{k+1}\:\textrm{and}\:\:\:\:\:\:\:\\ A_{\P\setminus\bigcup_{j<n}\P_j,A_{n}^{k}}(P)\in [2^{-n}, 2^{-n+1})\end{array}\right\}\:,
\eeq
and set
\beq\label{Pn}
\P_n:=\bigcup_{k\geq0}\p_n^k\:.
\eeq

Finally, remark that we have
\beq\label{P}
\P=\bigcup_{n\geq0}\P_n\:.
\eeq

Compared to the current partitioning presented in Section \ref{parttil}, the approach above has the advantage of being shorter and obeying a simpler algorithm. However it has the drawback that it does not preserve the \emph{convexity} of the family $\P_n$. In order to deal with this, one needs to work with so called ``generalized" trees: a family of tiles $\p$ is called a ``generalized" tree with top $P_0=[\vec{\a}_0, I_0]$ iff it obeys 1) and 2) in Definition \ref{tree} and $\p=\bigcup_{k}\p_k$ with each $\p_k$ tree with top $P_k=[\vec{\a}_k, I_k]$ such that the following Carleson packing condition is satisfied: $\sum_{I_k\subseteq I_j} |I_k|\lesssim |I_j|$ for any $j$. At the heuristic level, the distinction between these two ways of partitioning our family of tiles resides in the \emph{moment} in which one decides to perform a last stopping time argument. In order to avoid possible complications within the main proof we have chosen to do all the stopping time reasonings when defining the family $\P_n$ so that later we are able to avoid the discussion about ``generalized" trees and only work with the ``classical" tree structures.
\medskip

4) This remark is dedicated to an important feature regarding the behavior of the counting functions $\{\cN_{n}\}_n$ as defined in \eqref{gcntfn}. In \cite{lv9}, the author characterized the $L^{1}-$weak behavior of the so called lacunary Carleson operator (see the next remark). A key idea in that study was the understanding of the newly defined concept of \emph{grand maximal function}, which in our current context is defined as follows:

- fix $j\in\N$ and set
\beq\label{countj}
\cN(j):=\frac{1}{2^{j-1}}\sum_{n=2^{j-1}+1}^{2^j} \frac{1}{2^{n-1}} \cN_n\;,
\eeq
and define the grand maximal counting function of order $l\in\N$, ($l\geq 2$) by
\beq\label{grandmaxcount}
\cN^{[l]}:=\sup_{j\leq l}\cN(j)\,.
\eeq

With these we have the following key property:
\beq\label{grandmaxcountkey}
\|\cN^{[l]}\|_{1,\infty}\lesssim \log l\,,
\eeq
with the right bound being sharp.

It is precisely this last fact - that is, the possibility of existence of extremal configurations of tiles that realize the reverse inequality
\beq\label{grandmaxcountkeyy}
\|\cN^{[l]}\|_{1,\infty}\gtrsim \log l
\eeq
- that makes our tile partitioning process more delicate and suggests in \eqref{A1n1} the level set cut at height $c\,n$ in order to have a good control over \eqref{PPn1}.
\medskip

5) Finally, the previous remark connects with the celebrated theme regarding the behavior of the Carleson operator $C$ near $L^1$.

At the foundation of this theme resides the following heuristic question\footnote{In what follows we embrace the formalism from our paper \cite{lv9}}:

\textit{What is the behavior of the (almost everywhere) pointwise convergence of the
Fourier Series between the two known cases for the Lebesgue-scale spaces $L^p(\TT)$:
\begin{itemize}
\item $p=1$, divergence of the Fourier Series (Kolmogorov)
\item $p>1$, convergence of the Fourier Series (Carleson-Hunt) ?
\end{itemize}}

Using now the fact that the pointwise convergence of the Fourier Series is directly related to the $L^{1,\infty}$-behavior of the Carleson operator one can reformulate the above vague question into a precise problem
\medskip

\noindent\textbf{Problem:} 1) \textit{Let $Y\subseteq L^1(\TT)$ be a r.i. (quasi-)Banach space. Provide necessary and sufficient conditions for $Y$ in order to be a $C-$space, that is, $\exists\:c>0$ such that
\beq\label{prl1}
\|Cf\|_{1,\infty}\leq c\,\|f\|_{Y}\:.
\eeq}

2) \textit{In Lorentz space terminology, the above can be expressed as: Give a satisfactory description of the Lorentz spaces $Y\subseteq L^1(\TT)$ that are also $C-$spaces. If such exists, describe the maximal Lorentz $\mathcal{C}-$space $Y_0$.}

\medskip
In terms of known results we have two possible directions:
\begin{itemize}
\item on the \textit{negative} side (i.e. aiming for decreasingly smaller Banach rearrangement invariant spaces that are not $\mathcal{C}-$spaces): as mentioned above, the history of this direction starts with the result of Kolmogoroff, showing that $L^1(\TT)$ is not a $C-$space. The next results are due to Chen (\cite{ch}), Prohorenko (\cite{P}) and K\"orner (\cite{Kor}). The best up to date result belongs to Konyagin (\cite{koDivf}, \cite{koDivff}) who proved that for $\phi(u)=o(u\sqrt{\frac{\log u}{\log\log u}})$ as $u\rightarrow \infty$ the space $X=\phi(L)$ does not admit pointwise convergence.

\item on the \textit{positive} side (i.e. identifying increasingly larger $C$-spaces $Y$): historically, the topic starts with the results of Carleson and Hunt for $Y=L^p(\TT),\:p>1$. Next Sj\"olin (\cite{sj3}) showed that one can take
 $Y=L\log L \log\log L$  while F. Soria (\cite{So1}, \cite{So2}) increased $Y$ to a r.i. quasi-Banach space denoted $B^{*}_{\v}$. The best current results belong to Antonov (\cite{An}) for the Lorentz-space $Y=L\log L\log \log\log L$ and to Arias de Reyna (\cite{Ar}) for the quasi-Banch space $Y=QA$ (a  r.i. quasi-Banach space that essentially has as its largest possible Lorentz space precisely Antonov's space)
\end{itemize}

In this context, there are several things that are worth being mentioned:
\begin{itemize}
\item all the positive results rely in their proof on extrapolation methods;

\item in \cite{lv7} the author reproved all the above positive results via a unitary method relying entirely only on time-frequency tools.

\item currently, all the positive results can be explained entirely based on the behavior of the grand maximal counting function (of order $l$) \eqref{grandmaxcount}, more precisely on the fact that there are (extremizers) configurations of tiles for which inequality \eqref{grandmaxcountkeyy} holds. In order to explain our claim, at least at the heuristic level, we use duality and write our Carleson operator as a bilinear form given by
\beq\label{DF}
\Lambda(f,g):=<Cf,\,g>\approx\sum_{n}\sum_{{\p_k\subset\P_n}\atop{\p_k\:L^{\infty}-\textrm{forest}}}
\sum_{{\p\subseteq\p_k}\atop{\p\:\textrm{maximal tree}}} <C^{\p}f,\,g>\:.
\eeq
 The key issue is that currently there are no methods \emph{near $L^1$} to distinguish\footnote{There is a similar problem regarding the maximal boundedness range for the Bilinear Hilbert transform - see Section \ref{BHT}.} between the absolute and the conditional summation in the RHS of \eqref{DF}. One can argue that if one assumes absolute summation in the RHS of \eqref{DF}, then Antonov's result is the best possible result and, further on, it is a direct consequence of the logarithmic divergence of the $L^{1,\infty}-$norm of the grand maximal counting function (of order $l$).

\item there exists an old model problem for the problem stated above which has it's own history. This model problem regards the almost everywhere convergence of lacunary sequences of partial Fourier sums and goes back to early 20th century in works of Kolmogorov, (see \cite{Kat}),  Littlewood and Paley, (\cite{LP}), and  Zygmund (\cite{Zyg}). In the quest for identifying the largest possible Lorentz space (or r.i. quasi Banach space) for which one has almost everywhere convergence along lacunary sequences of partial Fourier sums some partial progress has been made - see the works of Chen (\cite{ch}), Prohorenko (\cite{P}), K\"orner (\cite{Kor}) and later Konyagin (\cite{koDivf}, \cite{koDivff}).

    More recently, (\cite{lv9}), the author succeeded in giving a definitive answer to this problem. Indeed, by defining the lacunary Carleson operator associated with an (arbitrary) lacunary sequence $\{n_j\}_j$ by
\beq\label{carlac}
C_{lac}^{\{n_j\}_j}f(x):=\sup_{j\in\N}\left|\int_{\TT}e^{i\,2\pi\,n_j\,(x-y)}\,\cot(\pi\,(x-y))\,f(y)\,dy\right|\:,
\eeq
  we showed that  $\:\exists\:C_1=C_1(\{n_j\}_j)>0$ such that
\beq\label{carlblac}
\|C_{lac}^{\{n_j\}_j}f\|_{1,\infty}\leq C_1\,\|f\|_{L\log\log L\log\log\log\log L}\:\:,
\eeq
  and moreover that this result is essentially sharp. The proof relies in a key fashion on the properties of the grand maximal counting function explained above.

  Also, very recently, (\cite{lvLac}), we provided the sharp result regarding the strong $L^1$ bound for the lacunary Carleson operator, that is, $\:\exists\:C_2=C_2(\{n_j\}_j)>0$ such that
\beq\label{carlblac1}
\|C_{lac}^{\{n_j\}_j}f\|_{1}\leq C_2\,\|f\|_{L\log L}\:\:.
\eeq
\end{itemize}

Returning now to the original problem of the pointwise convergence of the full sequence of partial sums, we mention that the recent works \cite{lv9} and \cite{lvLac}, unraveled several subtle key points:
\begin{itemize}
\item the \textit{structure of the frequencies} of the trees involved in the time-frequency decomposition of the Carleson operator plays a fundamental role in identifying larger classes of r.i. Banach spaces for which we have pointwise convergence.
\item the \textit{structure of the input function} creates certain ``resonances" with the structure of the frequencies from the above item.
\end{itemize}

As a consequence, we expect that structural theorems from additive combinatorics will play a fundamental role in any relevant advancement on the problem.

We end by listing the three relevant main conjectures in this subject:

\medskip
\noindent\textbf{Conjecture 1.} [\textsf{$L^{1,\infty}$-behavior}]  \textit{The largest Lorentz space $Y_0\subseteq L^1(\TT)$ such that  $\exists\:c=c_{Y_0}>0$ with
\beq\label{prlc1}
\|Cf\|_{1,\infty}\leq c\,\|f\|_{Y_0}\:,
\eeq
is $Y_0=L\sqrt{\log L}$.}
\medskip

\medskip
\noindent\textbf{Conjecture 2.} [\textsf{$L^{1}$-behavior}]  \textit{The largest Lorentz space $Y_1\subseteq L^1(\TT)$ such that  $\exists\:c=c_{Y_1}>0$ with
\beq\label{prlc1}
\|Cf\|_{1}\leq c\,\|f\|_{Y_1}\:,
\eeq
 is $Y_1 = L\log L$.}
\medskip

If true, Conjecture 1 above is essentially sharp due to the result of Konyagin (\cite{koDivf}, \cite{koDivff}) while if Conjecture 2 is true then this is definitely sharp due to the fact that both the Hardy-Littlewood maximal operator and the Hilbert transform map sharply $L \log L$ into $L^1$.

Finally, we present a last conjecture, that while weaker than both conjectures above, its resolution would still be a major breakthrough in the field of time-frequency analysis due to the new methods that one needs to develop:

\medskip
\noindent\textbf{Conjecture 3.} [\textsf{$L^{1,\infty}$-intermediate behavior}]  \textit{Prove that there exits a constant $c>0$ such that the following holds
\beq\label{prlc3}
\|Cf\|_{1,\infty}\leq c\,\|f\|_{L \log L}\:.
\eeq}
\medskip

\section{\bf Appendix - Results on the $L^{\infty}-$distribution of polynomials}\label{append}
$\newline$

\noindent{\bf Lemma A.} {\it If $q\in\Q_{d-1}$ and $I,\:J$ are some (not
necessarily dyadic) intervals obeying $I\supseteq J$, then there
exists a constant $c(d)\leq (2d)^d$ such that
$$ \left\|q\right\|_{L^{\infty}(I)}\leq
c(d)\left(\frac{|I|}{|J|}\right)^{d-1}\left\|q\right\|_{L^{\infty}(J)}\:.$$}
$\newline${\it Proof.} Let $\{x_J^k\}_{k\in\{1,\ldots , d\}}$ be
obtained as in the
 procedure described in Section 2. Then, since $q\in\Q_{d-1}$, for any $x\in I$ we have that
 \beq \label{Lag}
 q(x):=\sum_{j=1}^{d}\frac{\prod_{k=1\atop{k\not=j}}^{d}(x-x_J^k)}
 {\prod_{k=1\atop{k\not=j}}^{d}(x_J^j-x_J^k)}\:q(x_J^j)\:.
 \eeq
As a consequence,
$$\left\|q\right\|_{L^{\infty}(I)}\leq
d\:\left\|q\right\|_{L^{\infty}(J)}\sup_{j\atop{x\in I}}\left|
\frac{\prod_{k=1\atop{k\not=j}}^{d}(x-x_J^k)}{\prod_{k=1\atop{k\not=j}}^{d}(x_J^j-x_J^k)}\right|\leq
d\:\left\|q\right\|_{L^{\infty}(J)}\frac{|I|^{d-1}}{(\frac{|J|}{2d})^{d-1}}\:.$$

\begin{flushright}
$\Box$
\end{flushright}

\noindent{\bf Lemma B.} {\it If $q\in\Q_{d-1},\:\eta>0$ and $I\subset\TT$
some (dyadic) interval, then \beq\label{levels} |\{y\in
I\:|\:|q(y)|<\eta\}|\leq c(d)
\left(\frac{\eta}{\left\|q\right\|_{L^{\infty}(I)}}\right)^{\frac{1}{d-1}}|I|\:.
\eeq } $\newline${\it Proof.} The set $A_{\eta}=\{y\in
I\:|\:|q(y)|<\eta\}$ is the pre-image of $(-\eta,\eta)$ under a polynomial
of degree $d-1$, so it can be written as
$$A_{\eta}=\bigcup_{k=1}^r J_k(\eta)\:,$$
where $r\in\N,\:r\leq d-1$ and $\{J_k(\eta)\}_k$ are open intervals.
Now all that remains is to apply the previous lemma with
$J=J_k(\eta)$ for each $k$.
\begin{flushright}
$\Box$
\end{flushright}

\noindent{\bf Lemma C.} {\it If $P=[\a^1,\a^2,\ldots , \a^d,I]\in\P$ and $q\in
P$, then
$$\left\|q-q_P\right\|_{L^{\infty}(\tilde{I})}\leq
c(d)\:|I|^{-1}\:.$$} $\newline${\it Proof.}
Set $u:=q-q_P$; then, since both $q,\:q_P\in P$, we deduce (for all
$k\in\{1,\ldots , d\}$):
$$u(x_I^k)\in [-|I|^{-1},|I|^{-1}]\:.$$
On the other hand,
$$u(x):=\sum_{j=1}^{d}\frac{\prod_{k=1\atop{k\not=j}}^{d}(x-x_I^k)}{\prod_{k=1\atop{k\not=j}}^{d}(x_I^j-x_I^k)}\:u(x_I^j)\:\:\:\:\:\:\:\:\:\forall\:\:x\in I\:.$$
Then, proceeding as in Lemma A, we conclude
$$\left\|u\right\|_{L^{\infty}(I)}\leq d\:
|I|^{-1}\frac{|I|^{d-1}}{(\frac{|I|}{2\,d})^{d-1}}\leq (2\,d)^d\: |I|^{-1}\:.$$

\begin{flushright}
$\Box$
\end{flushright}

\end{document}